\documentclass[a4paper,11pt,leqno]{amsart}
\usepackage{amsmath,amssymb,cite,mathrsfs}

\newtheorem{theorem}{Theorem}[section]
\newtheorem{lemma}[theorem]{Lemma}
\newtheorem{proposition}[theorem]{Proposition}

\theoremstyle{definition}

\newtheorem{example}[theorem]{Example}

\theoremstyle{remark}
\newtheorem{remark}[theorem]{Remark}

\numberwithin{equation}{section}

\DeclareMathOperator{\Func}{\mathfrak{F}}
\DeclareMathOperator{\Vol}{\mathrm{Vol}}
\DeclareMathOperator{\Vrt}{\mathrm{Vert}}
\DeclareMathOperator{\Conv}{\mathrm{Conv}}
\def\bd#1{\mathbf{#1}}
\def\fs{\Psi}
\def\id{\mathrm{id}}
\def\prs{\overline{\Delta}}
\def\Aut{\mathrm{Aut}(\Delta)}
\def\GL{\mathrm{GL}}

\def\abs#1{\lvert#1\rvert}
\def\norm#1{\lVert#1\rVert}
\long\def\comment#1{}

\begin{document}
\title[Witten multiple zeta-functions]
   {On Witten multiple zeta-functions associated with semisimple Lie
algebras III}
\author{Yasushi Komori, Kohji Matsumoto \and Hirofumi Tsumura}
\thanks{\noindent{Y.\,Komori:} Graduate School of Mathematics, Nagoya University, Chikusa-ku, Nagoya 464-8602 Japan, email:\,\texttt{komori@math.nagoya-u.ac.jp}\\
\ \ {K.\,Matsumoto:} Graduate School of Mathematics, Nagoya University, Chikusa-ku, Nagoya 464-8602 Japan, email:\,\texttt{kohjimat@math.nagoya-u.ac.jp}\\
\ \ \ {H.\,Tsumura:} Department of Mathematics and Information Sciences, Tokyo Metropolitan University, 1-1, Minami-Ohsawa, Hachioji, Tokyo 192-0397 Japan, email:\,\texttt{tsumura@tmu.ac.jp}}

\begin{abstract}
  We prove certain general forms of functional relations among Witten
  multiple zeta-functions in several variables (or zeta-functions of
  root systems).  The structural background of those functional
  relations is given by the symmetry with respect to Weyl groups.
  From those relations we can deduce explicit expressions of values of
  Witten zeta-functions at positive even integers, which is written in
  terms of generalized Bernoulli numbers of root systems.  Furthermore
  we introduce generating functions of those Bernoulli numbers of root
  systems, by which we can give an algorithm of calculating Bernoulli
  numbers of root systems.
\end{abstract}

\maketitle

\baselineskip 13pt

\section{Introduction} \label{sec-1}

Let ${\frak g}$ be a complex semisimple Lie algebra with rank $r$.
The Witten zeta-function associated with ${\frak g}$ is defined by
\begin{align}
  \label{1-1}
    \zeta_W(s;{\frak g})=\sum_{\varphi}(\dim\varphi)^{-s},
\end{align}
where the summation runs over all finite dimensional irreducible
representations $\varphi$ of ${\frak g}$.

Let ${\Bbb N}$ be the set of positive integers, ${\Bbb N}_0={\Bbb
  N}\cup\{0\}$, ${\Bbb Z}$ the ring of rational integers, ${\Bbb Q}$
the rational number field, ${\Bbb R}$ the real number field, ${\Bbb
  C}$ the complex number field, respectively.  Witten's motivation
\cite{Wi} of introducing the above zeta-function is to express the
volumes of certain moduli spaces in terms of special values of
\eqref{1-1}.  This expression is called Witten's volume formula, which
especially implies that
\begin{align}
  \label{1-2}
    \zeta_W(2k;{\frak g})=C_{W}(2k,{\frak g})\pi^{2kn}
\end{align}
for any $k\in{\Bbb N}$, where 
$n$ is the number of all positive roots and
$C_{W}(2k,{\frak g})\in{\Bbb Q}$ (Witten \cite{Wi}, Zagier \cite{Za}). 

When $\mathfrak{g}=\mathfrak{sl}(2)$, the corresponding Witten zeta-function is nothing
but the Riemann zeta-function $\zeta(s)$.   It is classically known that
\begin{equation}
\label{eq:r_zeta}
  \zeta(2k)=\frac{-1}{2}\frac{(2\pi\sqrt{-1})^{2k}}{(2k)!}B_{2k}\qquad (k\in\mathbb{N}),
\end{equation}
where $B_{2k}$ is the $2k$-th Bernoulli number.   Formula \eqref{1-2} is a
generalization of \eqref{eq:r_zeta}, but the values $C_W(2k,\mathfrak{g})$ 
are not explicitly
determined in the work of Witten and of Zagier.
    In the present paper we introduce a root-system theoretic
generalization of Bernoulli numbers and periodic Bernoulli functions,
and express $C_W(2k,\mathfrak{g})$ explicitly in terms of those generalized
Bernoulli periodic functions $P(\mathbf{k},\mathbf{y};\Delta)$ (defined in Section \ref{sec-4}).
This result will be given in Theorem \ref{thm:W-Z}.    Note that
Szenes \cite{Sz98,Sz03} also studied generalizations of Bernoulli
polynomials from the viewpoint of the theory of arrangement of
hyperplanes, which include $P(\mathbf{k},\mathbf{y};\Delta)$ mentioned above.
However, our root-system theoretic approach enables us to show
that our $P(\mathbf{k},\mathbf{y};\Delta)$ and its generating function 
$F(\mathbf{k},\mathbf{y};\Delta)$ are
quite natural extensions of the classical ones
(see Theorems proved in Section \ref{sec-6}).

Our explicit expression of $C_W(2k,\mathfrak{g})$ is obtained as special cases
of a general form of functional relations, which is another main
result of the present paper.   To explain this, we first introduce the
multi-variable version of Witten zeta-functions.

Let $\Delta$ be the set of all roots of ${\frak g}$, $\Delta_+$ the
set of all positive roots of ${\frak g}$,
$\fs=\{\alpha_1,\ldots,\alpha_r\}$ the fundamental system of $\Delta$,
$\alpha_j^{\vee}$ the coroot associated with $\alpha_j$ ($1\leq j\leq
r$).  Let $\lambda_1,\ldots,\lambda_r$ be the fundamental weights
satisfying
$\langle\alpha_i^{\vee},\lambda_j\rangle=\lambda_j(\alpha_i^{\vee})
=\delta_{ij}$ (Kronecker's delta).  A more explicit form of
$\zeta_W(s;{\frak g})$ can be written down in terms of roots and
weights by using Weyl's dimension formula (see (1.4) of \cite{KM2}).
Inspired by that form, we introduced in \cite{KM2} the multi-variable
version of Witten zeta-function
\begin{align}
  \label{1-3}
    \zeta_r({\bf s};{\frak g})=\sum_{m_1=1}^{\infty}\cdots
      \sum_{m_r=1}^{\infty}\prod_{\alpha\in\Delta_+}\langle\alpha^{\vee},
      m_1\lambda_1+\cdots+ m_r\lambda_r\rangle ^{-s_{\alpha}},
\end{align}
where ${\bf s}=(s_{\alpha})_{\alpha\in\Delta_+}\in {\Bbb C}^n$. In the
case that ${\frak g}$ is of type $X_r$, we call (\ref{1-3}) the
zeta-function of the root system of type $X_r$, and also denote it by
$\zeta_r({\bf s};X_r)$, where $X=A,B,C,D,E,F,G$.  Note that from (1.5)
and (1.7) in \cite{KM2}, we have
\begin{align}
  \label{1-3-2}
    \zeta_W(s;{\frak g})=K({\frak g})^s\zeta_r(s,\ldots,s;{\frak g}),
\end{align}
where
\begin{align}
  \label{1-3-3}
   K({\frak g})=\prod_{\alpha\in\Delta_+}\langle\alpha^{\vee},
           \lambda_1+\cdots+\lambda_r\rangle.
\end{align}
More generally, in \cite{KM2}, we introduced multiple zeta-functions
associated with sets of roots.  The main body of \cite{KM2} is devoted
to the study of recursive structures in the family of those
zeta-functions, which can be described in terms of Dynkin diagrams of
underlying root systems.

Now we discuss what are functional relations. 

The Euler-Zagier $r$-fold sum is defined by the multiple series
\begin{align}
  \label{1-4}
    \zeta_r(s_1,\ldots,s_r)=\sum_{m_1=1}^{\infty}\cdots
      \sum_{m_r=1}^{\infty}m_1^{-s_1}(m_1+m_2)^{-s_2}\times\\
       \times\cdots\times(m_1+\cdots+m_r)^{-s_r}.\notag
\end{align}
The harmonic product formula
\begin{align}
  \label{1-5}
    \zeta(s_1)\zeta(s_2)=\zeta(s_1+s_2)+\zeta_2(s_1,s_2)+\zeta_2(s_2,s_1)
\end{align}
due to L. Euler (where $\zeta(s)$ denotes the Riemann zeta-function),
and its $r$-ple analogue, are classically known (cf.~Bradley
\cite{Br}).  However, no other functional relations among multiple
zeta-functions has been discovered for a long time.

Let
\begin{align}
  \label{1-6}
    \zeta_{MT,2}(s_1,s_2,s_3)=\sum_{m=1}^{\infty}\sum_{n=1}^{\infty}
      m^{-s_1}n^{-s_2}(m+n)^{-s_3}.
\end{align}
This series is called Tornheim's harmonic double sum or the
Mordell-Tornheim double zeta-function, after the work of Tornheim
\cite{To} and Mordell \cite{Mo} in 1950s.  But it is to be noted that
this sum actually coincides with the Witten zeta-function \eqref{1-3}
for ${\frak g}={\frak{sl}}(3)$, that is, the simple Lie algebra of
type $A_2$.  Recently the third-named author \cite{TsC} proved that
there are certain functional relations among
$\zeta_{MT,2}(s_1,s_2,s_3)=\zeta_2(s_1,s_2,s_3;A_2)$ and the Riemann
zeta-function.  Moreover he obtained the same type of functional
relations for various relatives of $\zeta_{MT,2}(s_1,s_2,s_3)$ (see
\cite{TsA,TsJ}).  The method in these papers can be called the
``$u$-method'', because an auxiliary parameter $u>1$ was introduced to
ensure the absolute convergence in the argument.

In \cite{MTF}, by the same ``$u$-method'', the second and the
third-named authors proved certain functional relations among
$\zeta_3({\bf s};A_3)$, $\zeta_2({\bf s};A_2)$ and the Riemann
zeta-function.  The papers \cite{MTI,MTR,MNO,MNT} are also devoted to
the study of some new functional relations for certain (mainly) double
and triple zeta-functions and their relatives.

The above papers give many examples of functional relations.
Therefore, now, it is time appropriate to investigate whether there is
some structural reason underlying those functional relations, or not.
The first hint on this question was supplied by Nakamura's paper
\cite{Na1}, in which he presented a new simple proof of the result of
the third-named author \cite{TsC}. (Nakamura's method has then been
applied in \cite{MNO,MNT}.)

It can be observed from Nakamura's proof that the functional relations
proved in \cite{Na1,TsC} are connected with the symmetricity with
respect to the symmetric group ${\frak S}_3$, which is the Weyl group
associated with the Lie algebra of type $A_2$.  This suggests the
possibility of formulating general functional relations in terms of
Weyl groups of Lie algebras.

One of the main purposes of the present paper is to show that such general forms of
functional relations can indeed be proved.  In Section \ref{sec-2}, we
prepare some notation and preliminary results about root systems, Weyl
groups and convex polytopes, which play essential roles in the study
of the structural background of functional relations.  We will give
general forms of functional relations in Sections \ref{sec-3} and
\ref{sec-4}.  The most general form of functional relations is Theorem
\ref{thm:gen_func_rel}, which is specialized to the case of ``Weyl
group symmetric'' linear combinations $S({\bf s},{\bf y};I;\Delta)$
(defined by \eqref{eq:def_S}) of zeta-functions of root systems with
exponential factors in Theorem \ref{thm:FR} and Theorem
\ref{thm:val_rel}.  A na{\"\i}ve form of Theorem \ref{thm:val_rel} has
been announced in Section \ref{sec-3} of \cite{KMT}, but we
consider the generalized form with exponential factors because
this form can be applied to evaluations of $L$-functions of root systems
(see \cite{KM5} for the details.)

The theorems mentioned above give expressions of linear combinations
of zeta-functions in terms of certain multiple integrals involving
Lerch zeta-functions.  Since the values of Lerch zeta-functions at
positive integers ($\geq 2$) can be written as Bernoulli polynomials,
we can show more explicit forms of those multiple integrals in some
special cases.  We will carry out this procedure by using generating
functions (Theorem \ref{thm:func_rel}).

In particular, we find that the value $S({\bf s},{\bf
  y};\Delta)=S({\bf s},{\bf y};\emptyset;\Delta)$ at ${\bf s}={\bf
  k}=(k_{\alpha})$, where all $k_{\alpha}$'s are positive integers
($\geq 2$), is essentially a generalization $P({\bf k},{\bf
  y};\Delta)$ of periodic Bernoulli functions.  The generating
function $F({\bf t},{\bf y};\Delta)$ of $P({\bf k},{\bf y};\Delta)$
will be evaluated in Theorem \ref{thm:gen_func}. Consequently, we
can prove a generalization of Witten's volume formula \eqref{1-2}
(Theorem \ref{thm:W-Z}).

In Section \ref{sec-5}, we will show the Weyl group symmetry of
$S({\bf s},{\bf y};\Delta)$, $F({\bf t},{\bf y};\Delta)$ and
$P({\bf k},{\bf y};\Delta)$.    For our purpose it is not sufficient to
consider the usual Weyl group only, and hence we will introduce the extended affine
Weyl group $\widehat{W}$, and will prove the symmetricity with respect to
$\widehat{W}$ (Theorems \ref{thm:sym_S}, \ref{thm:gen_sym} and \ref{thm:inv_qv}).    These results ensure that the
existence of functional relations is indeed based on the symmetry with respect
to $\widehat{W}$.
It is to be noted that Weyl groups
already played a role in Zagier's sketch \cite{Za} of the proof of \eqref{1-2}.

In Section \ref{sec-6} we will prove that $P({\bf k},{\bf y};\Delta)$
can be continued to polynomials in ${\bf y}$ (Theorem
\ref{thm:Bernoulli}).  This may be regarded as a (root-system
theoretic) generalization of Bernoulli polynomials.  Since $P({\bf
  k},{\bf y};\Delta)$ is essentially the same as $S({\bf k},{\bf
  y};\Delta)$, this gives explicit relations among special values of
zeta-functions of root systems.  Moreover in the same theorem we will
give the continuation of $F({\bf t},{\bf y};\Delta)$.  

As examples, in
Section \ref{sec-7} we will calculate $P$ and $F$ explicitly in the
cases of $A_1$, $A_2$, $A_3$ and $C_2(\simeq B_2)$.  In particular, from the
explicit expansion of generating functions $F$, we will determine the
value of $C_{W}(2,A_2)$, $C_{W}(2,A_3)$ and $C_{W}(2,C_2)$ in
\eqref{1-2}.  The cases of $C_{W}(2,A_2)$, $C_{W}(2,A_3)$ are included
in the results of Mordell \cite{Mo}, Subbarao et al. \cite{SS} and
Gunnells-Sczech \cite{GS}, while the determination of the value of
$C_{W}(2,C_2)$ seems to be a new result.  By the same method it is
possible to obtain the explicit values of $C_{W}(2k,A_2)$,
$C_{W}(2k,A_3)$, $C_{W}(2k,C_2)$ for any $k\in {\Bbb N}$.  More
generally, we can deduce explicit value relations and special values
of zeta-functions of any simple algebra ${\frak g}$ by the same
argument, at least in principle, though the actual procedure will
become quite complicated when the rank of ${\frak g}$ becomes higher.
We also provide an example of a functional relation in the $A_2$ case.

Some parts of the contents of \cite{KM2} and the present paper have been already 
announced briefly in \cite{KMT,KMT2,KMTpja}.

\section{Root systems, Weyl groups and convex polytopes}\label{sec-2}

In this preparatory section, we first fix notation and summarize basic facts about root systems
and Weyl groups. See \cite{Hum,Hum72,Kac,Bourbaki} for the details.
Let $V$ be an $r$-dimensional real vector space equipped with an inner product $\langle \cdot,\cdot\rangle$.
We denote the norm of $v\in V$ by $\norm{v}=\langle v,v\rangle^{1/2}$.
The dual space $V^*$ is identified with $V$ via the inner product of $V$.
Let $\Delta$ be a finite reduced root system in $V$ and
$\fs=\{\alpha_1,\ldots,\alpha_r\}$ its fundamental system.
Let 
$\Delta_+$ and $\Delta_-$ be the set of all positive roots and negative roots respectively.
Then we have a decomposition of the root system $\Delta=\Delta_+\coprod\Delta_-$.
Let $Q^\vee$ be the coroot lattice,
$P$ the weight lattice,
$P_+$ the set of integral dominant weights 
and
$P_{++}$ the set of integral strongly dominant weights
respectively defined by
\begin{equation}
  Q^\vee=\bigoplus_{i=1}^r\mathbb{Z}\,\alpha^\vee_i,\qquad
  P=\bigoplus_{i=1}^r\mathbb{Z}\,\lambda_i,\qquad
  P_+=\bigoplus_{i=1}^r\mathbb{N}_0\,\lambda_i,\qquad
  P_{++}=\bigoplus_{i=1}^r\mathbb{N}\,\lambda_i,
\end{equation}
where the fundamental weights $\{\lambda_j\}_{j=1}^r$
are a basis dual to $\fs^\vee$
satisfying $\langle \alpha_i^\vee,\lambda_j\rangle=\delta_{ij}$.
Let
\begin{equation}
\rho=\frac{1}{2}\sum_{\alpha\in\Delta_+}\alpha=\sum_{j=1}^r\lambda_j
\end{equation}
be the lowest strongly dominant weight.
Then $P_{++}=P_++\rho$.

We define the reflection $\sigma_\alpha$ with respect to a root $\alpha\in\Delta$ as
\begin{equation}
 \sigma_\alpha:V\to V, \qquad \sigma_\alpha:v\mapsto v-\langle \alpha^\vee,v\rangle\alpha
\end{equation}
and for a subset $A\subset\Delta$, let
$W(A)$ be the group generated by reflections $\sigma_\alpha$ for $\alpha\in A$.
Let $W=W(\Delta)$ be the Weyl group.
Then $\sigma_j=\sigma_{\alpha_j}$ ($1\leq j \leq r$) generates $W$. Namely we have $W=W(\fs)$.
Any two fundamental systems $\fs$, $\fs'$ are conjugate under $W$.

We denote the fundamental domain called
the fundamental Weyl chamber by
\begin{equation}
C=\{v\in V~|~\langle \fs^\vee,v\rangle\geq0\},
\end{equation}
where $\langle\fs^\vee,v\rangle$ means any of
$\langle\alpha^\vee,v\rangle$ for $\alpha^\vee\in\fs^\vee$.
Then $W$ acts on the set of Weyl chambers $WC=\{wC~|~w\in W\}$ simply transitively.
Moreover if $wx=y$ for $x,y\in C$, then $x=y$ holds.
The stabilizer $W_x$ of a point $x\in V$ is generated 
by the reflections which stabilize $x$. 
We see that $P_+=P\cap C$.

Let $\Aut$ be the subgroup of all the automorphisms $\GL(V)$ which stabilizes $\Delta$
(see \cite[\S 12.2]{Hum72}). 
Then the Weyl group $W$ is a normal subgroup of $\Aut$ and
there exists a subgroup $\Omega\subset \Aut$ such that
\begin{equation}
\Aut=\Omega\ltimes W.
\end{equation}
The subgroup $\Omega$ is isomorphic to the group $\mathrm{Aut}(\Gamma)$ of automorphisms
of the Dynkin diagram $\Gamma$.
The group $\Aut$ is called the extended Weyl group.
For $w\in \Aut$, we set 
\begin{equation}
\Delta_w=\Delta_+\cap w^{-1}\Delta_-
\end{equation}
and the length function $\ell(w)=\abs{\Delta_w}$ (see \cite[\S 1.6]{Hum}).
The subgroup $\Omega$ is characterized as $w\in \Omega$ if and only if $\ell(w)=0$.
Note that $w\Delta_w=\Delta_-\cap w\Delta_+=-\Delta_{w^{-1}}$ and $\ell(w)=\ell(w^{-1})$.

For $u\in V$, let $\tau(u)$ be the translation by $u$, that is,
\begin{equation}
  \tau(u):V\to V, \qquad \tau(u):v\mapsto v+u.
\end{equation}
Since $\Aut$ stabilizes the coroot lattice $Q^\vee$, we can define
\begin{equation}
 \widehat{W}=\Aut\ltimes \tau(Q^\vee).
\end{equation}
Then $\widehat{W}=(\Omega\ltimes W)\ltimes\tau(Q^\vee)
  \simeq \Omega\ltimes(W\ltimes\tau(Q^\vee))$.
We call $\widehat{W}$ the extended affine Weyl group in this paper
(see \cite{Kac,Bourbaki} for the details of affine Weyl groups).
It should be noted that there are some other groups which are, in some references,
called the extended affine Weyl group.

Let $\widehat{V}=V\times\mathbb{R}$ and $\delta=(0,1)\in\widehat{V}$. We embed $V$ in $\widehat{V}$ and we have
$\widehat{V}=V\oplus\mathbb{R}\,\delta$.
For $\gamma=\eta+c\delta\in\widehat{V}$ with $\eta\in V$ and $c\in\mathbb{R}$,
we associate an affine linear functional on $V$ as
$\gamma(v)=\langle \eta,v\rangle+c$.
Let $\widehat{Q}^\vee$ be the affine coroot lattice defined by
\begin{equation}
\widehat{Q}^\vee=Q^\vee\oplus\mathbb{Z}\,\delta
\end{equation}
(see \cite{Kac}).
For a set $X$, 
let $\Func(X)$ be the set of all functions $f:X\to\mathbb{C}$.
For a function $f\in \Func(P)$, we define a subset
\begin{equation}
 H_f=\{\lambda\in P~|~f(\lambda)=0\}
\end{equation}
and for a subset $A$ of $\Func(P)$, define
$H_A=\bigcup_{f\in A}H_f$.
One sees that an action of $W$ is induced on $\Func(P)$ 
as $(wf)(\lambda)=f(w^{-1}\lambda)$.
Note that $V\subset\widehat{V}\subset\Func(P)$, where the second inclusion is given by
the associated functional mentioned above.

Let $I\subset \{1,\ldots,r\}$ and $\fs_I=\{\alpha_i~|~i\in I\}\subset\fs$.
Let $V_I$ be the linear subspace spanned by $\fs_I$.
Then $\Delta_I=\Delta\cap V_I$
is a root system in $V_I$ whose fundamental system is $\fs_I$. 
For the root system $\Delta_I$, 
we denote the corresponding
coroot lattice (resp.~weight lattice etc.)
 by $Q_I^\vee=\bigoplus_{i\in I}\mathbb{Z}\,\alpha_i^\vee$ 
(resp.~$P_I=\bigoplus_{i\in I}\mathbb{Z}\,\lambda_i$ etc.).
We define 
\begin{equation}
  C_I
  =\{v\in C~|~\langle\fs_{I^c}^\vee,v\rangle=0,\quad
  \langle\fs_I^\vee,v\rangle>0\},
\end{equation}
where $I^c$ is the complement of $I$.
Then the dimension of the linear span of $C_I$ is $\abs{I}$ and
we have a disjoint union
\begin{equation}
\label{eq:Coxeter_comp_C}
 C=\coprod_{I\subset\{1,\ldots,r\}}C_I 
\end{equation}
and the collection of all sets $wC_I$ 
for $w\in W$ and $I\subset\{1,\ldots,r\}$
is called the Coxeter complex (see \cite[\S 1.15]{Hum};
it should be noted that we use a little different notation),
which partitions $V$ and we have a decomposition
\begin{equation}
\label{eq:Coxeter_comp_P}
  P_+=\coprod_{I\subset\{1,\ldots,r\}}P_{I++},
\end{equation}
where
\begin{equation}
  P_{I++}=P_+\cap C_I.
\end{equation}
In particular, $P_{\emptyset++}=\{0\}$ and $P_{\{1,\ldots,r\}++}=P_{++}$.

The natural embedding $\iota:Q_I^\vee\to Q^\vee$ induces the projection
$\iota^*:P\to P_I$. Namely for $\lambda\in P$,
$\iota^*(\lambda)$ is defined as a unique element of $P_I$ satisfying
 $\langle \iota(q),\lambda\rangle=\langle q,\iota^*(\lambda)\rangle$ for all $q\in Q_I^\vee$.
Let
\begin{equation}
\label{eq:def_W_I}
  W^I=\{w\in W~|~\Delta^\vee_{I+}\subset w\Delta^\vee_+\}.
\end{equation}
Then we have the following key lemmas to functional relations among zeta-functions.
Note that the statements hold trivially in the case $I=\emptyset$ and 
hence we deal with $I\neq\emptyset$ in their proofs.
\begin{lemma}
\label{lm:Z_mcr}
  The subset
  $W^I$ coincides with the minimal (right) coset representatives 
  $\{w\in W~|~\ell(\sigma_iw)>\ell(w)\text{ for all }i\in I\}$
  of the parabolic subgroup $W(\Delta_I)$ (see \cite[\S 1.10]{Hum}).
  Therefore $\abs{W^I}=\bigl(W(\Delta):W(\Delta_I)\bigr)$.
\end{lemma}
\begin{proof}
  Let $w\in W^I$. Then $\Delta^\vee_{I+}\subset w\Delta^\vee_+$,
which implies 
$\Delta^\vee_{I+}\cap w\Delta^\vee_-=\emptyset$.
In particular, $\alpha^\vee_i\not\in w\Delta^\vee_-$ for $i\in I$, 
which yields $\Delta^\vee_+\cap w\Delta^\vee_-=
(\Delta^\vee_+\setminus\{\alpha_i^\vee\})\cap w\Delta^\vee_-$.
Therefore
\begin{equation}
\sigma_i(\Delta_{+}^{\vee}\cap w\Delta_{-}^{\vee})= 
  \sigma_i\bigl((\Delta_+^\vee\setminus\{\alpha_i^\vee\})\cap w\Delta_-^\vee\bigr)=
(\Delta^\vee_+\setminus\{\alpha_i^\vee\})\cap\sigma_iw\Delta_-^\vee\subset
\Delta^\vee_+\cap\sigma_iw\Delta_-^\vee
\end{equation}
and
$\ell(\sigma_iw)\geq\ell(w)$.
Since $\ell(\sigma_iw)=\ell(w)\pm1$, we have
$\ell(\sigma_iw)=\ell(w)+1$ and $w$ is a minimal coset representative.

  Assume that $w\in W$ satisfies $\ell(\sigma_iw)>\ell(w)$ for all $i\in I$.
Then we have
\begin{equation}
\sigma_i(\Delta_{+}^{\vee}\cap w\Delta_{-}^{\vee})=
\bigl(
(\Delta^\vee_+\setminus\{\alpha_i^\vee\})\cap\sigma_iw\Delta_-^\vee
\bigr)
\cup
\bigl(
\{-\alpha_i\}\cap\sigma_iw\Delta_-^\vee
\bigr).
\end{equation}
Since 
$\abs{\sigma_i(\Delta_{+}^{\vee}\cap w\Delta_{-}^{\vee})}=\ell(w)$ and
$\abs{(\Delta^\vee_+\setminus\{\alpha_i^\vee\})\cap\sigma_iw\Delta_-^\vee}\geq\ell(\sigma_iw)-1=\ell(w)$,
we have $\abs{\{-\alpha_i\}\cap\sigma_iw\Delta_-^\vee}=0$.
It implies that no element of $\Delta^\vee_+\cap w\Delta^\vee_-$ is sent to $\Delta^\vee_-$
by $\sigma_i$ for $i\in I$, and hence $\fs^\vee_I\cap w\Delta^\vee_-=\emptyset$.
Since $0\not\in\langle\Delta^\vee,\rho\rangle$ and
$\alpha^\vee\in w\Delta^\vee_-$ if and only if 
$\langle\alpha^\vee,w\rho\rangle<0$,
we have
$\langle\fs^\vee_I,w\rho\rangle>0$
and hence $\langle\Delta^\vee_{I+},w\rho\rangle>0$.
It follows that $\Delta^\vee_{I+}\cap w\Delta^\vee_-=\emptyset$
and $w\in W^I$.
\end{proof}
\begin{lemma}
\label{lm:pp_wp}
\begin{equation}
\iota^{*-1}(P_{I+})=P_{I+}\oplus P_{I^c}=\bigcup_{w\in W^I}w P_+.
\end{equation}
\end{lemma}
\begin{proof}
The first equality is clear. Prove the second equality.

Assume $w\in W^I$.
Then for $\lambda\in P_+$,
we have 
$
 \langle \Delta^\vee_{I+},w\lambda\rangle
=\langle w^{-1}\Delta^\vee_{I+},\lambda\rangle
\subset\langle \Delta^\vee_+,\lambda\rangle
\geq0.   
$
Hence $w P_+\subset\iota^{*-1}(P_{I+})$.

Conversely, assume $\lambda\in\iota^{*-1}(P_{I+})$.
Since $\abs{\Delta^\vee}<\infty$, it is possible to fix 
a sufficiently small constant $c>0$ such that
 $0<\abs{\langle\Delta^\vee,c\rho\rangle}<1$. 
Then we see that
$\lambda+c\rho$ is regular (see \cite[\S10.1]{Hum72}), i.e.,
 $0\not\in\langle\Delta^\vee,\lambda+c\rho\rangle$
and the signs of
$\langle\alpha^\vee,\lambda\rangle$ and $\langle\alpha^\vee,\lambda+c\rho\rangle$
coincide if $\langle\alpha^\vee,\lambda\rangle\neq 0$,
because $\langle\Delta^\vee,\lambda\rangle\subset\mathbb{Z}$.
Let $\Tilde{\Delta}^\vee_+=
\{\alpha^\vee\in\Delta^\vee~|~\langle\alpha^\vee,\lambda+c\rho\rangle>0\}$.
Then $\Tilde{\Delta}^\vee_+$ is a positive system and
hence there exists an element $w\in W$ such that
$\Tilde{\Delta}^\vee_+=w\Delta^\vee_+$.
Since $\lambda\in\iota^{*-1}(P_{I+})$,
we have $\Delta^\vee_{I+}\subset\Tilde{\Delta}^\vee_+$.
Hence
$\Delta^\vee_{I+}\subset w\Delta^\vee_+$, that is,
$w\in W^I$.
Moreover
 $\langle\Tilde{\Delta}^\vee_+,\lambda+c\rho\rangle>0$ implies
$\langle \Delta^\vee_+,w^{-1}(\lambda+c\rho)\rangle>0$
and
$\langle \Delta^\vee_+,w^{-1}\lambda\rangle\geq0$
again due to the integrality.
Therefore $\lambda\in wP_+$.
\end{proof}
\begin{lemma}
\label{lm:wp_uniq}
For $\lambda\in\iota^{*-1}(P_{I+})$, an element $w\in W^I$ satisfying
$\lambda\in wP_+$ (whose existence is assured by Lemma \ref{lm:pp_wp})
is unique
if and only if $\lambda\not\in H_{\Delta^\vee\setminus\Delta^\vee_I}$. 
\end{lemma}
\begin{proof}
Assume $\alpha^\vee\in\Delta^\vee\setminus\Delta^\vee_I$ and $\lambda\in\iota^{*-1}(P_{I+})\cap H_{\alpha^\vee}$.
Let $w\in W^I$ satisfy $\lambda\in wP_+$.
Then $\sigma_\alpha\lambda=\lambda\in wP_+$ and hence
$w^{-1}\lambda=\sigma_{w^{-1}\alpha}w^{-1}\lambda\in P_+$, which further implies 
$w^{-1}\alpha^\vee\in\Delta'^\vee$,
where $\Delta'^\vee$ is a coroot system orthogonal to
 $w^{-1}\lambda$ whose fundamental system is given
by $\fs'^\vee=\{\alpha^\vee_i\in\fs^\vee~|~\langle\alpha^\vee_i,w^{-1}\lambda\rangle=0\}$ (see \cite[\S 1.12]{Hum}).
If $\fs'^\vee\subset w^{-1}\Delta^\vee_I$,
then $W(\fs'^\vee)\fs'^\vee=\Delta'^\vee\subset w^{-1}\Delta^\vee_I$,
and hence $w^{-1}\alpha^\vee\in w^{-1}\Delta^\vee_I$,
which contradicts to the assumption
$\alpha^\vee\not\in\Delta^\vee_I$. Therefore there exists a fundamental coroot 
$\alpha^\vee_i\in\fs'^\vee\setminus w^{-1}\Delta^\vee_I$, 
which satisfies $\sigma_iw^{-1}\lambda=w^{-1}\lambda\in P_+$ by construction.
Since $w\in W^I$, 
we have $w^{-1}\Delta^\vee_{I+}\subset\Delta^\vee_+\setminus\{\alpha_i^\vee\}$.
Hence $\sigma_iw^{-1}\Delta^\vee_{I+}\subset\Delta^\vee_+$, because $\sigma_i(\Delta^\vee_+\setminus\{\alpha_i^\vee\})\subset\Delta^\vee_+$.
Then putting $w'=w\sigma_i$, we have $W^I\ni w'\neq w$ such that $\lambda\in wP_+\cap w'P_+$.

Conversely, 
assume that there exist $w,w'\in W^I$ such that $w\neq w'$ and
 $\lambda\in wP_+\cap w'P_+$.
This implies that $w^{-1}\lambda=w'^{-1}\lambda$ is on a
wall of $C$ and hence $\lambda\in H_{\Delta^\vee}$.
Let $\Delta''^\vee=\{\alpha^\vee\in\Delta^\vee~|~\lambda\in H_{\alpha^\vee}\}$ be a coroot system orthogonal to $\lambda$
so that
$\lambda\in H_{\Delta''^\vee}$.
Assume $\Delta''^\vee\subset\Delta^\vee_I$. 
Then by $\lambda=ww^{\prime-1}\lambda$,
we have
 $ww^{\prime-1}\in W_\lambda$ and hence
 $ww^{\prime-1}\in W(\Delta_I)$, 
because $W_\lambda=W(\Delta''^\vee)\subset W(\Delta^\vee_I)$ 
by the assumption.
Since $\id\neq ww^{\prime-1}\in W(\Delta_I)$,
there exists a coroot $\alpha^\vee\in\Delta^\vee_{I+}$ such that
$\beta^\vee=ww^{\prime-1}\alpha^\vee\in\Delta^\vee_{I-}$.
Then,
since 
$w^{-1}(ww^{\prime-1})\Delta^\vee_{I+}\subset\Delta^\vee_+$ and
$w^{-1}\Delta^\vee_{I+}\subset\Delta^\vee_+$,
we have
$w^{-1}\beta^\vee\in\Delta^\vee_+$ from the first inclusion
and $w^{-1}(-\beta^\vee)\in\Delta^\vee_+$ from the second one, which leads to the contradiction.
Therefore $\lambda\in H_{\alpha^\vee}$ for $\alpha\in\Delta''^\vee\setminus\Delta^\vee_I$.
\end{proof}

Next we give some 
definitions and
facts about convex polytopes (see \cite{Zie95,Hibi92}) and their triangulation.
For a subset $X\subset\mathbb{R}^N$, we denote by $\Conv(X)$ the convex hull of $X$.
A subset $\mathcal{P}\subset\mathbb{R}^N$ is called a convex polytope 
if $\mathcal{P}=\Conv(X)$ for some finite subset $X\subset\mathbb{R}^N$.
Let $\mathcal{P}$ be a $d$-dimensional polytope.
Let $\mathcal{H}$ be a hyperplane in $\mathbb{R}^N$.
Then $\mathcal{H}$ divides $\mathbb{R}^N$ into two half-spaces.
If $\mathcal{P}$ is entirely contained in one of the two closed half-spaces
and $\mathcal{P}\cap\mathcal{H}\neq\emptyset$, then $\mathcal{H}$ is
called a supporting hyperplane.
For a supporting hyperplane $\mathcal{H}$,
a subset $\mathcal{F}=\mathcal{P}\cap\mathcal{H}\neq\emptyset$
is called a face of the polytope $\mathcal{P}$.
If the dimension of a face $\mathcal{F}$ is $j$, then we call it a $j$-face $\mathcal{F}$. 
A $0$-face is called a vertex and a $(d-1)$-face a facet.
For convenience, 
we regard $\mathcal{P}$ itself as its unique $d$-face.
Let $\Vrt(\mathcal{P})$ be the set of the vertices of $\mathcal{P}$. 
Then 
\begin{equation}
\label{eq:prop_F}
  \mathcal{F}=\Conv(\Vrt(\mathcal{P})\cap\mathcal{F}),
\end{equation}
for a face $\mathcal{F}$.

It is known that 
triangulation of a convex polytope can be executed without adding any vertices.
Here we give an explicit procedure of a triangulation of $\mathcal{P}$.
Number all the vertices of $\mathcal{P}$ as $\bd{p}_1,\ldots,\bd{p}_k$.
For a face $\mathcal{F}$, 
by $\mathcal{N}(\mathcal{F})$
we mean the vertex $\bd{p}_j$ 
whose index $j$ is the smallest
in the vertices belonging to $\mathcal{F}$.
A full flag $\Phi$ is defined by the sequence
  \begin{equation}
    \Phi:\mathcal{F}_0\subset\mathcal{F}_1\subset\cdots\subset\mathcal{F}_{d-1}\subset
\mathcal{F}_d=
\mathcal{P},
  \end{equation}
with $j$-faces $\mathcal{F}_j$ such that
$\mathcal{N}(\mathcal{F}_j)\not\in\mathcal{F}_{j-1}$.
\begin{theorem}[\cite{Stan80}]
\label{thm:simp_div}
All the collection of the simplexes
with vertices $\mathcal{N}(\mathcal{F}_0)$, \dots, $\mathcal{N}(\mathcal{F}_{d-1})$,
 $\mathcal{N}(\mathcal{F}_d)$
associated with full flags gives a triangulation.
\end{theorem}
\begin{remark}
\label{rem:face_poset}
  This procedure only depends on its face poset structure (see \cite[\S5]{Hibi92}).
\end{remark}

For $\bd{a}={}^t(a_1,\ldots,a_N),\bd{b}={}^t(b_1,\ldots,b_N)\in\mathbb{C}^N$
we define $\bd{a}\cdot \bd{b}=a_1b_1+\cdots+a_Nb_N$.
The definition of polytopes above is that of ``V-polytopes''.
We mainly deal with another representation of polytopes,
``H-polytopes'', instead. Namely, consider a bounded subset of the form
\begin{equation}
  \mathcal{P}=\bigcap_{i\in I}\mathcal{H}_i^+\subset\mathbb{R}^N,
\end{equation}
where $\abs{I}<\infty$ and $\mathcal{H}_i^+=\{\bd{x}\in\mathbb{R}^N~|~\bd{a}_i\cdot \bd{x}\geq h_i\}$
with $\bd{a}_i\in\mathbb{R}^N$ and $h_i\in\mathbb{R}$.
The following theorem is intuitively clear but nontrivial (see, for example, \cite[Theorem 1.1]{Zie95}).
\begin{theorem}[Weyl-Minkowski]
  H-polytopes are V-polytopes and vice versa.
\end{theorem}
We have a representation of $k$-faces in terms of hyperplanes
 $\mathcal{H}_i=\{\bd{x}\in\mathbb{R}^N~|~\bd{a}_i\cdot \bd{x}=h_i\}$.
\begin{proposition}
\label{prop:k_face_1}
  Let $J\subset I$.
  Assume that
  $\mathcal{F}=\mathcal{P}\cap\bigcap_{j\in J}\mathcal{H}_j\neq\emptyset$.
  Then $\mathcal{F}$ is a face.
\end{proposition}
\begin{proof}
Let $\bd{x}\in\mathcal{P}$. Then
$\bd{x}\in\bigcap_{j\in J}\mathcal{H}_j^+$ and hence
$\bd{a}_j\cdot \bd{x}\geq h_j$ for all $j\in J$.
Set $\bd{a}=\sum_{j\in J}\bd{a}_j$ and $h=\sum_{j\in J}h_j$.
Let $\mathcal{H}$ be a hyperplane defined by $\{\bd{x}\in\mathbb{R}^N~|~\bd{a}\cdot\bd{x}=h\}$.
Then $\bd{a}\cdot\bd{x}\geq h$ for $\bd{x}\in\mathcal{P}$ and
 $\mathcal{P}\subset\mathcal{H}^+=\{\bd{x}\in\mathbb{R}^N~|~\bd{a}\cdot\bd{x}\geq h\}$.

Let $\bd{x}\in\mathcal{P}\cap\mathcal{H}$.
Then $\bd{a}\cdot\bd{x}=h$.
Since $\bd{a}_j\cdot\bd{x}=h_j+c_j\geq h_j$ with some $c_j\geq0$ for all $j\in J$,
we have $c_j=0$ and $\bd{x}\in\mathcal{H}_j$
for all $j\in J$. Thus $\bd{x}\in\mathcal{F}$.
It is easily seen that $\bd{x}\in\mathcal{F}$ satisfies $\bd{x}\in\mathcal{P}\cap\mathcal{H}$.
Therefore $\mathcal{F}=\mathcal{P}\cap\mathcal{H}$ and
$\mathcal{H}$ is a supporting hyperplane.
\end{proof}
\begin{proposition}
\label{prop:k_face_2}
Let $\mathcal{H}$ be a supporting hyperplane and
$\mathcal{F}=\mathcal{P}\cap\mathcal{H}$ a $k$-face.
Then there exists a set of indices $J\subset I$ such that
$\abs{J}=(\dim\mathcal{P})-k$ and
$\mathcal{F}=\mathcal{P}\cap\bigcap_{j\in J}\mathcal{H}_j$.
\end{proposition}
\begin{proof}
Assume $d=N$ without loss of generality.
Let $\bd{x}\in\mathcal{F}$.
Then $\bd{x}\in\mathcal{H}$ and hence $\bd{x}\in\partial\mathcal{P}$ since $\mathcal{P}\subset\mathcal{H}^+$.
If $\bd{x}\in\mathcal{H}_i^+\setminus\mathcal{H}_i$ for all $i\in I$, then $\bd{x}$ is in the interior of
$\mathcal{P}$, which contradicts 
to the above.
Thus $\bd{x}\in\mathcal{H}_j$ for some $j\in I$.

First
we assume that $\mathcal{F}=\mathcal{P}\cap\mathcal{H}$ is a facet.
Then there exist a subset
$\{\bd{x}_1,\ldots,\bd{x}_N\}\subset\mathcal{F}$ such that $\bd{x}_2-\bd{x}_1,\ldots,\bd{x}_N-\bd{x}_1$ are 
linearly independent.
Let $\mathcal{C}$ be the convex hull of $\{\bd{x}_1,\ldots,\bd{x}_N\}$.
We consider that $\mathcal{C}\subset\mathcal{F}$ is equipped with the relative topology. 
Note that
for $\bd{x}\in\mathcal{C}$, we have $\bd{x}\in\mathcal{H}_{i(\bd{x})}$, where $i(\bd{x})\in I$.
We show that
there exist an open subset $\mathcal{U}\subset\mathcal{C}$
and $i\in I$ such that
$\mathcal{H}_i\cap \mathcal{U}$ is dense in $\mathcal{U}$.
Fix an order $\{i_1,i_2,\ldots\}=I$.
If $\mathcal{H}_{i_1}\cap\mathcal{C}$ is dense in $\mathcal{C}$, then we have done.
Hence we assume that it is false. 
Then there exists an open subset $\mathcal{U}_1\subset\mathcal{C}$
such that $\mathcal{H}_{i_1}\cap\mathcal{U}_1=\emptyset$.
Similarly we see that
 $\mathcal{H}_{i_2}\cap\mathcal{U}_1$ is dense in $\mathcal{U}_1$ unless
there exists an open subset $\mathcal{U}_2\subset\mathcal{U}_1$
such that $\mathcal{H}_{i_2}\cap\mathcal{U}_2=\emptyset$.
Since $\abs{I}<\infty$, repeated application of this argument yields
the assertion.
Thus there exists a subset
$\{\bd{x}_1',\ldots,\bd{x}_N'\}\subset\mathcal{H}_i\cap\mathcal{U}$ such that $\bd{x}_2'-\bd{x}_1',\ldots,\bd{x}_N'-\bd{x}_1'$ are 
linearly independent.
Hence we have $\mathcal{F}\subset\mathcal{H}_i$ and $\mathcal{H}=\mathcal{H}_i$.

For any $k$-face $\mathcal{F}$, there exists a sequence of faces such that
\begin{equation}
  \mathcal{F}=\mathcal{F}_k\subset\mathcal{F}_{k+1}\subset\cdots\subset\mathcal{F}_{N-1}\subset\mathcal{P},
\end{equation}
where $\mathcal{F}_j$ ($k\leq j\leq N-1$) is a $j$-face.
Since $\mathcal{F}_j$ is a facet of $\mathcal{F}_{j+1}$, by the induction on dimensions,
we have $\mathcal{F}=\mathcal{P}\cap\bigcap_{j\in J}\mathcal{H}_j$ for some $J\subset I$ with $\abs{J}=N-k$.
\end{proof}

\begin{lemma}
\label{lm:simplex}
Let $\bd{a}={}^t(a_1,\ldots,a_N)\in\mathbb{C}^N$ and $\sigma$ be a simplex with vertices
  $\bd{p}_0,\dots,\bd{p}_N\in\mathbb{R}^N$ in general position.
Then
\begin{equation}
  \label{eq:simplex}
  \begin{split}
    \int_\sigma e^{\bd{a}\cdot \bd{x}}d\bd{x}
    &=N!\Vol(\sigma)\sum_{m=0}^N\frac{e^{\bd{a}\cdot\bd{p}_m}}{\prod_{j\neq m}\bd{a}\cdot (\bd{p}_m-\bd{p}_j)}\\
&=N!\Vol(\sigma) T(\bd{a}\cdot \bd{p}_0,\ldots,\bd{a}\cdot \bd{p}_N),\\
  \end{split}
\end{equation}
where
\begin{gather}
\label{eq:det_form_simplex}
T(x_0,\ldots,x_N)=\det
\begin{pmatrix}
  1&\cdots&1\\
  x_0&\cdots& x_N\\
  \vdots & \ddots&\vdots \\
  x_0^{N-1}&\cdots& x_N^{N-1} \\
  e^{x_0}&\cdots& e^{x_N}
\end{pmatrix}
\Bigg/
\det
\begin{pmatrix}
  1&\cdots&1\\
  x_0&\cdots& x_N\\
  \vdots & \ddots&\vdots \\
  x_0^{N-1}&\cdots& x_N^{N-1} \\
  x_0^N&\cdots& x_N^N
\end{pmatrix}
,\\
\label{eq:vol_simplex} 
  \Vol(\sigma)=\int_\sigma 1d\bd{x}=\frac{1}{N!}\abs{\det\Tilde{P}},
\end{gather}
and 
\begin{equation}
  \Tilde{P}=
  \begin{pmatrix}
    1& \cdots& 1\\
    \bd{p}_0&\cdots& \bd{p}_N
  \end{pmatrix},
\end{equation}
with $\bd{p}_j$ regarded as column vectors.
\end{lemma}
\begin{proof}
Let
$\tau$ be the simplex whose vertices consist of the origin and
${}^t(0,\dots,\overset{j}{1},\dots,0)$, 
$1\leq j\leq N$. Then by changing variables
from $\bd{x}$ to $\bd{y}$ as
$\bd{x}=\bd{p}_0+P \bd{y}$ with the $N\times N$ matrix $P=(\bd{p}_1-\bd{p}_0,\dots,\bd{p}_N-\bd{p}_0)$,
we have
\begin{equation}
\label{eq:simplex_orig} 
  \int_\sigma e^{\bd{a}\cdot\bd{x}}d\bd{x}=
  e^{\bd{a}\cdot \bd{p}_0}\abs{\det P}
  \int_\tau e^{\bd{b}\cdot \bd{y}}d\bd{y},
\end{equation}
where
${}^t \bd{b}=(b_1,\ldots,b_N)={}^t \bd{a} P$.
Since
\begin{multline}
  \int_0^{1-y_1-\cdots-y_{m-1}}e^{c+(b_1-c)y_1+\cdots+(b_m-c)y_m}dy_m
=\\
\frac{1}{b_m-c}(
e^{b_m+(b_1-b_m)y_1+\cdots+(b_{m-1}-b_m)y_{m-1}}
-
e^{c+(b_1-c)y_1+\cdots+(b_{m-1}-c)y_{m-1}}
),
\end{multline}
we see that
\begin{equation}
    \int_\tau e^{\bd{b}\cdot \bd{y}}d\bd{y}=q_0(\bd{b})+\sum_{m=1}^N q_m(\bd{b})e^{b_m},
\end{equation}
where $q_j(\bd{b})$ for $0\leq j\leq N$ are rational functions in $b_1,\ldots,b_N$.
In particular, we have 
\begin{equation}
  q_0(\bd{b})=\frac{(-1)^N}{b_1\cdots b_N}
\end{equation}
and hence
\begin{equation}
  \int_\sigma e^{\bd{a}\cdot \bd{x}}d\bd{x}=
\frac{\abs{\det(\bd{p}_0-\bd{p}_1,\dots,\bd{p}_0-\bd{p}_N)}}
{\bigl(\bd{a}\cdot (\bd{p}_0-\bd{p}_1)\bigr)\cdots\bigl(\bd{a}\cdot(\bd{p}_0-\bd{p}_N)\bigr)}e^{\bd{a}\cdot \bd{p}_0}
+\sum_{m=1}^N \Tilde{q}_m(\bd{a})e^{\bd{a}\cdot \bd{p}_m},
\end{equation}
where $\Tilde{q}_m(\bd{a})$ are certain rational functions in $a_1,\ldots,a_N$.
Since
 $e^{\bd{a}\cdot \bd{p}_j}$ 
for $0\leq j\leq N$ are 
linearly independent over the field of 
rational functions in $a_1,\ldots,a_N$,
exchanging the roles of the indices $0$ and $j$ in the change of variables
in \eqref{eq:simplex_orig}
yields 
\begin{equation}
  \begin{split}
    \int_\sigma e^{\bd{a}\cdot \bd{x}}d\bd{x}&=\sum_{m=0}^N\frac{\abs{\det(\bd{p}_m-\bd{p}_0,\dots,\bd{p}_m-\bd{p}_{m-1},\bd{p}_m-\bd{p}_{m+1},\dots,\bd{p}_m-\bd{p}_N)}}{\prod_{j\neq m}\bd{a}\cdot (\bd{p}_m-\bd{p}_j)}e^{\bd{a}\cdot \bd{p}_m}\\
    &=N!\Vol(\sigma)\sum_{m=0}^N\frac{e^{\bd{a}\cdot \bd{p}_m}}{\prod_{j\neq m}\bd{a}\cdot (\bd{p}_m-\bd{p}_j)}.
  \end{split}
\end{equation}
Generally we have
\begin{equation}
   \sum_{m=0}^N\frac{e^{x_m}}{\prod_{j\neq m}(x_m-x_j)}=T(x_0,\ldots,x_N)
\end{equation}
by the Laplace expansion of the numerator of the right-hand side of \eqref{eq:det_form_simplex} with respect to the last row
and hence 
the result \eqref{eq:simplex}.
\end{proof}
Although
the following lemma is a direct consequence of the second expression of 
 \eqref{eq:simplex} with the definition of Schur polynomials and the Jacobi-Trudi formula
(see \cite{Mac95}),
we give a direct proof for convenience.
\begin{lemma}
\label{lm:simplex_expand}
Let $\bd{a}\in\mathbb{C}^N$, $\sigma$ and
  $\bd{p}_0,\dots,\bd{p}_N\in\mathbb{R}^N$ are the same as in Lemma \ref{lm:simplex}.
Then the Taylor expansion with respect to $\bd{a}$ is given by
\begin{multline}
    \int_\sigma e^{\bd{a}\cdot \bd{x}}d\bd{x}=\Vol(\sigma)
\Bigl(1
+\frac{1}{N+1}\sum_{0\leq i\leq N}\bd{a}\cdot \bd{p}_{i}
+\cdots\\
+
\frac{N!}{(N+k)!}
\sum_{\substack{k_0,\ldots,k_N\geq 0\\k_0+\cdots+k_N=k}}
(\bd{a} \cdot \bd{p}_0)^{k_0}
\cdots
(\bd{a} \cdot \bd{p}_N)^{k_N}
+\cdots
\Bigr).
\end{multline}
\end{lemma}
\begin{proof}
We recall Dirichlet's integral (see \cite[Chap.~XII 12.5]{WhitWats})
for nonnegative integers $k_j$ and a continuous function $g$, that is,
\begin{equation}
\label{eq:Dir_int}
  \int_{\Tilde{\tau}} y_0^{k_0}\cdots y_N^{k_N}g(y_0+\cdots+y_N)dy_0\cdots dy_N=
\frac{k_0!\cdots k_N!}{(N+k_0+\cdots+k_N)!}\int_0^1 g(t)t^{k_0+\cdots+k_N+N}dt,
\end{equation}
where $\Tilde{\tau}$ is the $(N+1)$-dimensional simplex 
with their vertices 
${}^t(0,\dots,\overset{j}{1},\dots,0)$,
 $0\leq j\leq N$ and the origin.
This formula is easily obtained by repeated application of the beta integral.

We calculate 
\begin{equation}
 f(\bd{a})=\frac{1}{k!}\int_\sigma(\bd{a}\cdot \bd{x})^k d\bd{x}.
\end{equation}
By multiplying
\begin{equation}
  1=(k+1)\int_0^1s^kds,
\end{equation}
and changing variables as $\bd{x}'=s\bd{x}$, we obtain
\begin{equation}
  f(\bd{a})=\frac{k+1}{k!}\int_{\Tilde{\sigma}} (\bd{a}\cdot \bd{x}')^ks^{-N}dsd\bd{x}',
\end{equation}
where 
$\Tilde{\sigma}=\bigcup_{0\leq s\leq1}\begin{pmatrix}s\\s\sigma\end{pmatrix}$
 is an $(N+1)$-dimensional simplex.
Again we change variables as
$\Tilde{P}\Tilde{\bd{y}}=\Tilde{\bd{x}}=\begin{pmatrix}s\\\bd{x}'\end{pmatrix}$.
Then
\begin{equation}
  \begin{split}
  f(\bd{a})
&=
\frac{k+1}{k!}\abs{\det\Tilde{P}}\int_{\Tilde{\tau}}
\Bigl(\sum_{j=0}^N (\bd{a}\cdot \bd{p}_j) y_j\Bigr)^k
\Bigl(\sum_{j=0}^N y_j\Bigr)^{-N}d\Tilde{\bd{y}}
\\
&=
\frac{k+1}{k!}\abs{\det\Tilde{P}}
\sum_{\substack{k_0,\ldots,k_N\geq 0\\k_0+\cdots+k_N=k}}
(\bd{a} \cdot \bd{p}_0)^{k_0}
\cdots
(\bd{a} \cdot \bd{p}_N)^{k_N}
\frac{k!}{k_0!\cdots k_N!}
\int_{\Tilde{\tau}}
y_0^{k_0}\cdots y_N^{k_N}
\Bigl(\sum_{j=0}^N y_j\Bigr)^{-N}d\Tilde{\bd{y}},
\end{split}
\end{equation}
where $\Tilde{\bd{y}}={}^t(y_0,\ldots,y_N)$.
Hence applying \eqref{eq:Dir_int}, we obtain
\begin{equation}
   \begin{split}
f(\bd{a})
&=
\frac{k+1}{k!}\abs{\det\Tilde{P}}
\sum_{\substack{k_0,\ldots,k_N\geq 0\\k_0+\cdots+k_N=k}}
(\bd{a} \cdot \bd{p}_0)^{k_0}
\cdots
(\bd{a} \cdot \bd{p}_N)^{k_N}
\frac{k!}{(N+k)!}\frac{1}{k+1}
\\
&=
\Vol(\sigma)
\frac{N!}{(N+k)!}
\sum_{\substack{k_0,\ldots,k_N\geq 0\\k_0+\cdots+k_N=k}}
(\bd{a} \cdot \bd{p}_0)^{k_0}
\cdots
(\bd{a} \cdot \bd{p}_N)^{k_N}.
\end{split}
\end{equation}
\end{proof}
It should be remarked that $\bd{a}\cdot \bd{p}_j$ for $0\leq j\leq N$ are not linearly independent.
Thus the coefficients with respect to them are not unique.
Lemma \ref{lm:simplex_expand} is a special and exact case of the following lemma.
\begin{lemma}
\label{lm:simplex_expand_b}
Let $\bd{a}\in\mathbb{C}^N$, $\sigma$ and
$\bd{p}_0,\dots,\bd{p}_N\in\mathbb{R}^N$ are the same as in Lemma \ref{lm:simplex}.
Then for $\bd{b}\in\mathbb{C}^N$,
the coefficients of total degree $k$ with respect to $\bd{a}$ of the Taylor expansion of
\begin{equation}
\label{eq:int_a_b} 
  \int_\sigma e^{(\bd{a}+\bd{b})\cdot \bd{x}}d\bd{x}
\end{equation}
are holomorphic functions in $\bd{b}$ of the form
\begin{equation}
\label{eq:coef_exp}
\Vol(\sigma)
\sum_{m=0}^N
\sum_{\substack{k_0,\ldots,k_N\geq 0\\k_0+\cdots+k_N=k}}
c_{m,k_0,\ldots,k_N}
  \frac{e^{\bd{b}\cdot \bd{p}_m}}{\prod_{j\neq m}\bigl(\bd{b}\cdot (\bd{p}_m-\bd{p}_j)\bigr)^{k_j+1}},
\end{equation}
where $c_{m,k_0,\ldots,k_N}\in\mathbb{Q}$.
\end{lemma}
\begin{proof}
We 
assume 
 $\abs{\bd{a}\cdot(\bd{p}_m-\bd{p}_j)}<\abs{\bd{b}\cdot(\bd{p}_m-\bd{p}_j)}$ for all $j\neq m$.
Then we have
\begin{multline}
\frac{e^{(\bd{a}+\bd{b})\cdot \bd{p}_m}}{\prod_{j\neq m}(\bd{a}+\bd{b})\cdot (\bd{p}_m-\bd{p}_j)}
\\
\begin{aligned}
&=
\frac{e^{\bd{b}\cdot \bd{p}_m}}{\prod_{j\neq m}\bd{b}\cdot (\bd{p}_m-\bd{p}_j)}
\frac{e^{\bd{a}\cdot \bd{p}_m}}{\prod_{j\neq m}\biggl(1+\dfrac{\bd{a}\cdot (\bd{p}_m-\bd{p}_j)}{\bd{b}\cdot (\bd{p}_m-\bd{p}_j)}\biggr)}
\\
&=
\frac{e^{\bd{b}\cdot \bd{p}_m}}{\prod_{j\neq m}\bd{b}\cdot (\bd{p}_m-\bd{p}_j)}
e^{\bd{a}\cdot \bd{p}_m}\prod_{j\neq m}\sum_{k_j=0}^\infty\biggl(-\dfrac{\bd{a}\cdot (\bd{p}_m-\bd{p}_j)}{\bd{b}\cdot (\bd{p}_m-\bd{p}_j)}\biggr)^{k_j}
\\
&=
\sum_{k=0}^\infty
\sum_{\substack{k_0,\ldots,k_N\geq 0\\k_0+\cdots+k_N=k}}
\frac{e^{\bd{b}\cdot \bd{p}_m}}{\prod_{j\neq m}\bigl(\bd{b}\cdot (\bd{p}_m-\bd{p}_j)\bigr)^{k_j+1}}
\frac{(\bd{a}\cdot \bd{p}_m)^{k_m}}{k_m!}
\prod_{j\neq m}
\bigl(\bd{a}\cdot (\bd{p}_j-\bd{p}_m)\bigr)^{k_j}.
\end{aligned}
\end{multline}
Applying this result to the second member of
\eqref{eq:simplex} with $\bd{a}$ replaced by $\bd{a}+\bd{b}$,
we see that the coefficients of total degree $k$ are of the form \eqref{eq:coef_exp}.
The holomorphy with respect to $\bd{a}$ and $\bd{b}$
follows from the original integral form \eqref{eq:int_a_b}.
Therefore \eqref{eq:coef_exp} is valid for the whole space with removable singularities.
\end{proof}

\section{General functional relations}\label{sec-3}
The purpose of this section is to give a very general formulation of functional relations.
For $f,g\in\Func(P)$ and
 $I,J\subset \{1,\ldots,r\}$,
we define
\begin{equation}
\label{eq:zeta}
  \zeta(f,g;J;\Delta)
=
  \sum_{\lambda\in P_{J++}\setminus H_g}
  \frac{f(\lambda)}{g(\lambda)},
\end{equation}
and
\begin{equation}
\label{eq:cf}
    S(f,g;I;\Delta)=\sum_{\lambda\in\iota^{*-1}(P_{I+})\setminus H_g}
    \frac{f(\lambda)}{g(\lambda)}.
\end{equation}
We assume that \eqref{eq:cf} is absolutely convergent for a fixed $I$.
By \eqref{eq:def_W_I} and Lemma \ref{lm:pp_wp}, we have $\id\in W^I$
and $P_+\subset\iota^{*-1}(P_{I+})$. Hence
\eqref{eq:zeta} is also absolutely convergent for all $J$
since $P_{J++}\subset P_+$ by \eqref{eq:Coxeter_comp_P}.

For $s\in\mathbb{C}$, $\Re s>1$ and $x,c\in\mathbb{R}$
let
\begin{equation}
  L_s(x,c)=-\frac{\Gamma(s+1)}{(2\pi\sqrt{-1})^s}\sum_{\substack{n\in\mathbb{Z}\\n+c\neq0}}\frac{e^{2\pi\sqrt{-1}(n+c)x}}{(n+c)^s}.
\end{equation}
Let $D_{I^c}$ be a finite subset of $(Q^\vee\setminus\{0\})\oplus\mathbb{R}\,\delta\subset\widehat{V}$.
Then any element of
$\gamma\in D_{I^c}$ can be written as $\gamma=\eta_\gamma+c_\gamma\delta$ 
($\eta_\gamma\in Q^\vee\setminus\{0\}$, $c_\gamma\in\mathbb{R}$).
We assume that $D_{I^c}$ contains $B_{I^c}=\{\gamma_i\}_{i\in I^c}$ where 
$\gamma_i=\eta_i+c_i\delta$ for $i\in I^c$ such that
$\{\eta_i\}_{i\in I^c}$ forms a basis of $Q^\vee_{I^c}$ and $c_i\in\mathbb{R}$.
Let $\{\mu_i\}_{i\in I^c}\subset P_{I^c}$ be a basis dual to $\{\eta_i\}_{i\in I^c}$.
\begin{theorem}
\label{thm:gen_func_rel}
Let $s_\gamma\in\mathbb{C}$ with $\Re s_\gamma>1$ 
for $\gamma\in D_{I^c}$ and let $\mathbf{y}\in V_{I^c}$.
We assume that
\begin{align}
\label{eq:def_f}
f(\lambda+\mu)&=f(\lambda)e^{2\pi\sqrt{-1}\langle \mathbf{y},\mu\rangle},\\
\label{eq:def_g}
g(\lambda+\mu)&=g^\sharp(\lambda)\prod_{\gamma\in D_{I^c}} \gamma(\lambda+\mu)^{s_\gamma},
\end{align}
for any $\lambda\in P_{I+}$ and any $\mu\in P_{I^c}$,
where $g^\sharp\in\Func(P_I)$.
(Hence $f$ depends on $\mathbf{y}$, and $g$ depends on $s_{\gamma}$'s.) 
Then 

(i) We have
\begin{equation}
\label{eq:gen_func_rel}
\begin{split}
S(f,g;I;\Delta)&=
      \sum_{w\in W^I}
      \sum_{J\subset\{1,\ldots,r\}}
      \frac{1}{N_{w,J}}
      \zeta(w^{-1}f,w^{-1}g;J;\Delta)
\\ &=\sum_{J\subset I}
      \zeta(f\cdot f^\sharp,g^\sharp;J;\Delta),
\end{split}
\end{equation}
where
\begin{equation}
  N_{w,J}=\abs{w W(\Delta_{J^c})\cap W^I}
\end{equation}
and $f^\sharp\in\Func(P_I)$ is defined by
\begin{equation}
\label{eq:gen_func_rel_f}
  \begin{split}
f^\sharp(\lambda)
&=(-1)^{\abs{D_{I^c}}}
  e^{-2\pi\sqrt{-1}\langle \mathbf{y},\nu\rangle}
  \biggl(\prod_{\gamma\in D_{I^c}}
  \frac{(2\pi\sqrt{-1})^{s_\gamma}}{\Gamma(s_\gamma+1)}\biggr)
\\
&\qquad\times
  \int_0^1\dots\int_0^1
  \exp\Bigl(
    -
    2\pi\sqrt{-1}
    \sum_{\gamma\in D_{I^c}\setminus B_{I^c}}\gamma(\lambda-\nu)x_\gamma
  \Bigr)
  \Bigl(\prod_{\gamma\in D_{I^c}\setminus B_{I^c}}
  L_{s_\gamma}(x_\gamma,c_\gamma)\Bigr)
\\
&\qquad\times
  \biggl(
  \prod_{i\in I^c}
  L_{s_{\gamma_i}}
  \Bigl(
  \langle \mathbf{y},\mu_i\rangle-\sum_{\gamma\in D_{I^c}\setminus B_{I^c}}x_\gamma\langle\eta_\gamma,\mu_i\rangle,c_i
  \Bigr)
  \biggr)
  \prod_{\gamma\in D_{I^c}\setminus B_{I^c}}
  dx_\gamma,
\end{split}
\end{equation}
and
\begin{equation}
\nu=\sum_{i\in I^c} c_i\mu_i\in \mathbb{R}\otimes P_{I^c}.
\end{equation}

(ii)
The second member of
\eqref{eq:gen_func_rel} consists of $2^r\bigl(W(\Delta):W(\Delta_I)\bigr)$ terms.

(iii)
If $H_{\Delta^\vee\setminus\Delta^\vee_I}\subset H_g$,
then $\zeta(w^{-1}f,w^{-1}g;J;\Delta)=0$ unless $N_{w,J}=1$.
\end{theorem}
\begin{proof} 
First we claim that
for $w,w'\in W^I$ and $\lambda\in w P_{J++}$,
we have
$\lambda\in w' P_{J++}$ 
if and only if
$w'\in w W(\Delta_{J^c})$.
In fact,
$\lambda\in w P_{J++}\cap w' P_{J++}$ 
implies
$w'^{-1}\lambda\in P_{J++}$ and
$(w^{-1}w')w'^{-1}\lambda\in P_{J++}$,
and hence $w^{-1}w'$ stabilizes $P_{J++}=P_+\cap C_J$.
Therefore $w^{-1}w'\in W(\Delta_{J^c})$.
The converse statement is shown by reversing the arguments
and we have the claim.
By using this claim, Lemma \ref{lm:pp_wp} and 
the decomposition \eqref{eq:Coxeter_comp_P},
we have
\begin{equation}
  \begin{split}
    S(f,g;I;\Delta)
    &=\sum_{\lambda\in \iota^{*-1}(P_{I+})\setminus H_g}
    \frac{f(\lambda)}{g(\lambda)}
\\
    &=\sum_{w\in W^I}\sum_{J\subset\{1,\ldots,r\}}
    \frac{1}{N_{w,J}}
    \sum_{\lambda\in wP_{J++}\setminus H_g}
    \frac{f(\lambda)}{g(\lambda)}.
  \end{split}
\end{equation}
Therefore
\begin{equation}
  \begin{split}
    S(f,g;I;\Delta)
    &=\sum_{w\in W^I}\sum_{J\subset\{1,\ldots,r\}}
    \frac{1}{N_{w,J}}
    \sum_{\lambda\in P_{J++}\setminus w^{-1}H_g}
    \frac{f(w\lambda)}{g(w\lambda)}
\\
    &=\sum_{w\in W^I}\sum_{J\subset\{1,\ldots,r\}}
    \frac{1}{N_{w,J}}
    \sum_{\lambda\in P_{J++}\setminus H_{w^{-1}g}}
    \frac{(w^{-1}f)(\lambda)}{(w^{-1}g)(\lambda)}
\\
    &=\sum_{w\in W^I}\sum_{J\subset\{1,\ldots,r\}}
    \frac{1}{N_{w,J}}
    \zeta(w^{-1}f,w^{-1}g;J;\Delta),
\end{split}
\end{equation}
where the last member consists of
$2^r\bigl(W(\Delta):W(\Delta_I)\bigr)$ terms since
the cardinality of the power set of $\{1,\ldots,r\}$ is $2^r$ and 
$\abs{W^I}=\bigl(W(\Delta):W(\Delta_I)\bigr)$ by Lemma \ref{lm:Z_mcr}, which implies the statement (ii).

Assume that $\gamma=\eta_\gamma+c_\gamma\delta\in(Q^\vee\setminus\{0\})\oplus\mathbb{R}\,\delta$ and $\lambda\in P$. 
Then $\gamma(\lambda)=\langle\eta_\gamma,\lambda\rangle+c_\gamma\in\mathbb{Z}+c_\gamma$ and
for $\gamma(\lambda)\neq 0$,
we have
\begin{equation}
  \begin{split}
\label{eq:1_over_gamma}
  \frac{1}{\gamma(\lambda)^{s_\gamma}}
&=\sum_{\substack{n_\gamma\in\mathbb{Z}\\n_\gamma+c_\gamma\neq0}}
\int_0^1\frac{e^{2\pi\sqrt{-1}(n_\gamma+c_\gamma-\gamma(\lambda))x_\gamma}}{(n_\gamma+c_\gamma)^{s_\gamma}}dx_\gamma
\\
&=\int_0^1\sum_{\substack{n_\gamma\in\mathbb{Z}\\n_\gamma+c_\gamma\neq0}}
\frac{e^{2\pi\sqrt{-1}(n_\gamma+c_\gamma-\gamma(\lambda))x_\gamma}}{(n_\gamma+c_\gamma)^{s_\gamma}}dx_\gamma
\\
&=-\frac{(2\pi\sqrt{-1})^{s_\gamma}}{\Gamma(s_\gamma+1)}\int_0^1L_{s_\gamma}(x_\gamma,c_\gamma)e^{-2\pi\sqrt{-1}\gamma(\lambda) x_\gamma}dx_\gamma,
\end{split} 
\end{equation}
where we have used the absolute convergence because of $\Re s_\gamma>1$.

By using \eqref{eq:def_g}, \eqref{eq:def_f} and Lemma \ref{lm:pp_wp}, we have
\begin{equation}
  S(f,g;I;\Delta)=
  \sum_{\lambda\in P_{I+}}
  \sum_{\substack{\mu\in P_{I^c}\\\lambda+\mu\not\in H_g}}
  \frac{f(\lambda)}
  {g^\sharp(\lambda)}
  e^{2\pi\sqrt{-1}\langle \mathbf{y},\mu\rangle}
\biggl(\prod_{\gamma\in D_{I^c}} \frac{1}{\gamma(\lambda+\mu)^{s_\gamma}}\biggr).
\end{equation}
Applying \eqref{eq:1_over_gamma} to the right-hand side of the above, we obtain
\begin{equation}
\label{eq:S_f_g_I_1}
  \begin{split}
  S(f,g;I;\Delta)
  &=
   \sum_{\lambda\in P_{I+}}
   \sum_{\substack{\mu\in P_{I^c}\\\lambda+\mu\not\in H_g}}
  \frac{f(\lambda)}{g^\sharp(\lambda)}
  e^{2\pi\sqrt{-1}\langle \mathbf{y},\mu\rangle}
  \biggl(
  \prod_{\gamma\in B_{I^c}}
  \frac{1}{\gamma(\lambda+\mu)^{s_\gamma}}
  \biggr)
\\
&\qquad\times
  \biggl(
  \prod_{\gamma\in D_{I^c}\setminus B_{I^c}}
  -\frac{(2\pi\sqrt{-1})^{s_\gamma}}{\Gamma(s_\gamma+1)}
  \int_0^1
  L_{s_\gamma}(x_\gamma,c_\gamma)
  e^{-2\pi\sqrt{-1}\gamma(\lambda+\mu) x_\gamma}dx_\gamma
  \biggr).
\end{split}
\end{equation}
Note that if $\gamma(\lambda)=0$,
then $c_\gamma\in\mathbb{Z}$ and
the last member of
\eqref{eq:1_over_gamma} vanishes.
Hence we may add the case $\gamma(\lambda+\mu)=0$ for $\gamma\in D_{I^c}\setminus B_{I^c}$ in the above.
Therefore by using $H_g=H_{g^\sharp}\cup H_{D_{I^c}\setminus B_{I^c}}\cup H_{B_{I^c}}$ and putting
$\mu=\sum_{i\in I^c}n_i\mu_i$ ($n_i\in\mathbb{Z}$), we have
\begin{multline}
\label{eq:S_f_g_I_2}
  S(f,g;I;\Delta)
\\
\begin{aligned}
  &=
  \biggl(\prod_{\gamma\in D_{I^c}\setminus B_{I^c}}
  -\frac{(2\pi\sqrt{-1})^{s_\gamma}}{\Gamma(s_\gamma+1)}\biggr)
   \sum_{\lambda\in P_{I+}\setminus H_{g^\sharp}}
  \frac{f(\lambda)}{g^\sharp(\lambda)}
\\&\qquad\times
  \sum_{\substack{n_i\in\mathbb{Z}\\n_i+c_i\neq0\\i\in I^c}}
  \int_0^1\dots\int_0^1
  \exp\Bigl(-2\pi\sqrt{-1}
    \sum_{\gamma\in D_{I^c}\setminus B_{I^c}}\gamma(\lambda)x_\gamma
  \Bigr)
  \Bigl(\prod_{\gamma\in D_{I^c}\setminus B_{I^c}}
  L_{s_\gamma}(x_\gamma,c_\gamma)\Bigr)
\\&\qquad\times
  \prod_{i\in I^c}
  \frac{
    \exp\Bigl(2\pi\sqrt{-1}(
      \langle \mathbf{y},\mu_i\rangle
      -
      \sum_{\gamma\in D_{I^c}\setminus B_{I^c}}x_\gamma\langle\eta_\gamma,\mu_i\rangle
)
n_i\Bigr)
  }{(n_i+c_i)^{s_{\gamma_i}}}
  \prod_{\gamma\in D_{I^c}\setminus B_{I^c}}  
  dx_\gamma
\\
  &=(-1)^{\abs{D_{I^c}}}
  e^{-2\pi\sqrt{-1}\langle \mathbf{y},\nu\rangle}
  \biggl(\prod_{\gamma\in D_{I^c}}
  \frac{(2\pi\sqrt{-1})^{s_\gamma}}{\Gamma(s_\gamma+1)}\biggr)
   \sum_{\lambda\in P_{I+}\setminus H_{g^\sharp}}
  \frac{f(\lambda)}{g^\sharp(\lambda)}
\\
&\qquad\times
  \int_0^1\dots\int_0^1
  \exp\Bigl(
    -
    2\pi\sqrt{-1}
    \sum_{\gamma\in D_{I^c}\setminus B_{I^c}}\gamma(\lambda-\nu)x_\gamma
  \Bigr)
  \Bigl(\prod_{\gamma\in D_{I^c}\setminus B_{I^c}}
  L_{s_\gamma}(x_\gamma,c_\gamma)\Bigr)
\\
&\qquad\times
  \biggl(
  \prod_{i\in I^c}
  L_{s_{\gamma_i}}
  \Bigl(
  \langle \mathbf{y},\mu_i\rangle-\sum_{\gamma\in D_{I^c}\setminus B_{I^c}}x_\gamma\langle\eta_\gamma,\mu_i\rangle,c_i
  \Bigr)
  \biggr)
  \prod_{\gamma\in D_{I^c}\setminus B_{I^c}}
  dx_\gamma,
\end{aligned}
\end{multline}
which is equal
to the third member of \eqref{eq:gen_func_rel}, because 
$P_{I+} = \bigcup_{J\subset I} P_{J++}$. 
Hence the statement (i).

As for the last claim (iii) of the theorem, we first note 
that $N_{w,J}>1$ if and only if
$wP_{J++}\subset H_{\Delta^\vee\setminus\Delta_I^\vee}$. 
This follows from Lemma \ref{lm:wp_uniq} and the definition of
$N_{w,J}$, with
noting the claim proved on the first line of the present proof.
Now assume 
$H_{\Delta^\vee\setminus\Delta_I^\vee}\subset H_g$.
Then, if $N_{w,J}>1$, we have $w P_{J++}\subset H_g$.
This implies $\zeta(w^{-1}f,w^{-1}g;J;\Delta)=0$, because 
the definition \eqref{eq:zeta} is an empty sum in this case.
\end{proof}
In the present paper we mainly discuss the case when
$D_{I^c}\subset Q^\vee$.    Nevertheless we give the above
generalized form of Theorem \ref{thm:gen_func_rel}, by which we can treat the
case of zeta-functions of Hurwitz-type.

For a real number $x$,
let $\{x\}$ denote its fractional part $x-[x]$.
If $s=k\geq2$ and $c$ are integers, then it is known (see \cite{Apos}) that
\begin{equation}
\label{eq:L_is_B}
  \begin{split}
  L_k(x,c)&=-\frac{k!}{(2\pi\sqrt{-1})^k}\sum_{n\in\mathbb{Z}\setminus\{0\}}\frac{e^{2\pi\sqrt{-1} nx}}{n^k}\\
&=B_k(\{x\}),
\end{split}
\end{equation}
where $B_k(x)$ is the $k$-th Bernoulli polynomial.
Thus if all $L_s(x,c)$ in the integrand on the right-hand side of 
\eqref{eq:gen_func_rel_f} are written in terms of Bernoulli polynomials,
then we have a chance to obtain an explicit form of the right-hand side of
\eqref{eq:gen_func_rel_f}.
We calculate the integral in question through a generating function instead of a direct calculation.
The result will be stated in Theorem \ref{thm:func_rel} below.
\begin{remark}
When $s_\gamma=1$,
the argument \eqref{eq:1_over_gamma}
is not valid, because
the series is conditionally convergent.
Hence on the right-hand side of \eqref{eq:S_f_g_I_2},
$s_{\gamma_i}=1$ for $i\in I^c$ is not allowed.
However for $n\in\mathbb{Z}\setminus\{0\}$, we have
\begin{equation}
  \frac{1}{n}=(-2\pi\sqrt{-1})\int_0^1 B_1(x)e^{-2\pi\sqrt{-1}nx}dx,
\end{equation}
where the right-hand side vanishes if $n=0$.
Thus the value $s_\gamma=1$ is also allowed 
for $\gamma\in(D_{I^c}\setminus B_{I^c})\cap\widehat{Q}^\vee$
on the right-hand side of \eqref{eq:S_f_g_I_1}
and hence in Theorem \ref{thm:gen_func_rel}.
\end{remark}

Let $M,N\in\mathbb{N}$, $\mathbf{k}=(k_l)_{1\leq l\leq M+N}\in\mathbb{N}_0^{M+N}$,
$\mathbf{y}=(y_i)_{1\leq i\leq M}\in\mathbb{R}^M$
and
$b_j\in\mathbb{C}$, $c_{ij}\in\mathbb{R}$
for $1\leq i\leq M$ and $1\leq j\leq N$.
Let
\begin{equation}
  \label{eq:def_P_gen}
P(\mathbf{k},\mathbf{y})=
  \int_0^1\dots\int_0^1
  \exp\Bigl(
    \sum_{j=1}^N b_jx_j
  \Bigr)
  \biggl(\prod_{j=1}^N
  B_{k_j}(x_j)\biggr)
  \biggl(
  \prod_{i=1}^M
  B_{k_{N+i}}
  \Bigl(\Bigl\{
  y_i-\sum_{j=1}^Nc_{ij}x_j\Bigr\}
  \Bigr)
  \biggr)
  \prod_{j=1}^N
  dx_j.
  \end{equation}
For $\mathbf{t}=(t_l)_{1\leq l\leq M+N}\in\mathbb{C}^{M+N}$,
we define a generating function of $P(\mathbf{k},\mathbf{y})$ by
\begin{equation}
\label{eq:def_F_gen} 
  F(\mathbf{t},\mathbf{y})=
  \sum_{\mathbf{k}\in\mathbb{N}_0^{M+N}}P(\mathbf{k},\mathbf{y})\prod_{j=1}^{M+N}\frac{t_j^{k_j}}{k_j!}.
\end{equation}
\begin{lemma}
\label{lm:gen_abs}
(i) The generating function $F(\mathbf{t},\mathbf{y})$ is 
absolutely convergent, uniformly on $\mathcal{D}_R\times\mathbb{R}^M$
where $\mathcal{D}_R=\{t\in\mathbb{C}~|~\abs{t}\leq R\}^{M+N}$ with $0<R<2\pi$.

(ii)  The function $F(\cdot,\mathbf{y})$ is 
analytically continued to a meromorphic function in $\mathbf{t}$ 
on the whole space $\mathbb{C}^{M+N}$ and 
we have
\begin{equation}
\label{eq:gen_func_h}
  \begin{split}
    F(\mathbf{t},\mathbf{y})
&=
  \biggl(\prod_{j=1}^{M+N}
  \frac{t_j}{e^{t_j}-1}\biggr)
  \int_0^1\dots\int_0^1
  \biggl(\prod_{j=1}^N
  \exp\bigl((b_j+t_j) x_j\bigr)\biggr)
\\&\qquad\times
  \biggl(
  \prod_{i=1}^M
  \exp\Bigl(t_{N+i}
\Bigl\{
y_i
-\sum_{j=1}^Nc_{ij}x_j
\Bigr\}
\Bigr)
  \biggr)
  \prod_{j=1}^N
  dx_j
\\
&=
  \biggl(\prod_{j=1}^{M+N}
  \frac{t_j}{e^{t_j}-1}\biggr)
  \sum_{m_1=c_1^-}^{c_1^+}
  \cdots
  \sum_{m_M=c_M^-}^{c_M^+}
  \exp\Bigl(\sum_{i=1}^Mt_{N+i}(\{y_i\}+m_i)\Bigr)
  \\&\qquad\times
  \int_{\mathcal{P}_{\mathbf{m},\mathbf{y}}}
  \exp\Bigl(
  \sum_{j=1}^N(a_j+b_j)x_j
  \Bigr)
  \prod_{j=1}^N
  dx_j,
  \end{split}
\end{equation}
where $c_i^+$ is the minimum integer 
satisfying $c_i^+\geq\sum_{\substack{1\leq j\leq N\\c_{ij}>0}} c_{ij}$ and
$c_i^-$ is the maximum integer 
satisfying $c_i^-\leq\sum_{\substack{1\leq j\leq N\\c_{ij}<0}} c_{ij}$ and
\begin{gather}
  a_j=t_j-\sum_{i=1}^M t_{N+i}c_{ij},\\
  \mathcal{P}_{\mathbf{m},\mathbf{y}}=
  \left\{\bd{x}=(x_j)_{1\leq j\leq N}~\left|~
  \begin{aligned}
  &0\leq x_j\leq1,
  \quad (1\leq j\leq N)
\\
  &\{y_i\}+m_i-1\leq
  \sum_{j=1}^Nc_{ij}x_j
  \leq \{y_i\}+m_i,
  \quad (1\leq i\leq M)
\end{aligned}
\right.
\right\},
\end{gather}
which is a convex polytope.
\end{lemma}
\begin{proof}
(i) 
Fix $R'\in\mathbb{R}$ such that $R<R'<2\pi$. Then by the Cauchy integral formula 
\begin{equation}
    \frac{B_k(x)}{k!}=\frac{1}{2\pi\sqrt{-1}}\int_{\abs{z}=R'}\frac{z e^{zx}}{e^z-1}\frac{dz}{z^{k+1}},
\end{equation}
we have for $0\leq x\leq1$ 
\begin{equation}
\Bigl\lvert\frac{B_k(x)}{k!}\Bigr\rvert
\leq\frac{1}{2\pi}\int_{\abs{z}=R'}\Bigl\lvert\frac{z e^{zx}}{e^z-1}\Bigr\rvert\frac{\abs{dz}}{R'^{k+1}}
\leq\frac{C_{R'}}{R'^k},
\end{equation}
where
\begin{equation}
  C_{R'}=\max\Bigl\{\Bigl\lvert\frac{ze^{zx}}{e^z-1}\Bigr\rvert~\Bigm|~\abs{z}=R',~0\leq x\leq1\Bigr\}.
\end{equation}
Therefore
\begin{equation}
  \begin{split}
  \biggl\lvert P(\mathbf{k},\mathbf{y})\prod_{j=1}^{M+N}\frac{t_j^{k_j}}{k_j!}
  \biggr\rvert
&\leq
  \biggl(\prod_{j=1}^{M+N}\abs{t_j}^{k_j}\biggr)
  \int_0^1\dots\int_0^1
  \biggl\lvert
  \exp\Bigl(
    \sum_{j=1}^N b_jx_j
  \Bigr)
  \biggl(
  \prod_{j=1}^N
  \frac{B_{k_j}(x_j)}{k_j!}
  \biggr)
\\&\qquad\times
  \biggl(
  \prod_{i=1}^M
  \frac{B_{k_{N+i}}
    \Bigl(
    \Bigl\{
    y_i
    -\sum_{j=1}^Nc_{ij}x_j
    \Bigr\}
    \Bigr)}
  {k_{N+i}!}
  \biggr)
  \biggr\rvert
  \prod_{j=1}^N
  dx_j
\\
&\leq
  C
  \biggl(\prod_{j=1}^{M+N}R^{k_j}\biggr)
  \int_0^1\dots\int_0^1
  C_{R'}^{M+N}
  \prod_{j=1}^{M+N}
  \frac{1}{R'^{k_j}}
  \prod_{j=1}^N
  dx_j
\\
&=
  C C_{R'}^{M+N}
  \prod_{j=1}^{M+N}
  \biggl(\frac{R}{R'}\biggr)^{k_j},
\end{split}
\end{equation}
where $C=\exp\Bigl(\sum_{j=1}^N \abs{\Re b_j}\Bigr)$.
Since
\begin{equation}
  \sum_{\mathbf{k}\in\mathbb{N}_0^{M+N}}
  C C_{R'}^{M+N}
  \prod_{j=1}^{M+N}
  \biggl(\frac{R}{R'}\biggr)^{k_j}
=
  C C_{R'}^{M+N}
  \prod_{j=1}^{M+N}\frac{1}{1-R/R'}<\infty,
\end{equation}
we have the uniform and absolute convergence of $F(\mathbf{t},\mathbf{y})$.

(ii) Noting the absolute convergence shown in (i), we obtain
\begin{equation}
  \begin{split}
    F(\mathbf{t},\mathbf{y})
&=
  \sum_{\mathbf{k}\in\mathbb{N}_0^{M+N}}P(\mathbf{k},\mathbf{y})\prod_{j=1}^{M+N}\frac{t_j^{k_j}}{k_j!}\\
&=\sum_{\mathbf{k}\in\mathbb{N}_0^{M+N}}
  \int_0^1\dots\int_0^1
\exp\Bigl(\sum_{j=1}^N b_jx_j\Bigr)
  \biggl(
  \prod_{j=1}^N
  B_{k_j}(x_j)
  \frac{t_j^{k_j}}{k_j!}
  \biggr)
\\&\qquad\times
  \biggl(
  \prod_{i=1}^M
B_{k_{N+i}}
    \Bigl(
    \Bigl\{
    y_i
    -\sum_{j=1}^Nc_{ij}x_j
    \Bigr\}
    \Bigr)
  \frac{t_{N+i}^{k_{N+i}}}{k_{N+i}!}
  \biggr)
  \prod_{j=1}^N
  dx_j
\\
&=
  \biggl(\prod_{j=1}^{M+N}
  \frac{t_j}{e^{t_j}-1}\biggr)
  \int_0^1\dots\int_0^1
\exp\Bigl(\sum_{j=1}^N b_jx_j\Bigr)
  \biggl(\prod_{j=1}^N
  \exp(t_j x_j)\biggr)
\\&\qquad\times
  \biggl(
  \prod_{i=1}^M
  \exp\Bigl(t_{N+i}
\Bigl\{
y_i
-\sum_{j=1}^Nc_{ij}x_j
\Bigr\}
\Bigr)
  \biggr)
  \prod_{j=1}^N
  dx_j
\\
&=
  \biggl(\prod_{j=1}^{M+N}
  \frac{t_j}{e^{t_j}-1}\biggr)
  \int_0^1\dots\int_0^1
  \exp\biggl(\sum_{j=1}^N
\Bigl(b_j+t_j-\sum_{i=1}^M t_{N+i}c_{ij}\Bigr) 
x_j\biggr)
\\&\qquad\times
  \exp\biggl(
\sum_{i=1}^M
t_{N+i}\Bigl(y_i-
\Bigl[
y_i
-\sum_{j=1}^Nc_{ij}x_j
\Bigr]
\Bigr)
  \biggr)
  \prod_{j=1}^N
  dx_j.
\end{split}
\end{equation}
Here the inequality
\begin{equation}
y_i-c_i^+
\leq
y_i-\sum_{j=1}^N c_{ij}x_j
\leq y_i-c_i^-
\end{equation}
implies
\begin{equation}
\{y_i\}+c_i^-
\leq
y_i-\Bigl[y_i-\sum_{j=1}^Nc_{ij}x_j\Bigr]
\leq
\{y_i\}+c_i^+,
\end{equation}
and for $m_i\in\mathbb{Z}$, the region of $\bd{x}$ satisfying
\begin{equation}
  y_i
-\Bigl[
y_i
-\sum_{j=1}^N c_{ij}x_j
\Bigr]
=\{y_i\}+m_i
\end{equation}
is given by
\begin{equation}
\{y_i\}+m_i-1<
\sum_{j=1}^N c_{ij}x_j
\leq\{y_i\}+m_i.
\end{equation}
Therefore
\begin{equation}
  \begin{split}
    F(\mathbf{t},\mathbf{y})
&=
  \biggl(\prod_{j=1}^{M+N}
  \frac{t_j}{e^{t_j}-1}\biggr)
  \sum_{m_1=c_1^-}^{c_1^+}
  \cdots
  \sum_{m_M=c_M^-}^{c_M^+}
  \exp\Bigl(\sum_{i=1}^Mt_{N+i}(\{y_i\}+m_i)\Bigr)
  \\&\qquad\times
  \int_{\mathcal{P}_{\mathbf{m},\mathbf{y}}}
  \exp\Bigl(
  \sum_{j=1}^N(a_j+b_j)x_j
  \Bigr)
  \prod_{j=1}^N
  dx_j,
\end{split}
\end{equation}
which  is 
 a meromorphic function in $\mathbf{t}$ 
on the whole space $\mathbb{C}^{M+N}$.
\end{proof}
\begin{lemma}
\label{lm:cont_PF}
The function
$P(\mathbf{k},\mathbf{y})$ is continuous with respect to $\mathbf{y}$ on $\mathbb{R}^M$.
The function
$F(\mathbf{t},\mathbf{y})$ is continuous on $\{t\in\mathbb{C}~|~\abs{t}<2\pi\}^{M+N}\times\mathbb{R}^M$
and
holomorphic in $\mathbf{t}$ for a fixed $\mathbf{y}\in\mathbb{R}^M$.
\end{lemma}
\begin{proof}
For $\mathbf{y}\in\mathbb{R}^M$ and
$\bd{x}=(x_j)_{1\leq j\leq N}\in\mathbb{R}^N$ let $h(\mathbf{y},\bd{x})$ 
be the integrand of \eqref{eq:def_P_gen}.
Let $\{\mathbf{y}_n\}_{n=1}^{\infty}$ be a sequence in $\mathbb{R}^M$ with $\lim_{n\to\infty}\mathbf{y}_n=\mathbf{y}_\infty$
and 
we set $h_n(\bd{x})=h(\mathbf{y}_n,\bd{x})$ and $h_\infty(\bd{x})=h(\mathbf{y}_\infty,\bd{x})$.
Then for $\bd{x}\in[0,1]^N$,
\begin{equation}
\label{eq:limit_P}
  \lim_{n\to\infty}h_n(\bd{x})=h_\infty(\bd{x})
\end{equation}
holds
if $\bd{x}$ satisfies
\begin{equation}
(\mathbf{y}_\infty)_i
-\sum_{j=1}^Nc_{ij}x_j
\not\in\mathbb{Z},
\end{equation}
for all $1\leq i\leq M$.
Hence \eqref{eq:limit_P} holds almost everywhere and 
we have
\begin{equation}
  \begin{split}
  \lim_{n\to\infty}P(\mathbf{k},\mathbf{y}_n)
&=\lim_{n\to\infty}\int_{[0,1]^N}h_n(\bd{x})\prod_{j=1}^Ndx_j\\
&=\int_{[0,1]^N}\lim_{n\to\infty}h_n(\bd{x})\prod_{j=1}^Ndx_j\\
&=\int_{[0,1]^N}h_\infty(\bd{x})\prod_{j=1}^Ndx_j\\
&=P(\mathbf{k},\mathbf{y}_\infty),
\end{split}
\end{equation}
where we have used the uniform boundedness $h_n(\bd{x})\leq C$ for some $C>0$ and for all $n\in\mathbb{N}$ and $\bd{x}\in[0,1]^N$.

Combining the continuity of $P(\mathbf{k},\mathbf{y})$, definition
\eqref{eq:def_F_gen} and Lemma \ref{lm:gen_abs},
we obtain the continuity and the holomorphy of $F(\mathbf{t},\mathbf{y})$.
\end{proof}

Now we return to the situation $\mathbf{y}\in V_{I^c}$.
Let
$y_i=\langle\mathbf{y},\mu_i\rangle$ for $i\in I^c$
and we identify $\mathbf{y}$ with $(y_i)_{i\in I^c}\in\mathbb{R}^{\abs{I^c}}$.
We set $\mathbb{Q}\,[\mathbf{y}]=\mathbb{Q}\,[(y_i)_{i\in I^c}]$,
$A(\mathbf{y})=\sum_{i\in I^c}\mathbb{Z}\,y_i+\mathbb{Z}$ 
and $\abs{\mathbf{k}}=\sum_{\gamma\in D_{I^c}}k_\gamma$. 
Let
\begin{equation}
  \overline{D}=
  \{\gamma\circ\tau(-\nu)~|~\gamma\in D_{I^c}\setminus B_{I^c}\},
\qquad  \nu=\sum_{i\in I^c} c_i\mu_i\in P_{I^c}.
\end{equation}
\begin{theorem}
\label{thm:func_rel}
Assume the same condition as in Theorem \ref{thm:gen_func_rel}.
Moreover we assume that
$D_{I^c}\subset \widehat{Q}^\vee$ and $s_\gamma=k_\gamma$ are integers for all $\gamma\in D_{I^c}$
such that $k_\gamma\geq2$ for $\gamma\in B_{I^c}$ and $k_\gamma\geq1$ otherwise.
Then $f^\sharp\in\Func(P_I)$ in \eqref{eq:gen_func_rel_f} is 
of the form
\begin{equation}
\label{eq:func_rel_main}
  f^\sharp(\lambda)=e^{-2\pi\sqrt{-1}\langle \mathbf{y},\nu\rangle}
\sum_{k=0}^{\abs{\mathbf{k}}}(\pi\sqrt{-1})^{\abs{\mathbf{k}}-(k+N)}
\sum_\eta
\frac{f_\eta^{(k)}(\lambda)}{g_\eta^{(k)}(\lambda)},
\end{equation}
where $\eta$ runs over a certain finite set of indices, 
$N=\abs{D_{I^c}\setminus B_{I^c}}$ and
\begin{equation}
  f_\eta^{(k)}\in\mathbb{Q}\,[\mathbf{y}]e^{\pi\sqrt{-1} \mathbb{Q}\,\otimes A(\mathbf{y})\overline{D}},\qquad
  g_\eta^{(k)}\in \bigl(A(\mathbf{y})\overline{D}\bigr)^{k+N}.
\end{equation}
\end{theorem}
\begin{proof}
From \eqref{eq:gen_func_rel_f}, we have
\begin{equation}
f^\sharp(\lambda)
=(-1)^{\abs{D_{I^c}}}
  e^{-2\pi\sqrt{-1}\langle \mathbf{y},\nu\rangle}
  \biggl(\prod_{\gamma\in D_{I^c}}
  \frac{(2\pi\sqrt{-1})^{k_\gamma}}{k_\gamma!}\biggr)P(\mathbf{k},\mathbf{y},\lambda),
\end{equation}
where for $\mathbf{k}=(k_\gamma)_{\gamma\in D_{I^c}}$,
\begin{equation}
\label{eq:h_l_k}
\begin{split}
P(\mathbf{k},\mathbf{y},\lambda)&=
  \int_0^1\dots\int_0^1
  \exp\Bigl(
    -2\pi\sqrt{-1}
    \sum_{\gamma\in D_{I^c}\setminus B_{I^c}}\gamma(\lambda-\nu)x_\gamma
  \Bigr)
  \Bigl(\prod_{\gamma\in D_{I^c}\setminus B_{I^c}}
  B_{k_\gamma}(x_\gamma)\Bigr)
\\
&\qquad\times
  \biggl(
  \prod_{i\in I^c}
  B_{k_{\gamma_i}}
  \Bigl(\Bigl\{
  y_i-\sum_{\gamma\in D_{I^c}\setminus B_{I^c}}x_\gamma\langle\eta_\gamma,\mu_i\rangle\Bigr\}
  \Bigr)
  \biggr)
  \prod_{\gamma\in D_{I^c}\setminus B_{I^c}}
  dx_\gamma.
\end{split}
\end{equation}
Hence
$P(\mathbf{k},\mathbf{y},\lambda)$ is of the form \eqref{eq:def_P_gen}.
Therefore, by applying Lemma \ref{lm:gen_abs},
we find that 
$P(\mathbf{k},\mathbf{y},\lambda)$ is obtained as the coefficient of the term
\begin{equation}
\label{eq:exp_var}
\prod_{\gamma\in D_{I^c}}t_\gamma^{k_\gamma}
\end{equation}
in the generating function
\begin{equation}
\label{eq:gen_func_h_l_k}
  \begin{split}
F(\mathbf{t},\mathbf{y},\lambda)
&=
  \sum_{\mathbf{k}\in\mathbb{N}_0^{\abs{D_{I^c}}}}P(\mathbf{k},\mathbf{y},\lambda)\prod_{\gamma\in D_{I^c}}\frac{t_\gamma^{k_\gamma}}{k_\gamma!}\\
&=
  \biggl(\prod_{\gamma\in D_{I^c}}
  \frac{t_\gamma}{e^{t_\gamma}-1}\biggr)
  \sum_{\substack{m_i=0\\i\in I^c}}^{c_i^+}
  \exp\Bigl(\sum_{i\in I^c}t_{\gamma_i}(\{y_i\}+m_i)\Bigr)
  \\&\qquad\times
  \int_{\mathcal{P}_{\mathbf{m},\mathbf{y}}}
  \exp\bigl((\bd{a}+\bd{b})\cdot\bd{x}
  \bigr)
  \prod_{\gamma\in D_{I^c}\setminus B_{I^c}}
  dx_\gamma,
  \end{split}
\end{equation}
where $\bd{a}=(a_\gamma)_{\gamma\in D_{I^c}\setminus B_{I^c}}\in\mathbb{R}^N,
\bd{b}=(b_\gamma)_{\gamma\in D_{I^c}\setminus B_{I^c}}
\in\mathbb{C}^{N}$ with
\begin{align}
\label{eq:def_a_g}
a_\gamma&=t_\gamma-\sum_{i\in I^c} t_{\gamma_i}\langle\eta_\gamma,\mu_i\rangle,\\
\label{eq:def_b_g}
b_\gamma&=-2\pi\sqrt{-1}\gamma(\lambda-\nu).
\end{align}
Since $\langle\eta_\gamma,\mu_i\rangle\in\mathbb{Z}$,
any vertex $\bd{p}_j$ of $\mathcal{P}_{\mathbf{m},\mathbf{y}}$ satisfies
\begin{equation}
\label{eq:p_rat}
(\bd{p}_j)_\gamma\in\sum_{i\in I^c}\mathbb{Q}\,y_i+\mathbb{Q}=\mathbb{Q}\otimes A(\mathbf{y}).
\end{equation}

The first two factors of the last member of \eqref{eq:gen_func_h_l_k}
is expanded as
\begin{equation}
\label{eq:coef_t_1}
  \biggl(\prod_{\gamma\in D_{I^c}}
  \frac{t_\gamma}{e^{t_\gamma}-1}\biggr)
  \biggl(
  \prod_{i\in I^c}
  \exp\bigl(t_{\gamma_i}(\{y_i\}+m_i)\bigr)
  \biggr)
  =
  \sum_{\mathbf{k}'\in\mathbb{N}_0^{\abs{D_{I^c}}}}P_{\mathbf{k}'}(\mathbf{y})
  \prod_{\gamma\in D_{I^c}}t_\gamma^{k_\gamma'},
\end{equation}
where
$P_{\mathbf{k}'}(\mathbf{y})\in\mathbb{Q}\,[\mathbf{y}]$ is
of total degree at most $\abs{\mathbf{k}'}$.

Next we calculate the contribution of $\mathbf{t}$ of
the integral part.
By a triangulation $\mathcal{P}_{\mathbf{m},\mathbf{y}}=\bigcup_{l=1}^{L(\mathbf{m})}\sigma_{l,\mathbf{m},\mathbf{y}}$ in Theorem \ref{thm:simp_div},
the integral on $\mathcal{P}_{\mathbf{m},\mathbf{y}}$ is reduced to those on $\sigma_{l,\mathbf{m},\mathbf{y}}$.
Since \eqref{eq:exp_var} is of total degree $\abs{\mathbf{k}}$, and 
$\bd{a}$ is of the same degree as $\bd{t}$ by
\eqref{eq:def_a_g},
the contribution of $\mathbf{t}$ comes from
terms of total degree $\kappa\leq\abs{\mathbf{k}}$ with respect to $\bd{a}$.
By Lemma \ref{lm:simplex_expand_b}, we see that
these terms are calculated as the special values of the functions $h(\mathbf{b}')$ at $\bd{b}$,
where $h(\bd{b}')$ is a holomorphic function on $\mathbb{C}^N$ of the form
\begin{equation}
\label{eq:coef_t_2}
\Vol(\sigma_{l,\mathbf{m},\mathbf{y}})
\sum_{q=0}^{N}
\sum_{\substack{\kappa_0,\ldots,\kappa_N\geq 0\\\kappa_0+\cdots+\kappa_N=\kappa}}
c_{q,\kappa_0,\ldots,\kappa_N}
  \frac{e^{\bd{b}'\cdot \bd{p}_{j_q}}}{\prod_{\substack{q'=0\\q'\neq q}}^N\bigl(\bd{b}'\cdot (\bd{p}_{j_q}-\bd{p}_{j_{q'}})\bigr)^{\kappa_{q'}+1}},
\end{equation}
and $\bd{p}_{j_q}$'s are the vertices of $\sigma_{l,\mathbf{m},\mathbf{y}}$ and
$\Vol(\sigma_{l,\mathbf{m},\mathbf{y}})\in\mathbb{Q}\,[\mathbf{y}]$ is of total degree at most $N$ due to \eqref{eq:p_rat}.
Since
\begin{equation}
\bd{b}\cdot \bd{p}_{j}=
-2\pi\sqrt{-1}
  \sum_{\gamma\in D_{I^c}\setminus B_{I^c}}
\gamma(\lambda-\nu)  
(\bd{p}_{j})_\gamma\in
\pi\sqrt{-1}\mathbb{Q}\otimes A(\mathbf{y})\overline{D},
\end{equation}
each term of \eqref{eq:coef_t_2} is an element of
\begin{equation}
\mathbb{Q}\,[\mathbf{y}]
\frac{e^{\pi\sqrt{-1}\mathbb{Q}\,\otimes A(\mathbf{y})\overline{D}}}{(\pi\sqrt{-1}\mathbb{Q}\otimes A(\mathbf{y})\overline{D})^{\kappa'+N}}
=
\frac{1}{(\pi\sqrt{-1})^{\kappa'+N}}
\mathbb{Q}\,[\mathbf{y}]
\frac{e^{\pi\sqrt{-1}\mathbb{Q}\,\otimes A(\mathbf{y})\overline{D}}}{\bigl(A(\mathbf{y})\overline{D}\bigr)^{\kappa'+N}},
\end{equation}
where $0\leq\kappa'\leq\kappa\leq\abs{\mathbf{k}}$.

Combining \eqref{eq:coef_t_1} and \eqref{eq:coef_t_2}
for all $\mathbf{m}$ and $l\in L(\mathbf{m})$
appearing in the sum, we see that
the coefficient of \eqref{eq:exp_var} is of the form
\eqref{eq:func_rel_main}.
\end{proof}

\begin{remark}
It may happen that the denominator of 
\eqref{eq:func_rel_main} vanishes. 
However the original form \eqref{eq:gen_func_rel_f}
implies that $f^\sharp$ is well-defined on $P_I$.
In fact, the values can be obtained
by use of analytic continuation in \eqref{eq:coef_t_2}.
\end{remark}

\begin{remark}
  It should be noted that 
  in Theorems \ref{thm:gen_func_rel} and \ref{thm:func_rel}
  we have treated $\mathbf{y}\in V_{I^c}$ as a fixed parameter.
  In general, as a function of $\mathbf{y}\in V_{I^c}$,
  \eqref{eq:func_rel_main} is not a real analytic function
  on the whole space $V_{I^c}$. 
  We study this fact in a special case in Section \ref{sec-6}.
\end{remark}

We conclude this section with the following proposition, whose
proof is a direct generalization of that of (2.1) in \cite{KM2}. 
\begin{proposition}
\label{prop:irr_decomp}
Let $f,g\in\Func(P)$ and $J\subset\{1,\ldots,r\}$.
Assume that $J=J_1\coprod J_2$ and that $f$ and $g$ are decomposed as
\begin{equation}
  f(\lambda_1+\lambda_2)=f_1(\lambda_1)f_2(\lambda_2),\qquad
  g(\lambda_1+\lambda_2)=g_1(\lambda_1)g_2(\lambda_2),
\end{equation}
for 
$f_1,g_1\in\Func(P_{J_1})$,
$f_2,g_2\in\Func(P_{J_2})$,
$\lambda_1\in P_{J_1}$,
$\lambda_2\in P_{J_2}$.
Then we have
\begin{equation}
  \zeta(f,g;J;\Delta)
=
  \zeta(f_1,g_1;J_1;\Delta)
  \zeta(f_2,g_2;J_2;\Delta).
\end{equation}
\end{proposition}

\section{Functional relations, value relations and generating functions}\label{sec-4}

Hereafter we deal with the special case
\begin{align}
\label{eq:def_sp_f}
f(\lambda)&=e^{2\pi\sqrt{-1}\langle \mathbf{y},\lambda\rangle},\\
\label{eq:def_sp_g}
g(\lambda)&=\prod_{\alpha\in\Delta_+}\langle\alpha^\vee,\lambda\rangle^{s_\alpha},
\end{align}
for $\mathbf{y}\in V$ and $\mathbf{s}=(s_\alpha)_{\alpha\in\prs}\in\mathbb{C}^{\abs{\Delta_+}}$,
where $\prs$ is the quotient of $\Delta$ obtained by identifying $\alpha$ and $-\alpha$.
We define an action of $\Aut$ by
\begin{equation}
\label{eq:A_act_on_s}
  (w\mathbf{s})_\alpha=s_{w^{-1}\alpha},
\end{equation}
and let
\begin{equation}
\label{eq:def_LZ}
  \zeta_r(\mathbf{s},\mathbf{y};\Delta)=\sum_{\lambda\in P_{++}}
e^{2\pi\sqrt{-1}\langle \mathbf{y},\lambda\rangle}
\prod_{\alpha\in\Delta_+}
\frac{1}{\langle\alpha^\vee,\lambda\rangle^{s_\alpha}}.
\end{equation}
Then we have
\begin{equation}
\label{eq:z_is_LZ}
\zeta(f,g;J;\Delta)
=\begin{cases}
\zeta_r(\mathbf{s},\mathbf{y};\Delta),\qquad &\text{if }J=\{1,\ldots,r\},\\
0,&\text{otherwise}
\end{cases}
\end{equation}
because for $J\neq\{1,\ldots,r\}$
we have $P_{J++}\subset H_{\Delta^\vee}=H_g$.
When $\mathbf{y}=0$, the function $\zeta_r(\mathbf{s};\Delta)=\zeta_r(\mathbf{s},0;\Delta)$
coincides with the zeta-function of the root system $\Delta$, defined
by (3.1) of \cite{KM2}. 
Therefore \eqref{eq:def_LZ} is the
Lerch-type generalization of zeta-functions of root systems. 
Also we have 
\begin{equation}
\label{eq:def_S}
  S(f,g;I;\Delta)
=\sum_{\lambda\in \iota^{*-1}(P_{I+})\setminus H_{\Delta^\vee}}
    e^{2\pi\sqrt{-1}\langle \mathbf{y},\lambda\rangle}
    \prod_{\alpha\in\Delta_+}
    \frac{1}{\langle\alpha^\vee,\lambda\rangle^{s_\alpha}},
\end{equation}
which we denote by $S(\mathbf{s},\mathbf{y};I;\Delta)$.
Note that $H_{\Delta^\vee}$ coincides with the
set of all walls of Weyl chambers.
Let 
\begin{equation}
  \mathcal{S}=\{\mathbf{s}=(s_\alpha)\in\mathbb{C}^{\abs{\Delta_+}}~|~
\Re s_\alpha>1\quad\text{for }\alpha\in\Delta_+\}.
\end{equation}
Note that $\mathcal{S}$ is an $\Aut$-invariant set.
\begin{lemma}
\label{lm:abs_conv}
$\zeta_r(\mathbf{s},\mathbf{y};\Delta)$ and
$S(\mathbf{s},\mathbf{y};I;\Delta)$
  are absolutely convergent, uniformly on 
$\mathcal{D}\times V$ where $\mathcal{D}$ is
any compact subset
of the set $\mathcal{S}$.
Hence they are continuous on $\mathcal{S}\times V$, and 
holomorphic in $\mathbf{s}$ for a fixed $\mathbf{y}\in V$.
\end{lemma}
\begin{proof}
Since for $\mathbf{s}\in\mathcal{D}$, $\alpha\in\Delta_+$ and $\lambda\in P\setminus H_{\Delta^\vee}$, 
  \begin{equation}
    \biggl\lvert\frac{1}{\langle\alpha^\vee,\lambda\rangle^{s_\alpha}}\biggr\rvert
\leq\frac{e^{\pi\abs{\Im s_\alpha}}}{\abs{\langle\alpha^\vee,\lambda\rangle}^{\Re s_\alpha}}
  \end{equation}
(the factor $e^{\pi\abs{\Im s_\alpha}}$ appears when $\langle\alpha^\vee,\lambda\rangle$ is negative),
we have
\begin{equation}
  \biggl\lvert\prod_{\alpha\in\Delta_+}
\frac{1}{\langle\alpha^\vee,\lambda\rangle^{s_\alpha}}\biggr\rvert\leq
h(\mathbf{s})
\prod_{\alpha\in\Delta_+}
\frac{1}{\abs{\langle\alpha^\vee,\lambda\rangle}^{\Re s_\alpha}}
\leq
h(\mathbf{s})
\prod_{\alpha\in\fs}
\frac{1}{\abs{\langle\alpha^\vee,\lambda\rangle}^{\Re s_\alpha}}
\leq
A
\prod_{\alpha\in\fs}
\frac{1}{\abs{\langle\alpha^\vee,\lambda\rangle}^{B}},
\end{equation}
where $h(\mathbf{s})=\prod_{\alpha\in\Delta_+}e^{\pi \abs{\Im s_\alpha}}$ and $A=\max_{\mathbf{s}\in\mathcal{D}}h(\mathbf{s})$,
 $B=\min_{\alpha\in\fs}(\min_{\mathbf{s}\in\mathcal{D}}\Re s_\alpha)>1$.
It follows that
\begin{equation}
  \sum_{\lambda\in P\setminus H_{\Delta^\vee}}
  \prod_{\alpha\in\fs}
  \frac{1}{\abs{\langle\alpha^\vee,\lambda\rangle}^{B}}
  =
  \abs{W}
  \sum_{\lambda\in P_{++}}
  \prod_{\alpha\in\fs}
  \frac{1}{\abs{\langle\alpha^\vee,\lambda\rangle}^{B}}
  =
  \abs{W}(\zeta(B))^r<\infty,
\end{equation}
and hence the uniform and absolute convergence on $\mathcal{D}\times V$.
\end{proof}
\begin{remark}
Although
the statements in Lemma \ref{lm:abs_conv}
and in the rest of this paper
hold for larger regions than $\mathcal{S}$,
we work with $\mathcal{S}$ for simplicity.
For instance, the above proof of Lemma \ref{lm:abs_conv}
holds for the region
$\{\mathbf{s}=(s_\alpha)\in\mathbb{C}^{\abs{\Delta_+}}~|~
\Re s_\alpha>1\quad\text{for }\alpha\in\fs,\quad
\Re s_\alpha>0\quad\text{otherwise}\}$.
\end{remark}

First 
we apply Theorem \ref{thm:gen_func_rel} to the case $I\neq\emptyset$.
Then Theorem \ref{thm:gen_func_rel} implies the following theorem:
\begin{theorem}
  \label{thm:FR}
When $I\neq\emptyset$, 
for $\mathbf{s}\in\mathcal{S}$ and $\mathbf{y}\in V$,
we have
\begin{multline}
\label{eq:func_eq}
S(\mathbf{s},\mathbf{y};I;\Delta)
\\
\begin{aligned}
&=
      \sum_{w\in W^I}
    \Bigl(\prod_{\alpha\in\Delta_{w^{-1}}}(-1)^{-s_{\alpha}}\Bigr)
\zeta_r(w^{-1}\mathbf{s},w^{-1}\mathbf{y};\Delta)\\
  &=(-1)^{\abs{\Delta_+\setminus\Delta_{I+}}}
  \biggl(\prod_{\alpha\in\Delta_+\setminus\Delta_{I+}}
  \frac{(2\pi\sqrt{-1})^{s_\alpha}}{\Gamma(s_\alpha+1)}\biggr)
  \sum_{\lambda\in P_{I++}}
e^{2\pi\sqrt{-1}\langle \mathbf{y},\lambda\rangle}
\prod_{\alpha\in\Delta_{I+}}
\frac{1}{\langle\alpha^\vee,\lambda\rangle^{s_\alpha}}
\\
&\qquad\times
  \int_0^1\dots\int_0^1
  \exp\Bigl(
    -
    2\pi\sqrt{-1}
    \sum_{\alpha\in\Delta_+\setminus(\Delta_{I+}\cup\fs)}x_\alpha\langle\alpha^\vee,\lambda\rangle
  \Bigr)
  \Bigl(\prod_{\alpha\in\Delta_+\setminus(\Delta_{I+}\cup\fs)}
  L_{s_\alpha}(x_\alpha,0)\Bigr)
\\
&\qquad\times
  \biggl(
  \prod_{i\in I^c}
  L_{s_{\alpha_i}}
  \Bigl(
  \langle \mathbf{y},\lambda_i\rangle-\sum_{\alpha\in\Delta_+\setminus(\Delta_{I+}\cup\fs)}x_\alpha\langle\alpha^\vee,\lambda_i\rangle,0
  \Bigr)
  \biggr)
  \prod_{\alpha\in\Delta_+\setminus(\Delta_{I+}\cup\fs)}
  dx_\alpha.
\end{aligned}
\end{multline}
The second member consists of $\bigl(W(\Delta):W(\Delta_I)\bigr)$ terms.
\end{theorem}
\begin{proof}
It is easy to check \eqref{eq:def_f} and
\eqref{eq:def_g} for \eqref{eq:def_sp_f} and
\eqref{eq:def_sp_g},
with 
$D_{I^c}=\Delta_+\setminus\Delta_{I+}$, $B_{I^c}=\fs_{I^c}$ 
(hence $c_\gamma$=0 for all $\gamma\in D_{I^c}$), 
and $g^\sharp(\lambda)=\prod_{\alpha\in\Delta_{I+}}\langle\alpha^\vee,\lambda\rangle^{s_\alpha}$ 
for $\lambda\in P_{I+}$.
Since $P_{J++}\subset H_{g^\sharp}$ for all $J\subsetneq I$,
we have the second equality.

Next we check the first equality.
From \eqref{eq:gen_func_rel}, \eqref{eq:z_is_LZ} and Theorem \ref{thm:gen_func_rel}
 (iii), we have
\begin{equation}
   S(\mathbf{s},\mathbf{y};I;\Delta)=\sum_{w\in W^I} \zeta(w^{-1}f, w^{-1}g;\{1,\ldots,r\};\Delta).
 \end{equation}
 Further,
  \begin{equation}
    \begin{split}
  \zeta(w^{-1}f;w^{-1}g;\{1,\ldots,r\};\Delta)
  &=
\sum_{\lambda\in P_{++}}
e^{2\pi\sqrt{-1}\langle \mathbf{y},w\lambda\rangle}
\prod_{\alpha\in\Delta_+}
\frac{1}{\langle\alpha^\vee,w\lambda\rangle^{s_\alpha}}\\
  &=
\sum_{\lambda\in P_{++}}
e^{2\pi\sqrt{-1}\langle w^{-1}\mathbf{y},\lambda\rangle}
\prod_{\alpha\in\Delta_+}
\frac{1}{\langle w^{-1}\alpha^\vee,\lambda\rangle^{s_\alpha}}\\
&=
\sum_{\lambda\in P_{++}}
e^{2\pi\sqrt{-1}\langle w^{-1}\mathbf{y},\lambda\rangle}
\prod_{\alpha\in w^{-1}\Delta_+}
\frac{1}{\langle\alpha^{\vee}, \lambda\rangle^{s_{w\alpha}}},
\end{split}
\end{equation}
by rewriting $\alpha$ as $w\alpha$. 
When
\begin{equation}
  \alpha \in w^{-1}\Delta_+ \cap \Delta_-
  =-(\Delta_+ \cap w^{-1}\Delta_-)
  =-\Delta_w,
\end{equation}
we further replace $\alpha$ by $-\alpha$.   Then
we have
\begin{multline}
  \zeta(w^{-1}f;w^{-1}g;\{1,\ldots,r\};\Delta)\\
  \begin{aligned}
  &=
\Bigl(\prod_{\alpha\in\Delta_w}(-1)^{-s_{w\alpha}}\Bigr)
\sum_{\lambda\in P_{++}}
e^{2\pi\sqrt{-1}\langle w^{-1}\mathbf{y},\lambda\rangle}
\prod_{\alpha\in\Delta_+}
\frac{1}{\langle\alpha^\vee,\lambda\rangle^{s_{w\alpha}}}\\
&=\Bigl(\prod_{\alpha\in\Delta_{w^{-1}}}(-1)^{-s_{\alpha}}\Bigr)
\zeta_r(w^{-1}\mathbf{s},w^{-1}\mathbf{y};\Delta),
\end{aligned}
\end{multline}
where we have used the fact that $w\Delta_w=-\Delta_{w^{-1}}$.
Hence the
first equality follows. 

Lemma \ref{lm:Z_mcr} implies that
the second member of \eqref{eq:func_eq} 
consists of $\bigl(W(\Delta):W(\Delta_I)\bigr)$ terms. 
\end{proof}

Next
we deal with the case $I=\emptyset$.
Let $S(\mathbf{s},\mathbf{y};\Delta)=S(\mathbf{s},\mathbf{y};\emptyset;\Delta)$.
Then we have
 the following theorem by Theorem \ref{thm:gen_func_rel}.
\begin{theorem}
For $\mathbf{s}\in\mathcal{S}$ and $\mathbf{y}\in V$,
\label{thm:val_rel}
\begin{equation}
    \begin{split}
S(\mathbf{s},\mathbf{y};\Delta)&=
      \sum_{w\in W}
    \Bigl(\prod_{\alpha\in\Delta_{w^{-1}}}(-1)^{-s_{\alpha}}\Bigr)
\zeta_r(w^{-1}\mathbf{s},w^{-1}\mathbf{y};\Delta)\\
  &=(-1)^{\abs{\Delta_+}}
  \biggl(\prod_{\alpha\in\Delta_+}
  \frac{(2\pi\sqrt{-1})^{s_\alpha}}{\Gamma(s_\alpha+1)}\biggr)
  \int_0^1\dots\int_0^1
  \Bigl(\prod_{\alpha\in\Delta_+\setminus\fs}
  L_{s_\alpha}(x_\alpha,0)\Bigr)
\\
&\qquad\times
  \biggl(
  \prod_{i=1}^r
  L_{s_{\alpha_i}}
  \Bigl(
  \langle \mathbf{y},\lambda_i\rangle-\sum_{\alpha\in\Delta_+\setminus\fs}x_\alpha\langle\alpha^\vee,\lambda_i\rangle,0
  \Bigr)
  \biggr)
  \prod_{\alpha\in\Delta_+\setminus\fs}
  dx_\alpha.
\end{split}
\end{equation}
\end{theorem}
The above two theorems are general functional relations among
zeta-functions of root systems with exponential factors. In some cases
it is possible to deduce, from these theorems, more explicit
functional relations among zeta-functions.  (See Example
\ref{eg:funcrel}).  However in general it is not easy to deduce
explicit forms of functional relations from \eqref{eq:func_eq} by
direct calculations.  Therefore, in our forthcoming paper \cite{KM5},
we will consider some structural background of our technique more
deeply and will present much improved versions of Theorem
\ref{thm:gen_func} and Theorem \ref{thm:Bernoulli}.  In fact by using
these results, we will give explicit forms of other concrete examples
which we do not treat in this paper.  On the other hand, in
\cite{KMT-CJ} we will introduce another technique of deducing explicit
forms. This can be regarded as a certain refinement of the
``$u$-method'' developed in our previous papers.  By using this
technique, we give explicit functional relations among zeta-functions
associated with root systems of types $A_3$, $C_2(\simeq B_2)$, $B_3$ and $C_3$.


Now we study special values of $S(\mathbf{s},\mathbf{y};\Delta)$. 
We recall that by \eqref{eq:L_is_B}
\begin{equation}
\label{eq:L_is_B-2}
L_k(x,0)=B_k(\{x\})
\end{equation}
for a real number $x$.
Motivated by this observation, 
for $\mathbf{k}=(k_\alpha)_{\alpha\in\prs}\in \mathbb{N}_0^{\abs{\Delta_+}}$ and $\mathbf{y}\in V$
we define
\begin{equation}
\label{eq:def_P}
  \begin{split}
  P(\mathbf{k},\mathbf{y};\Delta)
&=
  \int_0^1\dots\int_0^1
  \Bigl(\prod_{\alpha\in\Delta_+\setminus\fs}
  B_{k_\alpha}(x_\alpha)\Bigr)
\\
&\qquad\times
  \biggl(
  \prod_{i=1}^r
  B_{k_{\alpha_i}}
  \Bigl(\Bigl\{
  \langle \mathbf{y},\lambda_i\rangle-\sum_{\alpha\in\Delta_+\setminus\fs}x_\alpha\langle\alpha^\vee,\lambda_i\rangle
  \Bigr\}\Bigr)
  \biggr)
  \prod_{\alpha\in\Delta_+\setminus\fs}
  dx_\alpha,
\end{split}
\end{equation}
so that
\begin{equation}
\label{eq:S_B}
S(\mathbf{k},\mathbf{y};\Delta)=
(-1)^{\abs{\Delta_+}}
  \biggl(\prod_{\alpha\in\Delta_+}
  \frac{(2\pi\sqrt{-1})^{k_\alpha}}{k_\alpha!}\biggr)
P(\mathbf{k},\mathbf{y};\Delta)
\end{equation}
for $\mathbf{k}\in\mathcal{S}\cap\mathbb{N}_0^{\abs{\Delta_+}}$.
This function $P(\mathbf{k},\mathbf{y};\Delta)$
may be regarded as a generalization of the Bernoulli periodic functions 
and $B_\mathbf{k}(\Delta)=P(\mathbf{k},0;\Delta)$
the Bernoulli numbers
(see \cite{Apos}).
We define generating functions of $P(\mathbf{k},\mathbf{y};\Delta)$
and $B_\mathbf{k}(\Delta)$
as
\begin{align}
\label{eq:def_F}
  F(\mathbf{t},\mathbf{y};\Delta)&=
  \sum_{\mathbf{k}\in \mathbb{N}_0^{\abs{\Delta_+}}}P(\mathbf{k},\mathbf{y};\Delta)
  \prod_{\alpha\in\Delta_+}
  \frac{t_\alpha^{k_\alpha}}{k_\alpha!},
\\
  F(\mathbf{t};\Delta)&=
  \sum_{\mathbf{k}\in \mathbb{N}_0^{\abs{\Delta_+}}}B_{\mathbf{k}}(\Delta)
  \prod_{\alpha\in\Delta_+}
  \frac{t_\alpha^{k_\alpha}}{k_\alpha!},
\end{align}
where $\mathbf{t}=(t_\alpha)_{\alpha\in\prs}$ with $\abs{t_\alpha}<2\pi$.
Assume $\Delta$ is irreducible and not of type $A_1$.
Then
by Lemma \ref{lm:cont_PF}, we see that
$P(\mathbf{k},\mathbf{y};\Delta)$ is continuous in $\mathbf{y}$ on $V$ and
$F(\mathbf{t},\mathbf{y};\Delta)$ is continuous on $\{t\in\mathbb{C}~|~\abs{t}<2\pi\}^{\abs{\Delta_+}}\times V$ and
holomorphic in $\mathbf{t}$ for a fixed $\mathbf{y}\in V$.
Further 
by Lemma \ref{lm:gen_abs} (ii)
we see that
for a fixed $\mathbf{y}\in V$,
$F(\mathbf{t},\mathbf{y};\Delta)$ is
analytically continued to a meromorphic function
in $\mathbf{t}$ 
on the whole space $\mathbb{C}^{\abs{\Delta_+}}$.

\begin{theorem}
\label{thm:gen_func}
We have
\begin{equation}
\label{eq:gen_func}
  \begin{split}
    F(\mathbf{t},\mathbf{y};\Delta)&=
  \biggl(\prod_{\alpha\in\Delta_+}
  \frac{t_\alpha}{e^{t_\alpha}-1}\biggr)
  \int_0^1\dots\int_0^1
  \Bigl(\prod_{\alpha\in\Delta_+\setminus \fs}
  \exp(t_\alpha x_\alpha)\Bigr)
\\&\qquad\times
  \biggl(
  \prod_{i=1}^r
  \exp\Bigl(t_{\alpha_i}
\Bigl\{
\langle \mathbf{y},\lambda_i\rangle
-\sum_{\alpha\in\Delta_+\setminus \fs}x_\alpha\langle\alpha^\vee,\lambda_i\rangle
\Bigr\}
\Bigr)
  \biggr)
  \prod_{\alpha\in\Delta_+\setminus \fs}  
  dx_\alpha
\\
&=
  \biggl(\prod_{\alpha\in\Delta_+}
  \frac{t_\alpha}{e^{t_\alpha}-1}\biggr)
  \sum_{m_1=0}^{2\langle\rho^\vee,\lambda_1\rangle-1}
  \cdots
  \sum_{m_r=0}^{2\langle\rho^\vee,\lambda_r\rangle-1}
  \exp\Bigl(\sum_{i=1}^rt_{\alpha_i}(\{\langle \mathbf{y},\lambda_i\rangle\}+m_i)\Bigr)
  \\&\qquad\times
  \int_{\mathcal{P}_{\mathbf{m},\mathbf{y}}}
  \exp\Bigl(
  \sum_{\alpha\in\Delta_+\setminus \fs}t^*_\alpha x_\alpha
  \Bigr)
  \prod_{\alpha\in\Delta_+\setminus \fs}  
  dx_\alpha,
  \end{split}
\end{equation}
where
$\rho^\vee=\frac{1}{2}\sum_{\alpha\in\Delta_+}\alpha^\vee$
is the positive half sum,
$t^*_\alpha=t_\alpha-\sum_{i=1}^rt_{\alpha_i}\langle\alpha^\vee,\lambda_i\rangle$,
$\mathbf{m}=(m_1,\ldots,m_r)$
and
\begin{multline}
\label{eq:D_m_y}
  \mathcal{P}_{\mathbf{m},\mathbf{y}}=\\
  \left\{\bd{x}=(x_\alpha)_{\alpha\in\Delta_+\setminus \fs}~\left|~
  \begin{aligned}
  &0\leq x_\alpha\leq1,
  \quad (\alpha\in\Delta_+\setminus \fs)
\\
  &\{\langle \mathbf{y},\lambda_i\rangle\}+m_i-1\leq
  \sum_{\alpha\in\Delta_+\setminus \fs}x_\alpha\langle\alpha^\vee,\lambda_i\rangle
  \leq \{\langle \mathbf{y},\lambda_i\rangle\}+m_i,
  \quad (1\leq i\leq r)
\end{aligned}
\right.
\right\}
\end{multline}
is a convex polytope.
In particular, we have
\begin{equation}
  \begin{split}
    F(\mathbf{t};\Delta)
&=
  \biggl(\prod_{\alpha\in\Delta_+}
  \frac{t_\alpha}{e^{t_\alpha}-1}\biggr)
  \sum_{m_1=1}^{2\langle\rho^\vee,\lambda_1\rangle-1}
  \cdots
  \sum_{m_r=1}^{2\langle\rho^\vee,\lambda_r\rangle-1}
  \exp\Bigl(\sum_{i=1}^rt_{\alpha_i}m_i\Bigr)
  \\&\qquad\times
  \int_{\mathcal{P}_{\mathbf{m}}}
  \exp\Bigl(
  \sum_{\alpha\in\Delta_+\setminus \fs}t^*_\alpha x_\alpha
  \Bigr)
  \prod_{\alpha\in\Delta_+\setminus \fs}  
  dx_\alpha,
  \end{split}
\end{equation}
where
\begin{multline}
\label{eq:D_m}
  \mathcal{P}_{\mathbf{m}}=
  \mathcal{P}_{\mathbf{m},0}=\\
  \left\{\bd{x}=(x_\alpha)_{\alpha\in\Delta_+\setminus \fs}~\left|~
  \begin{aligned}
  &0\leq x_\alpha\leq1,
  \quad (\alpha\in\Delta_+\setminus \fs)
\\
  &m_i-1\leq
  \sum_{\alpha\in\Delta_+\setminus \fs}x_\alpha\langle\alpha^\vee,\lambda_i\rangle
  \leq m_i,
  \quad (1\leq i\leq r)
\end{aligned}
\right.
\right\}.
\end{multline}
\end{theorem}
\begin{proof}
Applying Lemma \ref{lm:gen_abs} to the case
$N=\abs{\Delta_+\setminus\fs}$, $M=\abs{\fs}$,
$b_j=0$, $k_j=k_\alpha$, 
$y_i=\langle\mathbf{y},\lambda_i\rangle$ and $c_{ij}=\langle\alpha^\vee,\lambda_i\rangle$,
we obtain \eqref{eq:gen_func}.
\end{proof}
From the above theorem we can deduce the following formula.
In the case when all $k_{\alpha}$'s are the same, this formula gives
a refinement of Witten's formula \eqref{1-2}.    In other words, it gives
a multiple generalization of the classical formula \eqref{eq:r_zeta}.
In \cite{Wi}, Witten showed that the volume of certain moduli spaces
can be written in terms of special values of series \eqref{1-1}.
Moreover he remarked that the volume
is rational in the orientable
case, which implies \eqref{1-2}. Zagier \cite{Za} gives a brief sketch of
a more number-theoretic demonstration of \eqref{1-2}.
Szenes \cite{Sz98} provides an algorithm of the evaluations
by use of iterated residues.
In our method, the rational number $C_W(2k,\mathfrak{g})$
is expressed in terms of generalized Bernoulli numbers,
which can be calculated by use of the generating functions.

\begin{theorem}\label{thm:W-Z}
Assume that $\Delta$ is an irreducible root system.
Let $k_\alpha=k_{\norm{\alpha}}\in\mathbb{N}$ and $\mathbf{k}=(k_\alpha)_{\alpha\in\prs}$.
Then we have
\begin{equation}
\label{eq:vol_formula}
  \zeta_r(2\mathbf{k};\Delta)
=\frac{(-1)^{\abs{\Delta_+}}}{\abs{W}}
  \biggl(\prod_{\alpha\in\Delta_+}
  \frac{(2\pi\sqrt{-1})^{2k_\alpha}}{(2k_\alpha)!}\biggr)
B_{2\mathbf{k}}(\Delta)
\in\mathbb{Q}\,\pi^{2\sum_l k_l\abs{(\Delta_+)_l}},
\end{equation}
where $l$ runs over the lengths of roots and $(\Delta_+)_l=\{\alpha\in\Delta_+~|~\norm{\alpha}=l\}$.
\end{theorem}
\begin{proof}
Since
the vertices $\bd{p}_j$ of $\mathcal{P_\mathbf{m}}$ satisfy
$(\bd{p}_j)_\alpha\in\mathbb{Q}$, by Theorem \ref{thm:gen_func} and Lemma \ref{lm:simplex_expand}, we have
\begin{equation}
  B_{2\mathbf{k}}(\Delta)=P(2\mathbf{k},0;\Delta)
\in\mathbb{Q},
\end{equation}
and hence by \eqref{eq:S_B},
\begin{equation}
  S(2\mathbf{k},0;\Delta)=
(-1)^{\abs{\Delta_+}}
  \biggl(\prod_{\alpha\in\Delta_+}
  \frac{(2\pi\sqrt{-1})^{2k_\alpha}}{(2k_\alpha)!}\biggr)
B_{2\mathbf{k}}(\Delta)
\in\mathbb{Q}\,\pi^{2\sum_l k_l\abs{(\Delta_+)_l}}.
\end{equation}
On the other hand, by Theorem \ref{thm:val_rel} 
\begin{equation}
 S(2\mathbf{k},0;\Delta)=\abs{W}\zeta_r(2\mathbf{k};\Delta),
\end{equation}
since roots of the same length form a single orbit.
Therefore we have \eqref{eq:vol_formula}.
\end{proof}
\begin{remark}
\label{rem:irr_decomp}
The assumption of the irreducibility of $\Delta$ in Theorem \ref{thm:W-Z} is not essential.
Since a reducible root system is decomposed into a direct sum of some irreducible root systems,
this assumption can be removed by use of Proposition \ref{prop:irr_decomp}.
\end{remark}
\begin{remark}
It is also to be stressed that our formula covers the case
when some of the $k_{\alpha}$'s are not the same.
For example, let $X_r=C_2(\simeq B_2)$. Then we can take the positive 
roots as $\{ \alpha_1,\alpha_2,2\alpha_1+\alpha_2,\alpha_1+\alpha_2\}$
with $\norm{\alpha_1^\vee}=\norm{\alpha_1^\vee+2\alpha_2^\vee}$, $\norm{\alpha_2^\vee}=\norm{\alpha_1^\vee+\alpha_2^\vee}$.
We see that 
\begin{equation}
  \begin{split}
\zeta_2(2,4,4,2;C_2)& =\sum_{m=1}^\infty\sum_{n=1}^\infty\frac{1}{m^2 n^4 (m+n)^4 (m+2n)^2} \\
& =\frac{53}{6810804000}\pi^{12}.
\end{split}
\end{equation}
Explicit forms of generating functions can be calculated with the aid of
Theorem \ref{thm:simp_div} and
Lemma \ref{lm:simplex}.
We give some more explicit examples in Section \ref{sec:examples}.
\end{remark}
\section{Actions of extended Weyl groups}\label{sec-5}
In this section, we study the action of $\widehat{W}$ on 
$S(\mathbf{s},\mathbf{y};\Delta)$,
$F(\mathbf{t},\mathbf{y};\Delta)$ and
$P(\mathbf{k},\mathbf{y};\Delta)$.
First consider the action of $\Aut\subset\widehat{W}$.
Note that $P\setminus H_{\Delta^\vee}$ is an $\Aut$-invariant set, because $H_{\Delta^\vee}$ is
$\Aut$-invariant.
An action of $\Aut$ is naturally induced on any function $f$ in $\mathbf{s}$ and $\mathbf{y}$ as follows:
For $w\in \Aut$, 
\begin{equation}
(wf)(\mathbf{s},\mathbf{y})=
  f(w^{-1}\mathbf{s},w^{-1}\mathbf{y}).
\end{equation}  
\begin{theorem}
\label{thm:sym_S}
For $\mathbf{s}\in\mathcal{S}$ and $\mathbf{y}\in V$,
and
for $w\in \Aut$, we have 
\begin{equation}
\label{eq:A_act_on_S}
(wS)(\mathbf{s},\mathbf{y};\Delta)
  =\Bigl(\prod_{\alpha\in\Delta_{w^{-1}}}(-1)^{-s_\alpha}\Bigr) S(\mathbf{s},\mathbf{y};\Delta),
  \end{equation}  
if $s_\alpha\in\mathbb{Z}$ for $\alpha\in\Delta_{w^{-1}}$.
\end{theorem}
\begin{proof}
From \eqref{eq:def_S}, we have
\begin{equation}
  (wS)(\mathbf{s},\mathbf{y};\Delta)
  =\sum_{\lambda\in P\setminus H_{\Delta^\vee}}
  e^{2\pi\sqrt{-1}\langle w^{-1}\mathbf{y},\lambda\rangle}
  \prod_{\alpha\in\Delta_+}
  \frac{1}{\langle\alpha^\vee,\lambda\rangle^{s_{w\alpha}}}.
\end{equation}
Rewriting $\lambda$ as $w^{-1}\lambda$ and noting
that $P\setminus H_{\Delta^\vee}$ is $\Aut$-invariant, we have
\begin{equation}
    \begin{split}
(wS)(\mathbf{s},\mathbf{y};\Delta)
&=\sum_{\lambda\in P\setminus H_{\Delta^\vee}}
    e^{2\pi\sqrt{-1}\langle \mathbf{y},\lambda\rangle}
    \prod_{\alpha\in\Delta_+}
    \frac{1}{\langle w\alpha^\vee,\lambda\rangle^{s_{w\alpha}}}\\
&=\sum_{\lambda\in P\setminus H_{\Delta^\vee}}
    e^{2\pi\sqrt{-1}\langle \mathbf{y},\lambda\rangle}
    \prod_{w^{-1}\alpha\in\Delta_+}
    \frac{1}{\langle\alpha^\vee,\lambda\rangle^{s_\alpha}}\\
&=\Bigl(\prod_{\alpha\in\Delta_{w^{-1}}}(-1)^{-s_\alpha}\Bigr) S(\mathbf{s},\mathbf{y};\Delta).
    \end{split}
  \end{equation}
\end{proof}
Thus we have
\begin{theorem}
\label{thm:S_vanish}
For $\mathbf{s}\in\mathcal{S}$ and $\mathbf{y}\in V$, we have
$S(\mathbf{s},\mathbf{y};\Delta)=0$ if
there exists an element $w\in \Aut_{\mathbf{s}}\cap \Aut_{\mathbf{y}}$ such that
$s_\alpha\in\mathbb{Z}$ for $\alpha\in\Delta_{w^{-1}}$ and
\begin{equation}
\label{eq:sum_not_even}
  \sum_{\alpha\in\Delta_{w^{-1}}}s_\alpha\not\in2\mathbb{Z},
\end{equation}
where $\Aut_{\mathbf{s}}$ and $\Aut_{\mathbf{y}}$ are the stabilizers
of $\mathbf{s}$ and $\mathbf{y}$ respectively
by regarding $\mathbf{y}\in V/Q^\vee$.
\end{theorem}
\begin{proof}
  Assume \eqref{eq:sum_not_even}. Then by Theorem \ref{thm:sym_S} 
  \begin{equation}
    \biggl(1-\Bigl(\prod_{\alpha\in\Delta_{w^{-1}}}(-1)^{-s_{\alpha}}\Bigr)\biggr)S(\mathbf{s},\mathbf{y};\Delta)=0,
  \end{equation}
  which implies $S(\mathbf{s},\mathbf{y};\Delta)=0$.
\end{proof}

\begin{lemma}
\label{lm:A_act_on_t}
  The extended Weyl group $\Aut$ acts on $\mathbb{C}^{\abs{\Delta_+}}$ by
\begin{equation}
\label{eq:A_act_on_t}
  (w\mathbf{t})_\alpha:=
  \begin{cases}
    t_{w^{-1}\alpha},\qquad &\text{if }\alpha\in\Delta_+\setminus\Delta_{w^{-1}},\\
    -t_{w^{-1}\alpha},\qquad &\text{if }\alpha\in\Delta_{w^{-1}},
  \end{cases}
\end{equation}
where $\mathbf{t}=(t_\alpha)_{\alpha\in\prs}\in\mathbb{C}^{\abs{\Delta_+}}$ and the representative $\alpha$ runs over $\Delta_+$.
\end{lemma}
\begin{proof}
What we have to check is that the definition \eqref{eq:A_act_on_t}
indeed defines an action. Since
\begin{equation}
  \bigl(v(w\mathbf{t})\bigr)_\alpha=
  \begin{cases}
    (w\mathbf{t})_{v^{-1}\alpha},\qquad &\text{if }\alpha\in\Delta_+\setminus\Delta_{v^{-1}},\\
    -(w\mathbf{t})_{v^{-1}\alpha},\qquad &\text{if }\alpha\in\Delta_{v^{-1}},
  \end{cases}
\end{equation}
we have
\begin{equation}
\label{eq:t_vw_is_pos}
      (v(w\mathbf{t}))_\alpha
  =
    t_{(vw)^{-1}\alpha},
\end{equation}
if and only if either
\begin{enumerate}
\item $\alpha\in\Delta_+\setminus\Delta_{v^{-1}}$ and $v^{-1}\alpha\in\Delta_+\setminus\Delta_{w^{-1}}$
\end{enumerate}
or
\begin{enumerate}
\setcounter{enumi}{1}
\item $\alpha\in\Delta_{v^{-1}}$ and $-v^{-1}\alpha\in\Delta_{w^{-1}}$
\end{enumerate}
holds. Here, the minus sign in the second case is caused by the
fact that if $\alpha\in \Delta_{v^{-1}}$ then $v^{-1}\alpha\in v^{-1}\Delta_{v^{-1}}=-\Delta_v\subset\Delta_-$.
Therefore \eqref{eq:t_vw_is_pos} is valid if and only if $\alpha\in\Delta_+$ and
\begin{equation}
  \begin{split}
\alpha
&\in
\bigl(v(\Delta_+\setminus\Delta_{w^{-1}})\cap v\Delta_+\bigr)\cup
\bigl(v(-\Delta_{w^{-1}})\cap v\Delta_-\bigr)
\\
&=
(v\Delta_+\setminus vw\Delta_-)\cup(v\Delta_-\cap vw\Delta_+)
\\
&=
(v\Delta_+\cap vw\Delta_+)\cup(v\Delta_-\cap vw\Delta_+)
\\
&=vw\Delta_+.
\end{split}
\end{equation}
This condition is equivalent to $\alpha\in\Delta_+\setminus\Delta_{(vw)^{-1}}$.
This implies
$v(w\mathbf{t})=(vw)(\mathbf{t})$.
\end{proof}
Note that we defined two types of actions of $\Aut$ on $\mathbb{C}^{\abs{\Delta_+}}$,
that is, \eqref{eq:A_act_on_s} and \eqref{eq:A_act_on_t}.
The action \eqref{eq:A_act_on_t} is used only on the variable $\mathbf{t}$
and should not be confused with the action \eqref{eq:A_act_on_s}.

If $\Delta$ is of type $A_1$, then $F(\mathbf{t},\mathbf{y};A_1)=te^{t\{y\}}/(e^t-1)$ 
(see Example \ref{eg:exam_A1} below)
is 
an even or, in other words, $\Aut$-invariant
function except for $y\in\mathbb{Z}$.
In the multiple cases, $F(\mathbf{t},\mathbf{y};\Delta)$ is revealed to be really an $\Aut$-invariant function. 
To show it, we need some notation and facts.
Fix $1\leq m\leq r$.
Note that $\sigma_m\Delta_+=(\Delta_+\setminus\{\alpha_m\})\coprod\{-\alpha_m\}$.
Let
$\Delta_1=(\Delta_+\setminus\fs)\cap \sigma_m(\Delta_+\setminus\fs)$ 
and
$\fs_1=\fs\cap \sigma_m(\Delta_+\setminus\fs)$ so that $\sigma_m(\Delta_+\setminus\fs)=\Delta_1\coprod\fs_1$.
Let
$\Delta_2=(\Delta_+\setminus\fs)\cap \sigma_m\fs$ 
and
$\fs_2=\fs\cap \sigma_m\fs$.
Then we have $\Delta_+\setminus\fs=\Delta_1\coprod\Delta_2$
and $\fs=\fs_1\coprod\fs_2\coprod\{\alpha_m\}$.
Moreover we see that $\sigma_m$ fixes $\fs_2$ pointwise and $\fs_1=\sigma_m\Delta_2$.
\begin{lemma}
\label{lm:lambda}
  \begin{equation}
\sum_{\alpha_i\in\fs_1}
\lambda_i\langle\alpha_i^\vee,\alpha_m\rangle=\alpha_m-2\lambda_m.
\end{equation}
\end{lemma}
\begin{proof}
Note that 
$\alpha_i\in\fs_1$ if and only if $\langle\alpha_i^\vee,\alpha_m\rangle\neq 0$ and $\alpha_i\neq\alpha_m$.
Let $v$ be the left-hand side.
Then
we have
\begin{equation}
  \langle \alpha_k^\vee,v\rangle=
  \begin{cases}
    \langle\alpha_k^\vee,\alpha_m\rangle,\qquad&\text{if }\alpha_k\in\fs_1,\\
    0,\qquad&\text{if }\alpha_k\in\fs_2\cup\{\alpha_m\},
  \end{cases}
\end{equation}
which determines the right-hand side uniquely.
\end{proof}

An action of $\Aut$ is naturally induced on any function $f$ in $\mathbf{t}$ and $\mathbf{y}$ as follows:
For $w\in \Aut$, 
\begin{equation}
(wf)(\mathbf{t},\mathbf{y})=
  f(w^{-1}\mathbf{t},w^{-1}\mathbf{y}).
\end{equation}  
\begin{theorem}
\label{thm:gen_sym}
Assume that $\Delta$ is an irreducible root system.
If $\Delta$ is not of type $A_1$, then
\begin{equation}
\label{eq:inv_F}
(wF)(\mathbf{t},\mathbf{y};\Delta)=F(\mathbf{t},\mathbf{y};\Delta)
\end{equation}
for $\mathbf{t}\in\mathbb{C}^{\abs{\Delta_+}}$ and $\mathbf{y}\in V$, and
for $w\in \Aut$.
Hence for $\mathbf{k}\in\mathbb{N}_0^{\abs{\Delta_+}}$ and $\mathbf{y}\in V$,
\begin{equation}
\label{eq:sym_B} 
(wP)(\mathbf{k},\mathbf{y};\Delta)
=
\Bigl(\prod_{\alpha\in\Delta_{w^{-1}}}(-1)^{-k_\alpha}\Bigr)
P(\mathbf{k},\mathbf{y};\Delta).
\end{equation}
\end{theorem}
\begin{remark}
If $\mathbf{k}$ is
in the region $\mathcal{S}$
of absolute convergence with respect to $\mathbf{s}$, the relation \eqref{eq:S_B}
and Theorem \ref{thm:sym_S} immediately imply \eqref{eq:sym_B},
while $\mathbf{k}\not\in\mathcal{S}$, it should be proved independently.
\end{remark}
\begin{remark}
  The assumption of the irreducibility is not essential by the same reason 
as in Remark \ref{rem:irr_decomp}.
\end{remark}
\begin{proof}
It is sufficient to show \eqref{eq:inv_F} for the cases $w=\sigma_m\in W$ and $w=\omega\in\Omega$
because $\Aut$ is generated by simple reflections and 
the subgroup $\Omega$.
Applying the simple reflection $\sigma_m$ to the second member of \eqref{eq:gen_func}, we have
\begin{multline}
  (\sigma_mF)(\mathbf{t},\mathbf{y};\Delta)
\\
  \begin{aligned}
&=
  \biggl(\prod_{\alpha\in\Delta_+}
  \frac{t_\alpha}{e^{t_\alpha}-1}\biggr)
  \int_0^1\dots\int_0^1
  \Bigl(\prod_{\alpha\in\Delta_+\setminus \fs}
  \exp(t_{\sigma_m\alpha} x_\alpha)\Bigr)
\\&\qquad\times
  \exp\Bigl(t_{\alpha_m}
\Bigl(1-
\Bigl\{
\langle \sigma_m\mathbf{y},\lambda_m\rangle
-\sum_{\alpha\in\Delta_+\setminus \fs}x_\alpha\langle\alpha^\vee,\lambda_m\rangle
\Bigr\}
\Bigr)
\Bigr)
\\&\qquad\times
  \biggl(
  \prod_{\substack{i=1\\i\neq m}}^r
  \exp\Bigl(t_{\sigma_m\alpha_i}
\Bigl\{
\langle \sigma_m\mathbf{y},\lambda_i\rangle
-\sum_{\alpha\in\Delta_+\setminus \fs}x_\alpha\langle\alpha^\vee,\lambda_i\rangle
\Bigr\}
\Bigr)
  \biggr)
  \prod_{\alpha\in\Delta_+\setminus \fs}  
  dx_\alpha,
\end{aligned}
\end{multline}
where we have used the fact that by the action of $\sigma_m$,
the factor
$\prod_{\alpha\in\Delta_+}t_\alpha/(e^{t_\alpha}-1)$ is sent to
\begin{equation}
  \begin{split}
  \frac{-t_{\sigma_m\alpha_m}}{e^{-t_{\sigma_m\alpha_m}}-1}
\prod_{\alpha\in\Delta_+\setminus\{\alpha_m\}}
  \frac{t_{\sigma_m\alpha}}{e^{t_{\sigma_m\alpha}}-1}
&=
  \frac{t_{\sigma_m\alpha_m}e^{t_{\sigma_m\alpha_m}}}{e^{t_{\sigma_m\alpha_m}}-1}
\prod_{\alpha\in\Delta_+\setminus\{\alpha_m\}}
  \frac{t_{\sigma_m\alpha}}{e^{t_{\sigma_m\alpha}}-1}\\
&=
  e^{t_{\alpha_m}}
\prod_{\alpha\in\Delta_+}
  \frac{t_\alpha}{e^{t_\alpha}-1}.
\end{split}
\end{equation}
Therefore, rewriting $x_\alpha$ as $x_{\sigma_m\alpha}$, we have
\begin{multline}
\label{eq:r_mF} 
  (\sigma_mF)(\mathbf{t},\mathbf{y};\Delta)
\\
  \begin{aligned}
&=
  \biggl(\prod_{\alpha\in\Delta_+}
  \frac{t_\alpha}{e^{t_\alpha}-1}\biggr)
  \int_0^1\dots\int_0^1
  \Bigl(\prod_{\alpha\in\Delta_+\setminus \fs}
  \exp(t_{\sigma_m\alpha} x_{\sigma_m\alpha})\Bigr)
\\&\qquad\times
  \exp\Bigl(t_{\alpha_m}
\Bigl(1-
\Bigl\{
\langle \mathbf{y},\sigma_m\lambda_m\rangle
-\sum_{\alpha\in\Delta_+\setminus \fs}x_{\sigma_m\alpha}\langle \sigma_m\alpha^\vee,\sigma_m\lambda_m\rangle
\Bigr\}
\Bigr)
\Bigr)
\\&\qquad\times
  \biggl(
  \prod_{\substack{i=1\\i\neq m}}^r
  \exp\Bigl(t_{\sigma_m\alpha_i}
\Bigl\{
\langle \mathbf{y},\sigma_m\lambda_i\rangle
-\sum_{\alpha\in\Delta_+\setminus \fs}x_{\sigma_m\alpha}\langle \sigma_m\alpha^\vee,\sigma_m\lambda_i\rangle
\Bigr\}
\Bigr)
  \biggr)
  \prod_{\alpha\in\Delta_+\setminus \fs}  
  dx_{\sigma_m\alpha}
\\
&=
  \biggl(\prod_{\alpha\in\Delta_+}
  \frac{t_\alpha}{e^{t_\alpha}-1}\biggr)
  \int_0^1\dots\int_0^1
  \Bigl(\prod_{\alpha\in \sigma_m(\Delta_+\setminus \fs)}
  \exp(t_\alpha x_\alpha)\Bigr)
\\&\qquad\times
  \exp\Bigl(t_{\alpha_m}
\Bigl(1-
\Bigl\{
\langle \mathbf{y},\lambda_m-\alpha_m\rangle
-\sum_{\alpha\in \sigma_m(\Delta_+\setminus \fs)}x_\alpha\langle\alpha^\vee,\lambda_m-\alpha_m\rangle
\Bigr\}
\Bigr)
\Bigr)
\\&\qquad\times
  \biggl(
  \prod_{\alpha_i\in\fs\setminus\{\alpha_m\}}
  \exp\Bigl(t_{\sigma_m\alpha_i}
\Bigl\{
\langle \mathbf{y},\lambda_i\rangle
-\sum_{\alpha\in \sigma_m(\Delta_+\setminus \fs)}x_\alpha\langle\alpha^\vee,\lambda_i\rangle
\Bigr\}
\Bigr)
  \biggr)
  \prod_{\alpha\in \sigma_m(\Delta_+\setminus \fs)}  
  dx_\alpha
\\
&=
  \biggl(\prod_{\alpha\in\Delta_+}
  \frac{t_\alpha}{e^{t_\alpha}-1}\biggr)
  \int_0^1\dots\int_0^1
  \Bigl(\prod_{\alpha\in\Delta_1}
  \exp(t_\alpha \{x_\alpha\})\Bigr)
  \Bigl(\prod_{\alpha_i\in\fs_1}
  \exp(t_{\alpha_i}\{x_{\alpha_i}\})\Bigr)
\\
&\qquad\times
  \exp\Bigl(t_{\alpha_m}
  \Bigl(1-
  \Bigl\{
  \langle \mathbf{y},\lambda_m-\alpha_m\rangle-\sum_{\alpha\in\Delta_1\cup\fs_1}x_\alpha
  \langle\alpha^\vee,\lambda_m-\alpha_m\rangle
  \Bigr\}\Bigr)\Bigr)
\\
&\qquad\times
  \biggl(
  \prod_{\alpha_j\in\fs_2}
  \exp\Bigl(t_{\alpha_j}
  \Bigl\{
  \langle \mathbf{y},\lambda_j\rangle-\sum_{\alpha\in\Delta_1\cup\fs_1}x_\alpha
  \langle\alpha^\vee,\lambda_j\rangle
  \Bigr\}\Bigr)
  \biggr)
\\
&\qquad\times
  \biggl(
  \prod_{\alpha_i\in\fs_1}
  \exp\Bigl(t_{\sigma_m\alpha_i}
  \Bigl\{
  \langle \mathbf{y},\lambda_i\rangle-\sum_{\alpha\in\Delta_1\cup\fs_1}x_\alpha
  \langle\alpha^\vee,\lambda_i\rangle
  \Bigr\}\Bigr)
  \biggr)
  \prod_{\alpha\in\Delta_1\cup\fs_1}
  dx_\alpha.
\end{aligned}
\end{multline}
Here we change variables from $\bd{x}=(x_\alpha)_{\alpha\in\Delta_1\cup\fs_1}$ to $\bd{z}=(z_\alpha)_{\alpha\in\Delta_1\cup\Delta_2}$ as 
\begin{equation}
\label{eq:chg_var_z_x}
  z_\alpha=
  \begin{cases}
    x_\alpha,\qquad&\text{if }\alpha\in\Delta_1,\\
  \langle \mathbf{y},\lambda_i\rangle-\sum_{\beta\in\Delta_1\cup\fs_1}x_\beta
  \langle\beta^\vee,\lambda_i\rangle,\qquad&\text{if }\alpha=\sigma_m\alpha_i\in\Delta_2,
  \end{cases}
\end{equation}
so that the Jacobian matrix is calculated as
\begin{equation}
\frac{\partial\bd{z}}{\partial\bd{x}}=
\begin{pmatrix}
  I_{\abs{\Delta_1}} & 0\\
  * & -I_{\abs{\Delta_2}}
\end{pmatrix},
\end{equation}
where $I_p$ is the $p\times p$ identity matrix,
since
\begin{equation}
  \begin{split}
    z_{\sigma_m\alpha_i}
    &=
    \langle \mathbf{y},\lambda_i\rangle-\sum_{\alpha\in\Delta_1\cup\fs_1}x_\alpha
    \langle\alpha^\vee,\lambda_i\rangle\\
    &=
    \langle \mathbf{y},\lambda_i\rangle-\sum_{\alpha\in\Delta_1}x_\alpha
    \langle\alpha^\vee,\lambda_i\rangle-x_{\alpha_i}.
  \end{split}
\end{equation}
Thus we have
$\abs{\det\partial\bd{x}/\partial\bd{z}}=1$.
For $\alpha=\sigma_m\alpha_k\in\Delta_2$ and $\alpha_i\in\fs_1$,
we have
\begin{equation}
  \langle\alpha^\vee,\lambda_i\rangle
=\langle\sigma_m\alpha^\vee_k,\lambda_i\rangle
=\langle\alpha^\vee_k,\sigma_m\lambda_i\rangle
=\langle\alpha^\vee_k,\lambda_i\rangle
=\delta_{ki},
\end{equation}
and hence
\begin{equation}
\label{eq:chg_var_x_z}
  \begin{split}
    x_{\alpha_i}
    &=
    \langle \mathbf{y},\lambda_i\rangle-\sum_{\alpha\in\Delta_1}z_\alpha
    \langle\alpha^\vee,\lambda_i\rangle-z_{\sigma_m\alpha_i}\\
    &=
    \langle \mathbf{y},\lambda_i\rangle-\sum_{\alpha\in\Delta_+\setminus\fs}z_\alpha
    \langle\alpha^\vee,\lambda_i\rangle.
  \end{split}
\end{equation}
For the fourth factor of the last integral in \eqref{eq:r_mF},
we have
\begin{equation}
  \langle\alpha^\vee,\lambda_j\rangle
=\langle\sigma_m\alpha^\vee_k,\lambda_j\rangle
=\langle\alpha^\vee_k,\sigma_m\lambda_j\rangle
=\langle\alpha^\vee_k,\lambda_j\rangle
=0,
\end{equation}
for $\alpha=\sigma_m\alpha_k\in\Delta_2$ 
and $\alpha_j\in\fs_2$,
and hence
\begin{equation}
\label{eq:chg_var_four}
  \begin{split}
  \langle \mathbf{y},\lambda_j\rangle-\sum_{\alpha\in\Delta_1\cup\fs_1}x_\alpha
\langle\alpha^\vee,\lambda_j\rangle
&=    
  \langle \mathbf{y},\lambda_j\rangle-\sum_{\alpha\in\Delta_1}z_\alpha
\langle\alpha^\vee,\lambda_j\rangle\\
&=    
  \langle \mathbf{y},\lambda_j\rangle-\sum_{\alpha\in\Delta_+\setminus\fs}z_\alpha
\langle\alpha^\vee,\lambda_j\rangle.
  \end{split}
\end{equation}
For the third factor, we have
\begin{multline}
  \langle \mathbf{y},\lambda_m-\alpha_m\rangle-\sum_{\alpha\in\Delta_1\cup\fs_1}x_\alpha
\langle\alpha^\vee,\lambda_m-\alpha_m\rangle\\
  \begin{aligned}
&=    
  \langle \mathbf{y},\lambda_m-\alpha_m\rangle-
\sum_{\alpha\in\Delta_1}x_\alpha
\langle\alpha^\vee,\lambda_m-\alpha_m\rangle
-\sum_{\alpha_i\in\fs_1}x_{\alpha_i}
\langle\alpha_i^\vee,\lambda_m-\alpha_m\rangle\\
&=    
  \langle \mathbf{y},\lambda_m-\alpha_m\rangle
-\sum_{\alpha\in\Delta_1}x_\alpha
\langle\alpha^\vee,\lambda_m-\alpha_m\rangle
+\sum_{\alpha_i\in\fs_1}x_{\alpha_i}
\langle\alpha_i^\vee,\alpha_m\rangle\\
&=    
  \langle \mathbf{y},\lambda_m-\alpha_m\rangle-
\sum_{\alpha\in\Delta_1}z_\alpha
\langle\alpha^\vee,\lambda_m-\alpha_m\rangle\\
&\qquad
+\sum_{\alpha_i\in\fs_1}
\biggl(
\langle \mathbf{y},\lambda_i\rangle-\sum_{\alpha\in\Delta_+\setminus\fs}z_\alpha
\langle\alpha^\vee,\lambda_i\rangle
\biggr)
\langle\alpha_i^\vee,\alpha_m\rangle
\end{aligned}
\end{multline}
by using \eqref{eq:chg_var_x_z}. Hence we have
\begin{multline}
\label{eq:chg_var_three}
  \langle \mathbf{y},\lambda_m-\alpha_m\rangle-\sum_{\alpha\in\Delta_1\cup\fs_1}x_\alpha
\langle\alpha^\vee,\lambda_m-\alpha_m\rangle\\
  \begin{aligned}
&=    
  \langle \mathbf{y},\lambda_m-\alpha_m\rangle-
\sum_{\alpha\in\Delta_1}z_\alpha
\langle\alpha^\vee,\lambda_m-\alpha_m\rangle\\
&\qquad
+\sum_{\alpha_i\in\fs_1}
\langle \mathbf{y},\lambda_i\rangle\langle\alpha_i^\vee,\alpha_m\rangle
-\sum_{\alpha\in\Delta_+\setminus\fs}
\sum_{\alpha_i\in\fs_1}
z_\alpha
\langle\alpha^\vee,\lambda_i\rangle\langle\alpha_i^\vee,\alpha_m\rangle
\\
&=
-\langle \mathbf{y},\lambda_m\rangle
+\sum_{\alpha\in\Delta_+\setminus\fs}z_\alpha
\langle\alpha^\vee,\lambda_m\rangle,
\end{aligned}
\end{multline}
where in the last line we have used Lemma \ref{lm:lambda} and the fact that
for $\alpha=\sigma_m\alpha_k\in\Delta_2$ 
we have
\begin{equation}
  \langle\alpha^\vee,\lambda_m-\alpha_m\rangle
=\langle\sigma_m\alpha^\vee_k,\sigma_m\lambda_m\rangle
=\langle\alpha^\vee_k,\lambda_m\rangle
=0.
\end{equation}
Since all the factors of the integrand 
of the right-hand side of \eqref{eq:r_mF} 
are periodic functions with its period $1$,
we integrate the interval $[0,1]$ with respect to $z_\alpha$. Therefore
using
\eqref{eq:chg_var_z_x}, 
\eqref{eq:chg_var_x_z}, \eqref{eq:chg_var_four} and
\eqref{eq:chg_var_three}
we have
\begin{multline}
  (\sigma_mF)(\mathbf{t},\mathbf{y};\Delta)
\\
  \begin{aligned}
&=
  \biggl(\prod_{\alpha\in\Delta_+}
  \frac{t_\alpha}{e^{t_\alpha}-1}\biggr)
  \int_0^1\dots\int_0^1
  \Bigl(\prod_{\alpha\in\Delta_1}
  \exp(t_\alpha \{z_\alpha\})\Bigr)
\\
&\qquad\times
  \biggl(\prod_{\alpha_i\in\fs_1}
  \exp\Bigl(t_{\alpha_i}
  \Bigl\{
  \langle \mathbf{y},\lambda_i\rangle-\sum_{\alpha\in\Delta_+\setminus\fs}z_\alpha
  \langle\alpha^\vee,\lambda_i\rangle
  \Bigr\}\Bigr)
  \biggr)
\\
&\qquad\times
  \exp\Bigl(t_{\alpha_m}\Bigl(1-
  \Bigl\{
  -\langle \mathbf{y},\lambda_m\rangle
  +\sum_{\alpha\in\Delta_+\setminus\fs}z_\alpha
  \langle\alpha^\vee,\lambda_m\rangle
  \Bigr\}\Bigr)\Bigr)
\\
&\qquad\times
  \biggl(
  \prod_{\alpha_j\in\fs_2}
  \exp\Bigl(t_{\alpha_j}
  \Bigl\{
  \langle \mathbf{y},\lambda_j\rangle-\sum_{\alpha\in\Delta_+\setminus\fs}z_\alpha
  \langle\alpha^\vee,\lambda_j\rangle
  \Bigr\}\Bigr)
  \biggr)
\\
&\qquad\times
  \Bigl(
  \prod_{\alpha_i\in\fs_1}
  \exp(t_{\sigma_m\alpha_i}\{z_{\sigma_m\alpha_i}\})
  \Bigr)
  \prod_{\alpha\in\Delta_+\setminus\fs}
  dz_\alpha
\\
&=
  \biggl(\prod_{\alpha\in\Delta_+}
  \frac{t_\alpha}{e^{t_\alpha}-1}\biggr)
  \int_0^1\dots\int_0^1
  \Bigl(\prod_{\alpha\in\Delta_+\setminus\fs}
  \exp(t_\alpha z_\alpha)\Bigr)
\\
&\qquad\times
  \biggl(\prod_{\alpha_i\in\fs\setminus\{\alpha_m\}}
  \exp\Bigl(t_{\alpha_i}
  \Bigl\{
  \langle \mathbf{y},\lambda_i\rangle-\sum_{\alpha\in\Delta_+\setminus\fs}z_\alpha
  \langle\alpha^\vee,\lambda_i\rangle
  \Bigr\}\Bigr)
  \biggr)
\\
&\qquad\times
  \exp\Bigl(t_{\alpha_m}\Bigl(1-
  \Bigl\{
  -\langle \mathbf{y},\lambda_m\rangle
  +\sum_{\alpha\in\Delta_+\setminus\fs}z_\alpha
  \langle\alpha^\vee,\lambda_m\rangle
  \Bigr\}\Bigr)\Bigr)
  \prod_{\alpha\in\Delta_+\setminus\fs}
  dz_\alpha
\\
&=
F(\mathbf{t},\mathbf{y};\Delta),
\end{aligned}
\end{multline}
where in the last line,
we have used the fact that for any $\alpha_m\in\fs$, there exists a root $\alpha\in\Delta_+\setminus\fs$ such that 
$\langle\alpha^\vee,\lambda_m\rangle\neq 0$ and thus
 $1-\{-x\}=\{x\}$ for $x\in\mathbb{R}\setminus\mathbb{Z}$  
implies that the integrand coincides with that of $F(\mathbf{t},\mathbf{y};\Delta)$
almost everywhere.

Lastly we check the invariance with respect to $\omega\in\Omega$.
Since $\omega\in\Omega$ permutes $\fs$ and leaves $\Delta_+$ and hence $\Delta_+\setminus\fs$ invariant,
we have
\begin{multline}
  (\omega F)(\mathbf{t},\mathbf{y};\Delta)
\\
  \begin{aligned}
&=
  \biggl(\prod_{\alpha\in\Delta_+}
  \frac{t_{\omega\alpha}}{e^{t_{\omega\alpha}}-1}\biggr)
  \int_0^1\dots\int_0^1
  \Bigl(\prod_{\alpha\in\Delta_+\setminus \fs}
  \exp(t_{\omega\alpha} x_{\omega\alpha})\Bigr)
\\&\qquad\times
  \biggl(
  \prod_{i=1}^r
  \exp\Bigl(t_{\omega\alpha_i}
\Bigl\{
\langle \omega^{-1}\mathbf{y},\lambda_i\rangle
-\sum_{\alpha\in\Delta_+\setminus \fs}x_{\omega\alpha}\langle\omega\alpha^\vee,\omega\lambda_i\rangle
\Bigr\}
\Bigr)
  \biggr)
  \prod_{\alpha\in\Delta_+\setminus \fs}  
  dx_{\omega\alpha}
\\
&=
  \biggl(\prod_{\alpha\in\Delta_+}
  \frac{t_\alpha}{e^{t_\alpha}-1}\biggr)
  \int_0^1\dots\int_0^1
  \Bigl(\prod_{\alpha\in\Delta_+\setminus \fs}
  \exp(t_\alpha x_\alpha)\Bigr)
\\&\qquad\times
  \biggl(
  \prod_{i=1}^r
  \exp\Bigl(t_{\omega\alpha_i}
\Bigl\{
\langle \mathbf{y},\omega\lambda_i\rangle
-\sum_{\alpha\in\Delta_+\setminus \fs}x_\alpha\langle\alpha^\vee,\omega\lambda_i\rangle
\Bigr\}
\Bigr)
  \biggr)
  \prod_{\alpha\in\Delta_+\setminus \fs}  
  dx_\alpha
\\
&=F(\mathbf{t},\mathbf{y};\Delta).
\end{aligned}
\end{multline}
\end{proof}

It is possible to extend the action of $\Aut$ to that of $\widehat{W}$ as follows:
For $q\in Q^\vee$,
\begin{equation}
  \begin{split}
  (\tau(q)\mathbf{s})_\alpha&=s_{\alpha}, \\
  (\tau(q)\mathbf{t})_\alpha&=t_{\alpha}, \\
  \tau(q)\mathbf{y}&=\mathbf{y}+q.
\end{split}
\end{equation}
We can observe the periodicity of $S$, $F$ and $P$ 
with respect to $\mathbf{y}$ from 
\eqref{eq:def_S}, \eqref{eq:def_P} and the first 
line of
\eqref{eq:gen_func}.
From this periodicity, we have
\begin{theorem}
\label{thm:inv_qv}
The action of $\Aut$ is extended to that of $\widehat{W}$ and is given by
\begin{equation}
  \begin{split}
    (\tau(q)S)(\mathbf{s},\mathbf{y};\Delta)&=S(\mathbf{s},\mathbf{y};\Delta),\\
    (\tau(q)F)(\mathbf{t},\mathbf{y};\Delta)&=F(\mathbf{t},\mathbf{y};\Delta),\\
    (\tau(q)P)(\mathbf{k},\mathbf{y};\Delta)&=P(\mathbf{k},\mathbf{y};\Delta),
  \end{split}
\end{equation}  
for $q\in Q^\vee$.
\end{theorem}

\begin{remark}
Some statements 
related with $S(\mathbf{s},\mathbf{y};\Delta)$ or $\zeta_r(\mathbf{s},\mathbf{y};\Delta)$
in Sections \ref{sec-4} and \ref{sec-5}
hold on any regions in $\mathbf{s}$ to which these functions are analytically continued.
In particular, in the case $\mathbf{y}=0$,
as we noticed at the beginning of Section \ref{sec-4},
the function $\zeta_r(\mathbf{s};\Delta)=\zeta_r(\mathbf{s},0;\Delta)$
coincides with the zeta-function defined in \cite{KM2},
and
its analytic continuation
is given in \cite[Theorem 6.1]{KM2}.
\end{remark}

\section{Generalization of Bernoulli polynomials}\label{sec-6}
In the previous sections, we have investigated
$P(\mathbf{k},\mathbf{y};\Delta)$ as a continuous function in $\mathbf{y}$.
In fact, this function is not real analytic in $\mathbf{y}$ in general.
However they are piecewise real analytic, 
and each piece is actually a polynomial in $\mathbf{y}$.   
In this section we will prove this fact, and will discuss
basic properties of those polynomials.

Let $\mathfrak{D}=\{\mathbf{y}\in V~|~0\leq\langle \mathbf{y},\lambda_i\rangle\leq1,\quad (1\leq i\leq r)\}$ be a period-parallelotope
 of $F(\mathbf{t},\cdot;\Delta)$ with its interior.
Let $\mathscr{R}$ be the set of all linearly independent subsets
$\mathbf{R}=\{\beta_1,\ldots,\beta_{r-1}\}\subset\Delta$, 
$\mathfrak{H}_{\mathbf{R}^\vee}=\bigoplus_{i=1}^{r-1}\mathbb{R}\,\beta_i^\vee$ the hyperplane passing through $\mathbf{R}^\vee\cup\{0\}$ and
\begin{equation}
\label{eq:def_H_R}
\mathfrak{H}_{\mathscr{R}}:=
    \bigcup_{\substack{\mathbf{R}\in\mathscr{R}\\q\in Q^\vee}}(\mathfrak{H}_{\mathbf{R}^\vee}+q).
\end{equation}
\begin{lemma}
\label{lm:loc_fin}
We have
\begin{equation}
\label{eq:H_R}
\mathfrak{H}_{\mathscr{R}}=
    \bigcup_{w\in W}w\Biggl(\bigcup_{j=1}^r
(\mathfrak{H}_{\fs^\vee\setminus\{\alpha_j^\vee\}}+\mathbb{Z}\,\alpha_j^\vee)\Biggr).
\end{equation}
The set $\{\mathfrak{H}_{\mathbf{R}^\vee}+q~|~\mathbf{R}\in\mathscr{R},q\in Q^\vee\}$ is locally finite, i.e., for any $\mathbf{y}\in V$, 
there exists a neighborhood $U(\mathbf{y})$ such that $U(\mathbf{y})$ intersects
finitely many of these hyperplanes.
\end{lemma}
\begin{proof}
  Fix $\mathbf{R}\in\mathscr{R}$. 
Then
$\Tilde{\Delta}^\vee=\Delta^\vee\cap\mathfrak{H}_{\mathbf{R}^\vee}$
is a coroot system so that $\mathbf{R}^\vee\subset\Tilde{\Delta}^\vee$.
Let $\mu$ be a nonzero vector normal to $\mathfrak{H}_{\mathbf{R}^\vee}$.
Then there exists an element $w\in W$ such that $w^{-1}\mu\in C$.
Put $w^{-1}\mu=\sum_{j=1}^r c_j\lambda_j$ with $c_j\geq 0$.
Then $\alpha^\vee=\sum_{j=1}^r a_j\alpha_j^\vee\in\Delta^\vee$ orthogonal to $w^{-1}\mu$
should satisfy $\sum_{j=1}^r a_jc_j=0$.
Since $a_j$ are all nonpositive or nonnegative,
we have $a_j=0$ for $j$ such that $c_j\neq 0$.  
Hence $c_j=0$ except for only one $j$,
because $w^{-1}\Tilde{\Delta}^\vee\subset\Delta^\vee$ is orthogonal to $w^{-1}\mu$ 
with codimension $1$.
That is,
 $w^{-1}\mu=c\lambda_j$
for some $c>0$.
Therefore $w(\fs^\vee\setminus\{\alpha_j^\vee\})$ 
is a fundamental system of $\Tilde{\Delta}^\vee$
and 
$\mathfrak{H}_{\mathbf{R}^\vee}
=
\mathfrak{H}_{w(\fs^\vee\setminus\{\alpha_j^\vee\})}
=
w\mathfrak{H}_{\fs^\vee\setminus\{\alpha_j^\vee\}}
$.
Moreover $Q^\vee=w Q^\vee=\bigoplus_{i=1}^r\mathbb{Z}\,w\alpha_i^\vee$, which implies
\begin{equation}
  \mathfrak{H}_{\mathbf{R}^\vee}+Q^\vee=
  w\mathfrak{H}_{\fs^\vee\setminus\{\alpha_j^\vee\}}+\mathbb{Z}\,w\alpha_j^\vee,
\end{equation}
since
$\bigoplus_{i=1,i\neq j}^r\mathbb{Z}\,\alpha_i^\vee
\subset\mathfrak{H}_{\fs^\vee\setminus\{\alpha_j^\vee\}}$.
This shows that 
$\mathfrak{H}_{\mathscr{R}}$ is contained in
the right-hand side of \eqref{eq:H_R}.
The opposite inclusion is clear.
The local finiteness follows from
the expression \eqref{eq:H_R} and
$\abs{W}<\infty$.
\end{proof}

Due to the local finiteness shown in
Lemma \ref{lm:loc_fin} and $\partial\mathfrak{D}\subset\mathfrak{H}_{\mathscr{R}}$,
we denote by $\mathfrak{D}^{(\nu)}$
each open connected component of $\mathfrak{D}\setminus\mathfrak{H}_{\mathscr{R}}$
so that
\begin{equation}
\label{eq:D_div}
  \mathfrak{D}\setminus\mathfrak{H}_{\mathscr{R}}=\coprod_{\nu\in\mathfrak{J}}
 \mathfrak{D}^{(\nu)},
\end{equation}
where $\mathfrak{J}$ is a set of indices.
Let $\mathscr{V}$ be the set of
all linearly independent subsets
$\mathbf{V}=\{\beta_1,\ldots,\beta_r\}\subset\Delta_+$ and
$\mathscr{A}=\{0,1\}^{n-r}$, where $n=\abs{\Delta_+}$.
Let $\mathscr{W}=\mathscr{V}\times\mathscr{A}$.
For subsets $u=\{u_1,\ldots,u_k\},v=\{v_1,\ldots,v_{k+1}\}\subset V$
and $c=\{c_1,\ldots,c_{k+1}\}\subset\mathbb{R}$,
let
\begin{equation}
  H(\mathbf{y};u,v,c)=
\det
\begin{pmatrix}
\langle u_1,v_1\rangle&\cdots&
\langle u_1,v_{k+1}\rangle
\\
\vdots&\ddots&\vdots
\\
\langle u_k,v_1\rangle&\cdots&
\langle u_k,v_{k+1}\rangle
\\
\langle \mathbf{y},v_1\rangle+c_1
&\cdots&
\langle \mathbf{y},v_{k+1}\rangle+c_{k+1}\\
  \end{pmatrix}.
\end{equation}
We give a simple description of the polytopes $\mathcal{P}_{\mathbf{m},\mathbf{y}}$ defined by \eqref{eq:D_m_y}.
For $\gamma\in\Delta_+$, $a\in\{0,1\}$ and $\mathbf{y}\in V$,
we define $\bd{u}(\gamma,a)\in\mathbb{R}^{n-r}$ by
\begin{equation}
  \bd{u}(\gamma,a)_\alpha
  =
  \begin{cases}
    (-1)^{1-a}\langle\alpha^\vee,\lambda_i\rangle,\qquad &\text{if }\gamma=\alpha_i\in\fs,\\
    (-1)^a\delta_{\alpha\gamma},\qquad &\text{if }\gamma\not\in\fs,
  \end{cases}
\end{equation}
where $\alpha$ runs over $\Delta_+\setminus\fs$,
and define $v(\gamma,a;\mathbf{y})\in\mathbb{R}$ by
\begin{equation}
  v(\gamma,a;\mathbf{y})
  =
  \begin{cases}
    (-1)^{1-a}(\{\langle \mathbf{y},\lambda_i\rangle\}+m_i-a)
    \qquad &\text{if }\gamma=\alpha_i\in\fs,\\
    (-1)^aa=-a,\qquad &\text{if }\gamma\not\in\fs.
  \end{cases}
\end{equation}
Further we define 
\begin{equation}
  \mathcal{H}_{\gamma,a}(\mathbf{y})=
    \{\bd{x}=(x_\alpha)_{\alpha\in\Delta_+\setminus\fs}\in\mathbb{R}^{n-r}~|~\bd{u}(\gamma,a)\cdot\bd{x}=v(\gamma,a;\mathbf{y})\},
\end{equation}
and
\begin{equation}
  \mathcal{H}_{\gamma,a}^+(\mathbf{y})=
    \{\bd{x}=(x_\alpha)_{\alpha\in\Delta_+\setminus\fs}\in\mathbb{R}^{n-r}~|~\bd{u}(\gamma,a)\cdot\bd{x}\geq v(\gamma,a;\mathbf{y})\},
\end{equation}
where for $\bd{w}=(w_\alpha),\bd{x}=(x_\alpha)\in\mathbb{C}^{n-r}$, we have set
\begin{equation}
  \bd{w}\cdot\bd{x}=\sum_{\alpha\in\Delta_+\setminus\fs}w_\alpha x_\alpha.
\end{equation}
Then we have
\begin{equation}
  \mathcal{P}_{\mathbf{m},\mathbf{y}}=
  \bigcap_{\substack{\gamma\in\Delta_+\\a\in\{0,1\}}}\mathcal{H}_{\gamma,a}^+(\mathbf{y}).
\end{equation}
We use the identification $\mathbb{C}\otimes V\simeq\mathbb{C}^r$
through $\mathbf{y}\mapsto(y_i)_{i=1}^r$ where
$y_i=\langle\mathbf{y},\lambda_i\rangle$ with $\langle\cdot,\cdot\rangle$
bilinearly extended over $\mathbb{C}$.
For $\mathbf{k}=(k_\alpha)_{\alpha\in\prs}\in\mathbb{N}_0^{n}$, we set
$\abs{\mathbf{k}}=\sum_{\alpha\in\prs}k_\alpha$.
\begin{theorem}
\label{thm:Bernoulli}
In each $\mathfrak{D}^{(\nu)}$,
the functions $F(\mathbf{t},\mathbf{y};\Delta)$ 
and $P(\mathbf{k},\mathbf{y};\Delta)$ are real analytic in $\mathbf{y}$.
Moreover,
$F(\mathbf{t},\mathbf{y};\Delta)$ is analytically 
continued 
to a meromorphic function $F^{(\nu)}(\mathbf{t},\mathbf{y};\Delta)$ 
from each $\mathbb{C}^{n}\times\mathfrak{D}^{(\nu)}$ to the whole space
$\mathbb{C}^{n}\times(\mathbb{C}\otimes V)$.
Similarly,
$P(\mathbf{k},\mathbf{y};\Delta)$ is analytically continued
to a polynomial function $B^{(\nu)}_{\mathbf{k}}(\mathbf{y};\Delta)\in\mathbb{Q}\,[\mathbf{y}]$
from each $\mathfrak{D}^{(\nu)}$ to the whole space $\mathbb{C}\otimes V$
with its total degree at most $\abs{\mathbf{k}}+n-r$.
\end{theorem}
\begin{proof}
Throughout this proof, we fix an index $\nu\in\mathfrak{J}$.
Note that $\{\langle\mathbf{y},\lambda_i\rangle\}=\langle\mathbf{y},\lambda_i\rangle$ holds for 
$\mathbf{y}\in\mathring{\mathfrak{D}}$.
We show this statement by several steps.
In the first three steps, 
we investigate the dependence of vertices of $\mathcal{P}_{\mathbf{m},\mathbf{y}}$ on $\mathbf{y}\in\mathfrak{D}^{(\nu)}$
and
in the last two steps, 
by use of this result and triangulation,
we show the analyticity of the generating function.
We fix $\mathbf{m}\in\mathbb{N}_0^r$ except in the last step.

(Step 1.)
Let $\mathbf{V}=\{\beta_1,\ldots,\beta_r\}\subset\Delta_+$ and $a_\gamma\in\{0,1\}$ for $\gamma\in\Delta_+\setminus\mathbf{V}$. 
Consider the intersection of $\abs{\Delta_+\setminus\mathbf{V}}(=n-r)$ hyperplanes
\begin{equation}
\label{eq:def_p}
  \bigcap_{\gamma\in\Delta_+\setminus\mathbf{V}}\mathcal{H}_{\gamma,a_\gamma}(\mathbf{y})
   =
    \{\bd{x}=(x_\alpha)~|~\bd{u}(\gamma,a_\gamma)\cdot\bd{x}=v(\gamma,a_\gamma;\mathbf{y})\text{ for }
\gamma\in\Delta_+\setminus\mathbf{V}\}.
\end{equation}
Then this set consists of the solutions of the system of the $(n-r)$ linear equations
\begin{equation}
\label{eq:eq_vert}
  \begin{cases}
    \sum_{\alpha\in\Delta_+\setminus \fs}x_\alpha\langle\alpha^\vee,\lambda_j\rangle
    =\langle \mathbf{y},\lambda_j\rangle+m_j-a_{\alpha_j},\qquad&\text{for }\gamma=\alpha_j\in\fs\setminus\mathbf{V},\\
    x_\gamma=a_\gamma,\qquad&\text{for }\gamma\in\Delta_+\setminus(\fs\cup \mathbf{V}).
  \end{cases}
\end{equation}
Let $I=\{i~|~\beta_i\in\mathbf{V}\setminus\fs\}$ and
$J=\{j~|~\alpha_j\in\fs\setminus\mathbf{V}\}$.
Note that $\abs{I}=\abs{J}=:k$ and $\{\beta_i~|~i\in I^c\}=\{\alpha_j~|~j\in J^c\}$.
The system of the linear equations \eqref{eq:eq_vert} has a unique solution
if and only if
\begin{equation}
\label{eq:det_nz}
\det
(\langle\beta_i^\vee,\lambda_j\rangle)_{j\in J}^{i\in I}
\neq0,
\end{equation}
and also if and only if
\begin{equation}
\label{eq:B_indept}
\mathbf{V}\in\mathscr{V},
\end{equation}
since
\begin{equation}
  \begin{split}
\bigl\lvert
\det
(\langle\beta_i^\vee,\lambda_j\rangle)_{1\leq j\leq r}^{1\leq i\leq r}
\bigr\rvert
&=
\Biggl\lvert
\det
\begin{pmatrix}
(\langle\beta_i^\vee,\lambda_j\rangle)_{j\in J}^{i\in I}
&
(\langle\beta_i^\vee,\lambda_j\rangle)_{j\in J^c}^{i\in I}
\\
(\langle\alpha_i^\vee,\lambda_j\rangle)_{j\in J}^{i\in J^c}
&
(\langle\alpha_i^\vee,\lambda_j\rangle)_{j\in J^c}^{i\in J^c}
\end{pmatrix}
\Biggr\rvert
\\
&=
\biggl\lvert
\det
\begin{pmatrix}
(\langle\beta_i^\vee,\lambda_j\rangle)_{j\in J}^{i\in I} & * \\
0 & I_{\abs{J^c}}
\end{pmatrix}
\biggr\rvert
\\
&=
\bigl\lvert
\det
(\langle\beta_i^\vee,\lambda_j\rangle)_{j\in J}^{i\in I}
\bigr\rvert,
\end{split}
\end{equation}
where $I_p$ is the $p\times p$ identity matrix.
We assume \eqref{eq:B_indept} and denote 
by $\bd{p}(\mathbf{y};\mathbf{W})$ the unique solution,
where $\mathbf{W}=(\mathbf{V},\mathbf{A})\in\mathscr{W}$ 
with the sequence
$\mathbf{A}=(a_\gamma)_{\gamma\in\Delta_+\setminus\mathbf{V}}$ regarded as an element of $\mathscr{A}$.
We see that $\bd{p}(\mathbf{y};\mathbf{W})$ depends on $\mathbf{y}$ affine-linearly.

(Step 2.)
We define $\xi_{\mathbf{y}}:\mathscr{W}\to\mathbb{R}^{n-r}$ by
\begin{equation}
\xi_{\mathbf{y}}:\mathbf{W}\mapsto \bd{p}(\mathbf{y};\mathbf{W}).
\end{equation}
Any vertex (that is, $0$-face) of $\mathcal{P}_{\mathbf{m},\mathbf{y}}$ is defined by the intersection of 
$(n-r)$
hyperplanes by Proposition \ref{prop:k_face_2}. Hence
$\Vrt(\mathcal{P}_{\mathbf{m},\mathbf{y}})\subset\xi_{\mathbf{y}}(\mathscr{W})$.
On the other hand, $\mathcal{P}_{\mathbf{m},\mathbf{y}}$ is defined 
by $n$ pairs of inequalities in
\eqref{eq:D_m_y}.   
The point $\bd{p}(\mathbf{y};\mathbf{W})$ is a vertex of $\mathcal{P}_{\mathbf{m},\mathbf{y}}$
if all of those inequalities hold.    
We see that $(n-r)$ pairs among them are satisfied,
because
\begin{equation}
\bd{p}(\mathbf{y};\mathbf{W})\in
  \bigcap_{\gamma\in\Delta_+\setminus\mathbf{V}}
\mathcal{H}_{\gamma,a_\gamma}(\mathbf{y}),
\end{equation}
and also it is easy to check
\begin{equation}
  \bd{p}(\mathbf{y};\mathbf{W})\in
  \bigcap_{\gamma\in\Delta_+\setminus\mathbf{V}}
  \bigl(\mathcal{H}_{\gamma,1-a_\gamma}^+(\mathbf{y})\setminus
  \mathcal{H}_{\gamma,1-a_\gamma}(\mathbf{y})\bigr).
\end{equation}
Therefore 
 $\bd{p}(\mathbf{y};\mathbf{W})\in\Vrt(\mathcal{P}_{\mathbf{m},\mathbf{y}})$ 
if and only if the remaining $r$ pairs of
inequalities are satisfied, that is,
\begin{equation}
\bd{p}(\mathbf{y};\mathbf{W})\in
  \bigcap_{\substack{\beta\in\mathbf{V}\\a\in\{0,1\}}}\mathcal{H}_{\beta,a}^+(\mathbf{y})
=
    \{\bd{x}=(x_\alpha)~|~\bd{u}(\beta,a)\cdot\bd{x}\geq v(\beta,a;\mathbf{y})\text{ for }
\beta\in\mathbf{V}, a\in\{0,1\}\},
\end{equation}
or equivalently $\bd{x}=\bd{p}(\mathbf{y};\mathbf{W})$ satisfies
 $r$ pairs of the linear inequalities
\begin{equation}
\label{eq:ineq_vert}
  \begin{cases}
    \langle \mathbf{y},\lambda_l\rangle+m_l-1\leq
    \sum_{\alpha\in\Delta_+\setminus \fs}x_\alpha\langle\alpha^\vee,\lambda_l\rangle
    \leq\langle \mathbf{y},\lambda_l\rangle+m_l,\qquad&\text{for }\beta=\alpha_l\in\mathbf{V}\cap\fs,\\
    0\leq x_\beta\leq1,\qquad&\text{for }\beta\in\mathbf{V}\setminus\fs.
  \end{cases}
\end{equation}
We see that it depends on $\mathbf{y}$
whether $\bd{p}(\mathbf{y};\mathbf{W})$ is a vertex,
or in other words,
whether the solution of \eqref{eq:eq_vert} satisfies \eqref{eq:ineq_vert}.
We will show in the next step that $\bd{p}(\mathbf{y};\mathbf{W})\in\mathcal{H}(\mathbf{y};\mathbf{W})$ implies $\mathbf{y}\in\mathfrak{H}_{\mathscr{R}}$,
where
\begin{equation}
\label{eq:H_y_W}
  \mathcal{H}(\mathbf{y};\mathbf{W})=
  \bigcup_{\substack{\beta\in\mathbf{V}\\a\in\{0,1\}}}\mathcal{H}_{\beta,a}(\mathbf{y}).
\end{equation}
Then for
$\mathbf{y}\in\mathfrak{D}\setminus\mathfrak{H}_{\mathscr{R}}$,
we can uniquely determine the $(n-r)$ hyperplanes
on which 
the point
$\bd{p}(\mathbf{y};\mathbf{W})$ lies;
they are
$\{\mathcal{H}_{\gamma,a_\gamma}(\mathbf{y})\}_{\gamma\in\Delta_+\setminus\mathbf{V}}$.
Therefore $\xi_{\mathbf{y}}$ is an injection.

For $\beta\in\mathbf{V}$ and $a\in\{0,1\}$, we define $f_{\beta,a}:\mathring{\mathfrak{D}}\to\mathbb{R}$ by
\begin{equation}
f_{\beta,a}:\mathbf{y}\mapsto\bd{u}(\beta,a)\cdot\bd{p}(\mathbf{y};\mathbf{W})-v(\beta,a;\mathbf{y}).
\end{equation}
Then for
$\mathbf{y}\in\mathfrak{D}\setminus\mathfrak{H}_{\mathscr{R}}$,
we have $f_{\beta,a}(\mathbf{y})\neq 0$ and hence we define
\begin{equation}
  f=(f_{\beta,a})_{\beta\in\mathbf{V},a\in\{0,1\}}:\mathfrak{D}\setminus\mathfrak{H}_{\mathscr{R}}\to(\mathbb{R}\setminus\{0\})^{2r}.
\end{equation}
Therefore
for
$\mathbf{y}\in\mathfrak{D}\setminus\mathfrak{H}_{\mathscr{R}}$,
the point
$\bd{p}(\mathbf{y};\mathbf{W})$
is a vertex if and only if
$f(\mathbf{y})$ is an element of the connected component $(0,\infty)^{2r}$.
Since each $f_{\beta,a}$ is continuous
and hence $f(\mathfrak{D}^{(\nu)})$ is connected,
we see that
for a fixed $\mathbf{W}\in\mathscr{W}$,
the point $\bd{p}(\mathbf{y};\mathbf{W})$ is always a vertex, 
or never a vertex, 
on $\mathfrak{D}^{(\nu)}$.
Thus
\begin{equation}
\mathscr{W}_{\mathbf{m}}:=\xi_{\mathbf{y}}^{-1}(\Vrt(\mathcal{P}_{\mathbf{m},\mathbf{y}}))\subset\mathscr{W}
\qquad(\mathbf{y}\in\mathfrak{D}^{(\nu)})
\end{equation}
has one-to-one correspondence with
$\Vrt(\mathcal{P}_{\mathbf{m},\mathbf{y}})$
and
is independent of $\mathbf{y}$
on $\mathfrak{D}^{(\nu)}$.

(Step 3.)
Now we prove the claim announced just before \eqref{eq:H_y_W}.
First we show that the condition
\begin{equation}
\label{eq:v_on_H} 
\bd{p}(\mathbf{y};\mathbf{W})\in\mathcal{H}_{\beta,a_\beta}(\mathbf{y})
\end{equation}
for some $\beta=\alpha_l\in\mathbf{V}\cap\fs$ and $a_{\beta}\in\{0,1\}$
implies $\mathbf{y}\in\mathfrak{H}_{\mathscr{R}}$.
For $\bd{x}=\bd{p}(\mathbf{y};\mathbf{W})$,
condition \eqref{eq:v_on_H} is equivalent to 
\begin{equation}
\label{eq:v_on_H_equiv}
    \sum_{\alpha\in\Delta_+\setminus \fs}x_\alpha\langle\alpha^\vee,\lambda_l\rangle
    =\langle \mathbf{y},\lambda_l\rangle+m_l-a_{\alpha_l}.
\end{equation}
From \eqref{eq:eq_vert} and \eqref{eq:v_on_H_equiv},
we have an overdetermined system with the $\abs{\mathbf{V}\setminus\fs}=k$ variables $x_\beta$ for $\beta\in\mathbf{V}\setminus\fs$
and the $\abs{(\fs\setminus\mathbf{V})\cup\{\alpha_l\}}=(k+1)$ equations
\begin{equation}
  \sum_{\beta\in\mathbf{V}\setminus\fs}x_{\beta}\langle\beta^\vee,\lambda_j\rangle
  =\langle \mathbf{y},\lambda_j\rangle+c_j,
\end{equation}
for $j\in J\cup\{l\}$,
where
\begin{equation}
  c_j=m_j-a_{\alpha_j}
  -\sum_{\gamma\in\Delta_+\setminus(\fs\cup\mathbf{V})}a_\gamma\langle\gamma^\vee,\lambda_j\rangle\in\mathbb{Z}.
\end{equation}
Hence we have
\begin{equation}
  \begin{pmatrix}
  x_{\beta_{i_1}}&\cdots&x_{\beta_{i_k}}&-1
  \end{pmatrix}
  \begin{pmatrix}
  \langle\beta_{i_1}^\vee,\lambda_{j_1}\rangle&\cdots&
  \langle\beta_{i_1}^\vee,\lambda_{j_k}\rangle&
  \langle\beta_{i_1}^\vee,\lambda_{l}\rangle\\
  \vdots & \ddots & \vdots & \vdots\\
  \langle\beta_{i_k}^\vee,\lambda_{j_1}\rangle&\cdots&
  \langle\beta_{i_k}^\vee,\lambda_{j_k}\rangle&
  \langle\beta_{i_k}^\vee,\lambda_{l}\rangle\\
  \langle \mathbf{y},\lambda_{j_1}\rangle+c_{j_1}&\cdots&
  \langle \mathbf{y},\lambda_{j_k}\rangle+c_{j_k}&
  \langle \mathbf{y},\lambda_{l}\rangle+c_{l}
  \end{pmatrix}
  =
  \begin{pmatrix}
    0 &\cdots& 0
  \end{pmatrix},
\end{equation}
where we have put $I=\{i_1,\ldots,i_k\}$ and $J=\{j_1,\ldots,j_k\}$.
As the consistency for these equations, we get 
\begin{equation}
\label{eq:hp_1}
  H(\mathbf{y};\{\beta_i^\vee\}_{i\in I},\{\lambda_j\}_{j\in J\cup\{l\}},\{c_j\}_{j\in J\cup\{l\}})=0.
\end{equation}
By direct substitution,
we see that each element of
\begin{equation}
  \{\beta^\vee-q\}_{\beta\in \mathbf{V}\setminus\{\alpha_l\}}\cup\{-q\},\qquad q=\sum_{j\in J\cup\{l\}} c_j \alpha_j^\vee
\end{equation}
satisfies \eqref{eq:hp_1}, while
$\alpha_l^\vee-q$ does not.
In fact, if $\mathbf{y}=-q$ or $\mathbf{y}=\beta^{\vee}-q$
($\beta\in(\mathbf{V}\cap \fs)\setminus\{\alpha_l\}$), then the last row of the
matrix is $(0,\ldots,0)$, and if $\mathbf{y}=\beta_{i_p}^\vee-q$ ($\beta_{i_p}\in\mathbf{V}\setminus \fs$), 
then the last
row is equal to the $p$-th row, and hence \eqref{eq:hp_1}
 follows, while if
$\mathbf{y}=\alpha_l^{\vee}-q$, then
the last row of the
matrix is $(0,\ldots,0,1)$ and hence
\begin{equation}
H(\mathbf{y};\{\beta_i^\vee\}_{i\in I},\{\lambda_j\}_{j\in J\cup\{l\}},\{c_j\}_{j\in J\cup\{l\}})
=\det(\langle\beta_i^{\vee}, \lambda_j\rangle)_{j\in J}^{i\in I}\neq 0,
\end{equation}
because of \eqref{eq:det_nz}.
By \eqref{eq:B_indept}, 
we see that
$\mathbf{V}\setminus\{\alpha_l\}\subset\Delta$
is a linearly independent subset 
and hence 
$(\mathbf{V}\setminus\{\alpha_l\})\in\mathscr{R}$.
It follows that \eqref{eq:hp_1} represents the hyperplane 
$\mathfrak{H}_{\mathbf{V}^\vee\setminus\{\alpha_l^\vee\}}-q\subset\mathfrak{H}_{\mathscr{R}}$.
Therefore
\eqref{eq:v_on_H} 
implies $\mathbf{y}\in\mathfrak{H}_{\mathscr{R}}$.

Similarly we see that
the condition 
 $\bd{p}(\mathbf{y};\mathbf{W})\in\mathcal{H}_{\beta,a_\beta}(\mathbf{y})$ for some $\beta=\beta_l\in\mathbf{V}\setminus\fs$ and $a_\beta\in\{0,1\}$
yields a hyperplane contained in $\mathfrak{H}_{\mathscr{R}}$ defined by
\begin{equation}
\label{eq:hp_2}
  H(\mathbf{y};\{\beta_i^\vee\}_{i\in I\setminus\{l\}},\{\lambda_j\}_{j\in J},\{d_j\}_{j\in J})=0,
\end{equation}
which
passes through $r$ points in general position
\begin{equation}
  \{\beta^\vee-q\}_{\beta\in \mathbf{V}\setminus\{\beta_l\}}\cup\{-q\},
\end{equation}
where
\begin{gather}
 q=\sum_{j\in J} d_j \alpha_j^\vee,\\
  d_j=m_j-a_{\alpha_j}-\sum_{\gamma\in\Delta_+\setminus(\fs\cup(\mathbf{V}\setminus\{\beta_l\}))}a_\gamma\langle\gamma^\vee,\lambda_j\rangle\in\mathbb{Z}.
\end{gather}
This completes the proof of our claim.

(Step 4.)
We have checked that 
on $\mathfrak{D}^{(\nu)}$,
the vertices
$\Vrt(\mathcal{P}_{\mathbf{m},\mathbf{y}})$ 
neither increase nor decrease and are indexed by $\mathscr{W}_{\mathbf{m}}$.
By numbering $\mathscr{W}_{\mathbf{m}}$ as $\{\mathbf{W}_1,\mathbf{W}_2,\ldots\}$,
we denote $\bd{p}_i(\mathbf{y})=\bd{p}(\mathbf{y};\mathbf{W}_i)$.
We see that on $\mathfrak{D}^{(\nu)}$, the polytopes
 $\mathcal{P}_{\mathbf{m},\mathbf{y}}$ keep $(n-r)$-dimensional or empty
because each vertex is determined by unique $(n-r)$ hyperplanes.
Assume that $\mathcal{P}_{\mathbf{m},\mathbf{y}}$ is not empty.
Next we will show that the 
face poset structure of
$\mathcal{P}_{\mathbf{m},\mathbf{y}}$ is independent of $\mathbf{y}$ on $\mathfrak{D}^{(\nu)}$.

Fix $\mathbf{y}_0\in\mathfrak{D}^{(\nu)}$. 
Consider a face $\mathcal{F}(\mathbf{y}_0)$ of $\mathcal{P}_{\mathbf{m},\mathbf{y}_0}$
and let
\begin{equation}
\label{eq:Vrt_F_y_0}
 \Vrt(\mathcal{P}_{\mathbf{m},\mathbf{y}_0})\cap\mathcal{F}(\mathbf{y}_0)=\{\bd{p}_{i_1}(\mathbf{y}_0),\ldots,\bd{p}_{i_h}(\mathbf{y}_0)\}. 
\end{equation}
Then by Proposition \ref{prop:k_face_2},
there exists a subset $\mathcal{J}_0=\{(\gamma_1,a_{\gamma_1}),(\gamma_2,a_{\gamma_2}),\ldots\}
\subset\Delta_+\times\{0,1\}$
such that $\abs{\mathcal{J}_0}=n-r-\dim\mathcal{F}(\mathbf{y}_0)$ and
\begin{equation}
\label{eq:p_contained}
\{\bd{p}_{i_1}(\mathbf{y}_0),\ldots,\bd{p}_{i_h}(\mathbf{y}_0)\}
\subset
\mathcal{F}(\mathbf{y}_0)=\mathcal{P}_{\mathbf{m},\mathbf{y}_0}\cap\bigcap_{(\gamma,a_\gamma)\in\mathcal{J}_0}\mathcal{H}_{\gamma,a_\gamma}(\mathbf{y}_0).
\end{equation}
By the definition of $\bd{p}_{i_1}(\mathbf{y}_0), \ldots,\bd{p}_{i_h}(\mathbf{y}_0)$ (see \eqref{eq:def_p}), we find
that for $(\gamma,a_\gamma)\in\Delta_+\times\{0,1\}$, the condition
\begin{equation}
\{\bd{p}_{i_1}(\mathbf{y}_0),\ldots,\bd{p}_{i_h}(\mathbf{y}_0)\}
\subset
\mathcal{H}_{\gamma,a_{\gamma}}(\mathbf{y}_0)
\end{equation}
is equivalent to 
\begin{equation}
\label{eq:common_cond}
\gamma\in\bigcap_{j=1}^h(\Delta_+\setminus\mathbf{V}_{i_j}),\qquad
a_{\gamma}=(\mathbf{A}_{i_1})_\gamma=\cdots=(\mathbf{A}_{i_h})_\gamma,
\end{equation}
where
$\mathbf{W}_{i_j}=(\mathbf{V}_{i_j},\mathbf{A}_{i_j})=\xi_{\mathbf{y}_0}^{-1}(\bd{p}_{i_j}(\mathbf{y}_0))$. 
Hence each $(\gamma,a_\gamma)\in\mathcal{J}_0$ satisfies \eqref{eq:common_cond}.
Assume that there exists a pair $(\gamma',a'_{\gamma'})\not\in\mathcal{J}_0$ satisfying \eqref{eq:common_cond},
Then 
\begin{equation}
\{\bd{p}_{i_1}(\mathbf{y}_0),\ldots,\bd{p}_{i_h}(\mathbf{y}_0)\}
\subset
\bigcap_{(\gamma,a_\gamma)\in\mathcal{J}'_0}\mathcal{H}_{\gamma,a_\gamma}(\mathbf{y}_0),
\end{equation}
where $\mathcal{J}'_0=\mathcal{J}_0\cup\{(\gamma',a'_{\gamma'})\}$.
Hence 
by \eqref{eq:prop_F}, we have
\begin{equation}
\label{eq:F_contr}
\mathcal{F}(\mathbf{y}_0)=\Conv\{\bd{p}_{i_1}(\mathbf{y}_0),\ldots,\bd{p}_{i_h}(\mathbf{y}_0)\}
\subset
\bigcap_{(\gamma,a_\gamma)\in\mathcal{J}'_0}\mathcal{H}_{\gamma,a_\gamma}(\mathbf{y}_0),
\end{equation}
and in particular,
\begin{equation}
\bd{p}_{i_1}(\mathbf{y}_0)\in
\bigcap_{(\gamma,a_\gamma)\in\mathcal{J}'_0}\mathcal{H}_{\gamma,a_\gamma}(\mathbf{y}_0).
\end{equation}
Since 
$(n-r)$ hyperplanes on which $\bd{p}_{i_1}(\mathbf{y}_0)$ lies
are uniquely determined and their intersection consists of only $\bd{p}_{i_1}(\mathbf{y}_0)$,
their normal vectors $\{\bd{u}(\gamma,a_\gamma)\}_{(\gamma,a_\gamma)\in\mathcal{J}'_0}$ must be linearly independent.
It follows that the dimension of the right-hand side of \eqref{eq:F_contr} is
$n-r-\abs{\mathcal{J}'_0}<n-r-\abs{\mathcal{J}_0}=\dim\mathcal{F}(\mathbf{y}_0)$,
which contradicts.
Hence \eqref{eq:common_cond} is also a sufficient condition for
$(\gamma,a_\gamma)\in\mathcal{J}_0$.

By definition,
we have
\begin{equation}
\{\bd{p}_{i_1}(\mathbf{y}),\ldots,\bd{p}_{i_h}(\mathbf{y})\}
\subset
\bigcap_{(\gamma,a_\gamma)\in\mathcal{J}_0}
\mathcal{H}_{\gamma,a_\gamma}(\mathbf{y}),
\end{equation}
for all $\mathbf{y}\in\mathfrak{D}^{(\nu)}$.
Define
\begin{equation}
\label{eq:def_Fy}
 \mathcal{F}(\mathbf{y})=\mathcal{P}_{\mathbf{m},\mathbf{y}}\cap\bigcap_{(\gamma,a_\gamma)\in\mathcal{J}_0}\mathcal{H}_{\gamma,a_\gamma}(\mathbf{y}).  
\end{equation}
Then $\{\bd{p}_{i_1}(\mathbf{y}),\ldots,\bd{p}_{i_h}(\mathbf{y})\}\subset\mathcal{F}(\mathbf{y})$ and by Proposition \ref{prop:k_face_1},
we see that $\mathcal{F}(\mathbf{y})$ is a face.
Fix another $\mathbf{y}_1\in\mathfrak{D}^{(\nu)}$.
Then by the argument at the beginning of this step,
there exists a subset $\mathcal{J}_1\subset\Delta_+\times\{0,1\}$ such that
$\abs{\mathcal{J}_1}=n-r-\dim\mathcal{F}(\mathbf{y}_1)$ and
\begin{equation}
\label{eq:def_J1}
  \begin{split}
 \Vrt(\mathcal{P}_{\mathbf{m},\mathbf{y}_1})\cap\mathcal{F}(\mathbf{y}_1)
&=\{\bd{p}_{i_1}(\mathbf{y}_1),\ldots,\bd{p}_{i_h}(\mathbf{y}_1),\bd{p}_{i_{h+1}}(\mathbf{y}_1),\ldots,\bd{p}_{i_{h'}}(\mathbf{y}_1)\}
\\
&\subset\mathcal{F}(\mathbf{y}_1)
=\mathcal{P}_{\mathbf{m},\mathbf{y}_1}\cap\bigcap_{(\gamma,a_\gamma)\in\mathcal{J}_1}
\mathcal{H}_{\gamma,a_\gamma}(\mathbf{y}_1), 
\end{split}
\end{equation}
with $h'\geq h$
(because all of $\bd{p}_{i_1}(\mathbf{y}_1),\ldots ,\bd{p}_{i_h}(\mathbf{y}_1)$ are
vertices of $\mathcal{F}(\mathbf{y}_1)$, while so far we cannot exclude the possibility of 
the existence of other vertices on $\mathcal{F}(\mathbf{y}_1)$).
Since each $(\gamma,a_\gamma)\in\mathcal{J}_1$ satisfies
\begin{equation}
\label{eq:common_cond_s}
\gamma\in\bigcap_{j=1}^{h'}(\Delta_+\setminus\mathbf{V}_{i_j}),\qquad
a_{\gamma}=(\mathbf{A}_{i_1})_\gamma=\cdots=(\mathbf{A}_{i_{h'}})_\gamma,
\end{equation}
which is equal to or stronger than condition \eqref{eq:common_cond},
we have $\mathcal{J}_1\subset\mathcal{J}_0$.
On the other hand,
by comparing \eqref{eq:def_Fy} and \eqref{eq:def_J1},
we see that each
$(\gamma,a_\gamma)\in\mathcal{J}_0$ satisfies \eqref{eq:common_cond_s}.
As shown in the previous paragraph, 
condition \eqref{eq:common_cond_s} is sufficient for $(\gamma,a_\gamma)\in\mathcal{J}_1$,
which implies $\mathcal{J}_0=\mathcal{J}_1$ and hence $\dim\mathcal{F}(\mathbf{y}_1)=\dim\mathcal{F}(\mathbf{y}_0)$.
If $h'>h$, then \eqref{eq:common_cond_s} implies
\begin{equation}
  \{\bd{p}_{i_1}(\mathbf{y}_0),\ldots,\bd{p}_{i_h}(\mathbf{y}_0),\bd{p}_{i_{h+1}}(\mathbf{y}_0),\ldots,\bd{p}_{i_{h'}}(\mathbf{y}_0)\}
\subset\mathcal{F}(\mathbf{y}_0),
\end{equation}
which contradicts to \eqref{eq:Vrt_F_y_0} and hence $h'=h$.
Therefore
for all $\mathbf{y}\in\mathfrak{D}^{(\nu)}$,
we see that all faces of $\mathcal{P}_{\mathbf{m},\mathbf{y}}$ 
are determined at $\mathbf{y}_0$ and
are described in the form \eqref{eq:def_Fy},
and we have
\begin{equation}
 \Vrt(\mathcal{P}_{\mathbf{m},\mathbf{y}})\cap\mathcal{F}(\mathbf{y})=\{\bd{p}_{i_1}(\mathbf{y}),\ldots,\bd{p}_{i_h}(\mathbf{y})\}.
\end{equation}

Assume that $\mathcal{F}'(\mathbf{y}_0)\subset\mathcal{F}(\mathbf{y}_0)$ 
for faces 
$\mathcal{F}'(\mathbf{y}_0),\mathcal{F}(\mathbf{y}_0)$ 
of $\mathcal{P}_{\mathbf{m},\mathbf{y}_0}$.
Then by \eqref{eq:prop_F}, it is equivalent to
\begin{equation}
  \Vrt(\mathcal{P}_{\mathbf{m},\mathbf{y}_0})\cap\mathcal{F}'(\mathbf{y}_0)
\subset
  \Vrt(\mathcal{P}_{\mathbf{m},\mathbf{y}_0})\cap\mathcal{F}(\mathbf{y}_0).
\end{equation}
By applying $\xi_{\mathbf{y}_0}^{-1}$, we obtain an equivalent condition 
independent of $\mathbf{y}$ and hence
$\mathcal{F}'(\mathbf{y}_1)\subset\mathcal{F}(\mathbf{y}_1)$.
Therefore the face poset structure is indeed 
independent of $\mathbf{y}$ on $\mathfrak{D}^{(\nu)}$.

(Step 5.)
By Theorem \ref{thm:simp_div}, we have a triangulation of $\mathcal{P}_{\mathbf{m},\mathbf{y}}$
with $(n-r)$-dimensional simplexes $\sigma_{l,\mathbf{m},\mathbf{y}}$ as
\begin{equation}
\label{eq:sigma_kmy}
\mathcal{P}_{\mathbf{m},\mathbf{y}}=\bigcup_{l=1}^{L(\mathbf{m},\mathbf{y})} \sigma_{l,\mathbf{m},\mathbf{y}},
\end{equation}
where $L(\mathbf{m},\mathbf{y})$ is the number of the simplexes.  
From the previous step and by Remark \ref{rem:face_poset}, we see that
this triangulation does not depend on $\mathbf{y}$ up to the order of simplexes, i.e.,
\begin{equation}
  \{\mathcal{I}(1,\mathbf{m},\mathbf{y}),\ldots,\mathcal{I}(L(\mathbf{m},\mathbf{y}),\mathbf{m},\mathbf{y})\}
\end{equation}
is independent of $\mathbf{y}$,
where $\mathcal{I}(l,\mathbf{m},\mathbf{y})$
is the set of all indices of the vertices of $\sigma_{l,\mathbf{m},\mathbf{y}}$.
Reordering $\sigma_{l,\mathbf{m},\mathbf{y}}$ with respect to $l$ if necessary,
we assume that 
each $\mathcal{I}(l,\mathbf{m})=\mathcal{I}(l,\mathbf{m},\mathbf{y})$ is independent of $\mathbf{y}$ on $\mathfrak{D}^{(\nu)}$.
Note that $\abs{\mathcal{I}(l,\mathbf{m})}=n-r+1$. 
Let $L(\mathbf{m})=L(\mathbf{m},\mathbf{y})$ if $\mathcal{P}_{\mathbf{m},\mathbf{y}}$ is not empty, and $L(\mathbf{m})=0$ otherwise.

By \eqref{eq:sigma_kmy} and Lemma \ref{lm:simplex},
we find that the integral in the third expression of \eqref{eq:gen_func} is
\begin{equation}
  (n-r)!\sum_{l=1}^{L(\mathbf{m})}
 \Vol(\sigma_{l,\mathbf{m},\mathbf{y}})
\sum_{i\in\mathcal{I}(l,\mathbf{m})}
\frac{e^{\bd{t}^*\cdot\bd{p}_i(\mathbf{y})}}
{\prod_{\substack{j\in\mathcal{I}(l,\mathbf{m})\\j\neq i}}\bd{t}^*\cdot (\bd{p}_i(\mathbf{y})-\bd{p}_j(\mathbf{y}))},
\end{equation}
where $\bd{t}^*=(t^*_\alpha)$ with $t^*_\alpha=t_\alpha-\sum_{i=1}^rt_{\alpha_i}\langle\alpha^\vee,\lambda_i\rangle$.
We see that $\Vol(\sigma_{l,\mathbf{m},\mathbf{y}})$ is a polynomial function
in $\mathbf{y}=(y_i)_{i=1}^r$
with rational coefficients and
its total degree at most $n-r$
on $\mathfrak{D}^{(\nu)}$
 due to \eqref{eq:vol_simplex}, because in Step 1 we have 
shown that $\bd{p}_i(\mathbf{y})$ depends on $\mathbf{y}$ affine-linearly.
Therefore from \eqref{eq:gen_func},
we have the generating function
\begin{equation}
\label{eq:F_final}
F(\mathbf{t},\mathbf{y};\Delta)=
  \biggl(\prod_{\alpha\in\Delta_+}
  \frac{t_\alpha}{e^{t_\alpha}-1}\biggr)\sum_{j=1}^{(n-r+1)\sum_{\mathbf{m}}L(\mathbf{m})}
  \frac{f_j(\mathbf{y})e^{h_j(\mathbf{t},\mathbf{y})}}{g_j(\mathbf{t},\mathbf{y})},
\end{equation}
which is valid for all $\mathbf{y}\in\mathfrak{D}^{(\nu)}$,
where $f_j\in \mathbb{Q}\,[\mathbf{y}]$ 
with
its total degree at most $n-r$ and $g_j\in\mathbb{Z}\,[\mathbf{t},\mathbf{y}]$, $h_j\in\mathbb{Q}\,[\mathbf{t},\mathbf{y}]$ are of the form
\begin{align}
g_j(\mathbf{t},\mathbf{y})&=
\sum_{\alpha\in\Delta_+}(\langle \phi_\alpha, \mathbf{y}\rangle +c_\alpha)t_\alpha,
\qquad \phi_\alpha\in P,\quad c_\alpha\in\mathbb{Z},
\\
h_j(\mathbf{t},\mathbf{y})&=
\sum_{\alpha\in\Delta_+}(\langle \varphi_\alpha, \mathbf{y}\rangle +d_\alpha)t_\alpha,
\qquad \varphi_\alpha\in \mathbb{Q}\otimes P,\quad d_\alpha\in\mathbb{Q}.
\end{align}
We see from \eqref{eq:F_final} that
$F(\mathbf{t},\mathbf{y};\Delta)$ 
is meromorphically continued
from
 $\mathbb{C}^{n}\times\mathfrak{D}^{(\nu)}$
 to
the whole space $\mathbb{C}^{n}\times(\mathbb{C}\otimes V)$,
and
that $P(\mathbf{k},\mathbf{y};\Delta)$ is analytically continued to a polynomial function
in $\mathbf{y}=(y_i)_{i=1}^r$
with rational coefficients 
and its total degree at most $\abs{\mathbf{k}}+n-r$
by \eqref{eq:def_F} and Lemma \ref{lm:simplex_expand}.
\end{proof}

\begin{theorem}
The function 
$P(\mathbf{k},\mathbf{y};\Delta)$ is not real analytic in $\mathbf{y}$ on $V$
unless $P(\mathbf{k},\mathbf{y};\Delta)$ is a constant.
\end{theorem}
\begin{proof}
  By Theorem 
  \ref{thm:inv_qv}, we see that for $\mathbf{k}\in\mathbb{N}_0^n$,
  $P(\mathbf{k},\mathbf{y};\Delta)$ is a periodic function in $\mathbf{y}$ with its periods $\fs^\vee$,
  while by Theorem \ref{thm:Bernoulli},
  $P(\mathbf{k},\mathbf{y};\Delta)$ is a polynomial function in $\mathbf{y}$
  on some open region.
  Therefore such a polynomial expression
  cannot be extended to the whole space unless $P(\mathbf{k},\mathbf{y};\Delta)$ is
  a constant. This implies that there are some points on $\mathfrak{H}_{\mathscr{R}}$,
  at which $P(\mathbf{k},\mathbf{y};\Delta)$ is not real analytic.
\end{proof}

The polynomials $B^{(\nu)}_{\mathbf{k}}(\mathbf{y};\Delta)$
may be regarded 
as (root-system theoretic) generalizations
of Bernoulli polynomials.
For instance, they possess the following property.
\begin{theorem}
Assume that $\Delta$ is an irreducible root system and
is not of type $A_1$.
For $\mathbf{k}\in\mathbb{N}_0^{n}$,
$\mathbf{y}\in\partial\mathfrak{D}^{(\nu)}$ and $\mathbf{y}'\in\partial\mathfrak{D}^{(\nu')}$
with $\mathbf{y}\equiv\mathbf{y}'\pmod{Q^\vee}$,
we have
  \begin{equation}
   B^{(\nu)}_{\mathbf{k}}(\mathbf{y};\Delta)
   =
   B^{(\nu')}_{\mathbf{k}}(\mathbf{y}';\Delta).
  \end{equation}
\end{theorem}
\begin{proof}
If $\mathbf{y}=\mathbf{y}'$, then the result follows from the
continuity proved in Lemma \ref{lm:cont_PF}.   
If $\mathbf{y}\neq\mathbf{y}'$, but
$\mathbf{y}\equiv\mathbf{y}'\pmod{Q^{\vee}}$, we also use the
periodicity.  
\end{proof}
This theorem also holds in the $A_1$ case with $k\neq 1$ and 
can be regarded as a multiple analogue of the formula for the classical Bernoulli polynomials
\begin{equation}
  B_k(0)=B_k(1),
\end{equation}
for $k\neq1$.
Moreover 
the formula
\begin{equation}
\label{eq:Bernoulli_property}
  B_k(1-y)=(-1)^k B_k(y)
\end{equation}
is well-known. In the rest of this section we will show the results
analogous to the above formula 
for $B^{(\nu)}_{\mathbf{k}}(\mathbf{y};\Delta)$
(Theorem \ref{thm:A_act_on_B}), 
and 
its vector-valued version
(Theorem \ref{thm:A_act_on_VB}).
The latter gives a finite dimensional representation of Weyl groups.
In this framework, \eqref{eq:Bernoulli_property} can be
interpreted as an action of the Weyl group of type $A_1$ (Example \ref{eg:exam_A1}).
These results
will not be used in the present paper, but we insert this topic because
of its own interest.

\begin{lemma}
\label{lm:A_act_on_D}
Fix $w\in\Aut$ and $\nu\in \mathfrak{J}$. Then there exist unique $q_{w,\nu}\in Q^\vee$ and $\kappa\in \mathfrak{J}$ such that
\begin{equation}
  \label{eq:A_act_on_D}
  \tau(q_{w,\nu})w\mathfrak{D}^{(\nu)}=\mathfrak{D}^{(\kappa)}.
\end{equation}
Thus $\Aut$ acts on $\mathfrak{J}$
as $w(\nu)=\kappa$. Moreover $q_{v,w(\nu)}+v q_{w,\nu}=q_{vw,\nu}$ for $v,w\in\Aut$.
\end{lemma}
\begin{proof}
It can be easily seen from the definition \eqref{eq:def_H_R}
that $\mathfrak{H}_{\mathscr{R}}$ is $\widehat{W}$-invariant.
Therefore $\widehat{W}$ acts on $V\setminus\mathfrak{H}_{\mathscr{R}}$ as homeomorphisms and
a connected component is mapped to
another one.

Fix $\mathbf{y}\in w\mathfrak{D}^{(\nu)}$.
There exists a unique $q\in Q^\vee$ such that
 $0\leq\langle\tau(q)\mathbf{y},\lambda_j\rangle<1$ for $1\leq j\leq r$, that is,
 $q=-\sum_{j=1}^r a_j\alpha_j^\vee\in Q^\vee$,
where $a_j=[\langle\mathbf{y},\lambda_j\rangle]$ is
the integer part of $\langle\mathbf{y},\lambda_j\rangle$.
Denote this $q$ by $q_{w,\nu}$.
Then $\tau(q_{w,\nu})\mathbf{y}\in\mathfrak{D}^{(\kappa)}$ for some $\kappa\in \mathfrak{J}$ and
thus $\tau(q_{w,\nu})w\mathfrak{D}^{(\nu)}=\mathfrak{D}^{(\kappa)}$.

Let $w,w'\in\Aut$.
Assume
\begin{equation}
  \tau(q)w\mathfrak{D}^{(\nu)}=\mathfrak{D}^{(\kappa)},\qquad
  \tau(q')w'\mathfrak{D}^{(\kappa)}=\mathfrak{D}^{(\kappa')},\qquad
  \tau(q'')w'w\mathfrak{D}^{(\nu)}=\mathfrak{D}^{(\kappa'')},
\end{equation}
for $q,q',q''\in Q^\vee$ and $\nu,\kappa,\kappa',\kappa''\in \mathfrak{J}$.
Then we have
\begin{equation}
  \begin{split}
    \mathfrak{D}^{(\kappa')}
    &=
    \tau(q')w'\tau(q)w\mathfrak{D}^{(\nu)}\\
    &=
    \tau(q'+w'q)w'w\mathfrak{D}^{(\nu)}.
\end{split}
\end{equation}
By the uniqueness, we have $q'+w'q=q''$ and $\kappa'=\kappa''$.
\end{proof}
\begin{theorem}
\label{thm:A_act_on_B}
For $w\in\Aut$,
\begin{equation}
\label{eq:A_act_on_B}
  B^{(\nu)}_{\mathbf{k}}(\tau(q_{w^{-1},w(\nu)})w^{-1}\mathbf{y};\Delta)
  =
  \Bigl(\prod_{\alpha\in\Delta_w}(-1)^{-k_\alpha}\Bigr)
  B^{(w(\nu))}_{w\mathbf{k}}(\mathbf{y};\Delta).
\end{equation}
\end{theorem}
\begin{proof}
By
Theorems \ref{thm:gen_sym} and \ref{thm:inv_qv},
we have
\begin{equation}
  \begin{split}
    P(\mathbf{k},\tau(q_{w^{-1},w(\nu)})w^{-1}\mathbf{y};\Delta)
    &=
    P(\mathbf{k},w^{-1}\mathbf{y};\Delta)
    \\
    &=
    \Bigl(\prod_{\alpha\in\Delta_{w^{-1}}}(-1)^{-k_{w^{-1}\alpha}}\Bigr)
    P(w\mathbf{k},\mathbf{y};\Delta)
    \\
    &=
    \Bigl(\prod_{\alpha\in\Delta_w}(-1)^{-k_\alpha}\Bigr)
    P(w\mathbf{k},\mathbf{y};\Delta),
  \end{split}
\end{equation}
where we have used $w^{-1}\Delta_{w^{-1}}=-\Delta_w$.
By \eqref{eq:A_act_on_D} 
in Lemma \ref{lm:A_act_on_D}
with replacing $w$, $\nu$ by
$w^{-1}$, $w(\nu)$ respectively,
we have
\begin{equation}
  \tau(q_{w^{-1},w(\nu)})w^{-1}\mathfrak{D}^{(w(\nu))}
=
  \mathfrak{D}^{(\nu)}.
\end{equation}
Hence $\mathbf{y}\in\mathfrak{D}^{(w(\nu))}$
implies $\tau(q_{w^{-1},w(\nu)})w^{-1}\mathbf{y}\in\mathfrak{D}^{(\nu)}$.   Therefore
we obtain
\begin{align}
  P(\mathbf{k},\tau(q_{w^{-1},w(\nu)})w^{-1}\mathbf{y};\Delta)
  &=
  B^{(\nu)}_{\mathbf{k}}(\tau(q_{w^{-1},w(\nu)})w^{-1}\mathbf{y};\Delta),
  \\
  P(\mathbf{k},\mathbf{y};\Delta)
  &=
  B^{(w(\nu))}_{\mathbf{k}}(\mathbf{y};\Delta),
\end{align}
by Theorem \ref{thm:Bernoulli}.
The theorem of identity implies \eqref{eq:A_act_on_B}
\end{proof}

Let $\mathfrak{P}$ be the $\mathbb{Q}\,$-vector space
 of all vector-valued polynomial functions of the form
$f=(f_\nu)_{\nu\in\mathfrak{J}}:V\to\mathbb{R}^{\abs{\mathfrak{J}}}$ with $f_\nu\in\mathbb{Q}\,[\mathbf{y}]$.
We define
a linear map $\phi(w):\mathfrak{P}\to\mathfrak{P}$
for $w\in\Aut$ by
\begin{equation}
\label{eq:act_phi}
  (\phi(w)f)_\nu(\mathbf{y})=
f_{w^{-1}(\nu)}(\tau(q_{w^{-1},\nu})w^{-1}\mathbf{y}).
\end{equation}
\begin{theorem}
\label{thm:rep_on_poly}
The pair $(\phi,\mathfrak{P})$ is a representation of $\Aut$.
\end{theorem}
\begin{proof}
For $v,w\in\Aut$,
we have
\begin{equation}
  \begin{split}
  (\phi(v)\phi(w)f)_\nu(\mathbf{y})
&=
(\phi(w)f)_{v^{-1}(\nu)}(\tau(q_{v^{-1},\nu})v^{-1}\mathbf{y})
\\
&=
f_{w^{-1}v^{-1}(\nu)}(\tau(q_{w^{-1},v^{-1}(\nu)})w^{-1}\tau(q_{v^{-1},\nu})v^{-1}\mathbf{y})
\\
&=
f_{(vw)^{-1}(\nu)}(\tau(q_{w^{-1},v^{-1}(\nu)}+w^{-1}q_{v^{-1},\nu})(vw)^{-1}\mathbf{y}).
\end{split}
\end{equation}
Since
by Lemma \ref{lm:A_act_on_D} we have
 $q_{w^{-1},v^{-1}(\nu)}+w^{-1}q_{v^{-1},\nu}=q_{(vw)^{-1},\nu}$, we obtain $\phi(vw)=\phi(v)\phi(w)$.
\end{proof}
Define $\mathbf{B}^{(\nu)}_{\mathbf{k}}\in\mathfrak{P}$ for $\mathbf{k}\in\mathbb{N}_0^{n}$ and $\nu\in\mathfrak{J}$ by
\begin{equation}
(\mathbf{B}^{(\nu)}_{\mathbf{k}})_\kappa(\mathbf{y})=
(\mathbf{B}^{(\nu)}_{\mathbf{k}})_\kappa(\mathbf{y};\Delta)=
\begin{cases}
B^{(\nu)}_{\mathbf{k}}(\mathbf{y};\Delta),\qquad&\text{if }\nu=\kappa,\\ 
0,&\text{otherwise},
\end{cases}
\end{equation}
and 
let
\begin{equation}
  \mathfrak{B}_{\overline{(\mathbf{k},\nu)}}=\sum_{(\mathbf{k}',\nu')\in\overline{(\mathbf{k},\nu)}}\mathbb{Q}\,\mathbf{B}^{(\nu')}_{\mathbf{k}'}\subset\mathfrak{P},
\end{equation}
where $\overline{(\mathbf{k},\nu)}$ is an element of the orbit space 
$(\mathbb{N}_0^{n}\times\mathfrak{J})/\Aut$.
\begin{theorem}
  \label{thm:A_act_on_VB} 
The vector subspace $\mathfrak{B}_{\overline{(\mathbf{k},\nu)}}$ is a finite 
dimensional $\Aut$-invariant subspace and the action is
\begin{equation}
  \label{eq:A_act_on_VB} 
  \phi(w)\mathbf{B}^{(\nu)}_{\mathbf{k}}=
    \Bigl(\prod_{\alpha\in\Delta_w}(-1)^{-k_{\alpha}}\Bigr)
    \mathbf{B}^{(w(\nu))}_{w\mathbf{k}}.
\end{equation}
\end{theorem}
\begin{proof}
If $\kappa=w(\nu)$, then we have
\begin{equation}
  \begin{split}
    (\phi(w)\mathbf{B}^{(\nu)}_{\mathbf{k}})_\kappa(\mathbf{y};\Delta)
    &=
    B^{(\nu)}_{\mathbf{k}}(\tau(q_{w^{-1},w(\nu)})w^{-1}\mathbf{y};\Delta)
    \\
    &=
    \Bigl(\prod_{\alpha\in\Delta_w}(-1)^{-k_{\alpha}}\Bigr)
    B^{(w(\nu))}_{w\mathbf{k}}(\mathbf{y};\Delta),
    \end{split}
  \end{equation}
by Theorem \ref{thm:A_act_on_B}
and otherwise
\begin{equation}
  (\phi(w)\mathbf{B}^{(\nu)}_{\mathbf{k}})_\kappa(\mathbf{y})=0.
\end{equation}
Thus we obtain \eqref{eq:A_act_on_VB}.
\end{proof}

\section{Examples}\label{sec-7}

\label{sec:examples}
\begin{example}
\label{eg:exam_A1}
The set of positive roots of type $A_1$ consists of only one root $\alpha_1$.
Hence we have $\Delta_+=\fs=\{\alpha_1\}$ and $2\rho^\vee=\alpha_1^\vee$. We set $t=t_{\alpha_1}$ and $y=\langle\mathbf{y},\lambda_1\rangle$.
Then by Theorem \ref{thm:gen_func}, 
we obtain the generating function $F(\mathbf{t},\mathbf{y};A_1)$ 
as
\begin{equation}
  F(\mathbf{t},\mathbf{y};A_1)=\frac{t}{e^t-1}e^{t\{y\}}.
\end{equation}
Since $\mathfrak{D}=\{\mathbf{y}~|~0\leq y\leq1\}$ and $\mathscr{R}=\{\emptyset\}$, we have $\mathfrak{H}_{\mathscr{R}}=Q^\vee$
and hence $\mathfrak{D}\setminus\mathfrak{H}_{\mathscr{R}}$ consists of only one connected 
component $\mathfrak{D}\setminus\mathfrak{H}_{\mathscr{R}}=\mathring{\mathfrak{D}}=\mathfrak{D}^{(1)}$.
Therefore we have
\begin{equation}
  F^{(1)}(\mathbf{t},\mathbf{y};A_1)=\frac{te^{ty}}{e^t-1}
\end{equation}
for $\mathbf{y}\in\mathfrak{D}^{(1)}$,
which coincides with the generating function of the classical Bernoulli polynomials.
Since $\mathfrak{J}=\{1\}$,
we see that
$\mathbf{B}^{(1)}_{\mathbf{k}}(\mathbf{y};A_1)$ consists of only one component, that is,
the classical Bernoulli polynomial $B_k(y)$.
The extended Weyl group $\Aut$ is $\{\id,\sigma_1\}$. 
For $\mathbf{y}\in\mathfrak{D}^{(1)}$, we have $\sigma_1\mathbf{y}=-\mathbf{y}$ and hence
$\tau(\alpha_1^\vee)\sigma_1\mathbf{y}=\alpha_1^\vee-\mathbf{y}\in\mathfrak{D}^{(1)}$ due to $0<\langle\alpha_1^\vee-\mathbf{y},\lambda_1\rangle<1$,
which implies that
the action on $\mathbb{Q}\,[\mathbf{y}]$ is given by $\bigl(\phi(\sigma_1)f\bigr)(\mathbf{y})=f(\alpha_1^\vee-\mathbf{y})$.
Therefore we have the well-known property
\begin{equation}
  (\phi(\sigma_1)\mathbf{B}^{(1)}_{\mathbf{k}})(\mathbf{y};A_1)=
  B_k(1-y)
=(-1)^kB_k(y)
=(-1)^k\mathbf{B}^{(1)}_{\mathbf{k}}(\mathbf{y};A_1).
\end{equation}
\end{example}
\begin{example}
\label{eg:exam_A2}
In the root system of type $A_2$, 
we have $\Delta_+=\{\alpha_1, \alpha_2, \alpha_1+\alpha_2\}$
and $\fs=\{\alpha_1,\alpha_2\}$.
Then $\mathscr{R}=\Delta_+$ and from Lemma \ref{lm:loc_fin} we have
\begin{equation}
\mathfrak{H}_{\mathscr{R}}=
(\mathbb{R}\,\alpha_1^\vee+\mathbb{Z}\,\alpha_2^\vee)
\cup
(\mathbb{R}\,\alpha_2^\vee+\mathbb{Z}\,\alpha_1^\vee)
\cup
(\mathbb{R}\,(\alpha_1^\vee+\alpha_2^\vee)+\mathbb{Z}\,\alpha_1^\vee).
\end{equation}
We see that $\mathfrak{D}\setminus\mathfrak{H}_{\mathscr{R}}=\mathfrak{D}^{(1)}\coprod\mathfrak{D}^{(2)}$, where
\begin{gather}
\mathfrak{D}^{(1)}=\{\mathbf{y}~|~0<y_2<y_1<1\},\\
\mathfrak{D}^{(2)}=\{\mathbf{y}~|~0<y_1<y_2<1\},
\end{gather}
with $y_1=\langle \mathbf{y},\lambda_1\rangle$ and $y_2=\langle \mathbf{y},\lambda_2\rangle$.
Let $t_1=t_{\alpha_1}$, $t_2=t_{\alpha_2}$ and $t_3=t_{\alpha_1+\alpha_2}$. 
For $\mathbf{y}\in\mathfrak{D}^{(1)}$,
the vertices of the polytopes $\mathcal{P}_{m_1,m_2,\mathbf{y}}$ in
\eqref{eq:D_m_y} are given by
\begin{equation}
  \begin{split}
    0,y_2 &\qquad\text{for }\mathcal{P}_{0,0,\mathbf{y}},\\
    y_2,y_1 &\qquad\text{for }\mathcal{P}_{0,1,\mathbf{y}},\\
    y_1,1 &\qquad\text{for }\mathcal{P}_{1,1,\mathbf{y}}.
  \end{split}
\end{equation}
Then by Theorem \ref{thm:gen_func} and Lemma \ref{lm:simplex}, 
we have
\begin{multline}
\label{eq:F_A_2}
  F^{(1)}(\mathbf{t},\mathbf{y};A_2)
\\
\begin{aligned}
  &=
  \frac{t_1t_2t_3e^{t_1y_1+t_2y_2}}{(e^{t_1}-1)(e^{t_2}-1)(e^{t_3}-1)}
\\
&\qquad
\times
  \biggl(
y_2\frac{e^{(t_3-t_1-t_2)y_2}-1}{(t_3-t_1-t_2)y_2}
+
  e^{t_2}(y_1-y_2)\frac{(e^{(t_3-t_1-t_2)y_1}-e^{(t_3-t_1-t_2)y_2})}{(t_3-t_1-t_2)(y_1-y_2)}
\\&\qquad\qquad
+e^{t_1+t_2}(1-y_1)\frac{(e^{t_3-t_1-t_2}-e^{(t_3-t_1-t_2)y_1})}{(t_3-t_1-t_2)(1-y_1)}
  \biggr)
  \\
  &=
  \frac{t_1t_2t_3e^{t_1y_1+t_2y_2}}{(e^{t_1}-1)(e^{t_2}-1)(e^{t_3}-1)(t_3-t_1-t_2)}
  \bigl(
  e^{(t_3-t_1-t_2)y_2}-1
\\
&\qquad+
  e^{t_2}(e^{(t_3-t_1-t_2)y_1}-e^{(t_3-t_1-t_2)y_2})
  +e^{t_1+t_2}(e^{t_3-t_1-t_2}-e^{(t_3-t_1-t_2)y_1})  
  \bigr).
\end{aligned}
\end{multline}
Hence by the Taylor expansion of \eqref{eq:F_A_2}, we have
\begin{equation}
\label{eq:B_1_2_2_2}
  \begin{split}
  B^{(1)}_{2,2,2}(\mathbf{y};A_2)=
&\frac{1}{30240}
+\frac{1}{360} (y_1 y_2 
- y_1^2
- y_2^2)
+\frac{1}{144} (3 y_1 y_2^2 
- 3 y_1^2 y_2 
+ 2 y_1^3)
\\
&
+\frac{1}{72} (-2 y_1 y_2^3 
- 3 y_1^2 y_2^2 
+ 4 y_1^3 y_2 
- 2 y_1^4
+ y_2^4)
\\
&
+\frac{1}{240} (-5 y_1 y_2^4
+ 10 y_1^2 y_2^3 
+ 10 y_1^3 y_2^2 
- 15 y_1^4 y_2 
+ 6 y_1^5)
\\
&
+\frac{1}{240}( 6 y_1 y_2^5 
- 5 y_1^2 y_2^4 
- 5 y_1^4 y_2^2 
+ 6 y_1^5 y_2 
- 2 y_1^6
- 2 y_2^6).
\end{split}
\end{equation}
Similarly we can calculate $B^{(2)}_{2,2,2}(\mathbf{y};A_2)$ for $\mathbf{y}\in\mathfrak{D}^{(2)}$, which coincides with
\eqref{eq:B_1_2_2_2} with $y_1$ and $y_2$ exchanged.
In the case $\mathbf{y}=0$, from Lemma \ref{lm:cont_PF}, we have
\begin{equation}
  \begin{split}
  F(\mathbf{t};A_2)
  &=
  \frac{t_1t_2t_3e^{t_1+t_2}(e^{t_3-t_1-t_2}-1)}{(e^{t_1}-1)(e^{t_2}-1)(e^{t_3}-1)(t_3-t_1-t_2)}
  \\
  &=
1
+\frac{1}{12}(
t_1 t_2
-t_1 t_3
-t_2 t_3)
+\frac{1}{360}(
t_1 t_2 t_3^2
-t_1^2 t_2 t_3
-t_1 t_2^2 t_3 
) \\
&\qquad
+\frac{1}{720}(t_1^2 t_2^2 
+t_1^2 t_3^2 
+t_2^2 t_3^2)
+\frac{1}{30240}t_1^2 t_2^2 t_3^2+\cdots,
\end{split}
\end{equation}
by letting $\mathbf{y}\to 0$ in \eqref{eq:F_A_2}.
Therefore by Theorem \ref{thm:W-Z}, we recover Mordell's formula \cite{Mo}:
\begin{equation}
  \label{WVF-A2}
  \zeta_2(2,2,2;A_2)=(-1)^3\frac{(2\pi\sqrt{-1})^6}{3!}\frac{1}{30240}=\frac{\pi^6}{2835}.
\end{equation}

We also discuss the action of the extended Weyl group.
Note that $\Aut$ is generated by $\{\sigma_1,\sigma_2,\omega\}$ 
where $\omega$ is a unique element of $\Omega$ such that $\omega\neq\id$,
and hence $\omega\alpha_1=\alpha_2$, $\omega\lambda_1=\lambda_2$ and $\omega^2=\id$.
Also note that $\alpha_1=2\lambda_1-\lambda_2$ and $\alpha_2=2\lambda_2-\lambda_1$.
For $\mathbf{y}\in\mathfrak{D}^{(1)}$, we have
\begin{equation}
  \begin{aligned}
    \langle\alpha_1^\vee+\sigma_1\mathbf{y},\lambda_1\rangle
    &=    
    1+\langle\mathbf{y},\sigma_1\lambda_1\rangle
    =
    1+\langle\mathbf{y},\lambda_1-\alpha_1\rangle=1-y_1+y_2,
    \\
    \langle\alpha_1^\vee+\sigma_1\mathbf{y},\lambda_2\rangle
    &=    
    \langle\mathbf{y},\sigma_1\lambda_2\rangle
    =
    \langle\mathbf{y},\lambda_2\rangle=y_2,
  \end{aligned}
\end{equation}
which implies
\begin{equation}
0<\langle\alpha_1^\vee+\sigma_1\mathbf{y},\lambda_2\rangle
<\langle\alpha_1^\vee+\sigma_1\mathbf{y},\lambda_1\rangle
<1.
\end{equation}
Therefore we have $\tau(\alpha_1^\vee)\sigma_1\mathfrak{D}^{(1)}=\mathfrak{D}^{(1)}$ and in a similar way $\sigma_1\mathfrak{D}^{(2)}=\mathfrak{D}^{(2)}$, and so on.
Thus
from \eqref{eq:act_phi}
 we see that
for $f=(f_1,f_2)$ with $f_1,f_2\in\mathbb{Q}\,[\mathbf{y}]$, the action of $\Aut$ is given by
\begin{align}
  (\phi(\sigma_1) f)(\mathbf{y})&=\bigl(f_1(\alpha_1^\vee+\sigma_1\mathbf{y}),f_2(\sigma_1\mathbf{y})\bigr),\notag\\
  (\phi(\sigma_2) f)(\mathbf{y})&=\bigl(f_1(\sigma_2\mathbf{y}),f_2(\alpha_2^\vee+\sigma_2\mathbf{y})\bigr),\\
  (\phi(\omega) f)(\mathbf{y})&=\bigl(f_2(\omega\mathbf{y}),f_1(\omega\mathbf{y})\bigr),\notag
\end{align}
or in terms of coordinates,
\begin{align}
  (\phi(\sigma_1) f)(y_1,y_2)&=\bigl(f_1(1-y_1+y_2,y_2),f_2(-y_1+y_2,y_2)\bigr),\notag\\
  (\phi(\sigma_2) f)(y_1,y_2)&=\bigl(f_1(y_1,y_1-y_2),f_2(y_1,1+y_1-y_2)\bigr),\\
  (\phi(\omega) f)(y_1,y_2)&=\bigl(f_2(y_2,y_1),f_1(y_2,y_1)\bigr).\notag
\end{align}
Then for
\begin{align}
\mathbf{B}^{(1)}_{2,2,2}(\mathbf{y};A_2)&=(B^{(1)}_{2,2,2}(\mathbf{y};A_2),0),\\
\mathbf{B}^{(2)}_{2,2,2}(\mathbf{y};A_2)&=(0,B^{(2)}_{2,2,2}(\mathbf{y};A_2)),
\end{align}
we have
\begin{align}
  \phi(\sigma_1)\mathbf{B}^{(1)}_{2,2,2}&=\mathbf{B}^{(1)}_{2,2,2},\notag\\
  \phi(\sigma_2)\mathbf{B}^{(1)}_{2,2,2}&=\mathbf{B}^{(1)}_{2,2,2},\\
  \phi(\omega)\mathbf{B}^{(1)}_{2,2,2}&=\mathbf{B}^{(2)}_{2,2,2}.\notag
\end{align}
More generally we can show 
\begin{align}
  \phi(\sigma_1)\mathbf{B}^{(1)}_{k_1,k_2,k_3}&=\mathbf{B}^{(1)}_{k_1,k_3,k_2},\notag\\
  \phi(\sigma_2)\mathbf{B}^{(1)}_{k_1,k_2,k_3}&=\mathbf{B}^{(1)}_{k_3,k_2,k_1},\\
  \phi(\omega)\mathbf{B}^{(1)}_{k_1,k_2,k_3}&=\mathbf{B}^{(2)}_{k_2,k_1,k_3}.\notag
\end{align}

\end{example}
\begin{example}
\label{eg:exam_B2}
The set of positive roots of type $C_2(\simeq B_2)$ consists of $\alpha_1$, $\alpha_2$,
 $2\alpha_1+\alpha_2$ and $\alpha_1+\alpha_2$.
Let $t_1=t_{\alpha_1}$, $t_2=t_{\alpha_2}$, $t_3=t_{2\alpha_1+\alpha_2}$ and $t_4=t_{\alpha_1+\alpha_2}$.
The vertices $(x_{2\alpha_1+\alpha_2},x_{\alpha_1+\alpha_2})$ of the polytopes $\mathcal{P}_{m_1,m_2}$ in \eqref{eq:D_m} are given by
\begin{equation}
  \begin{split}
    (0,0),(0,1/2),(1,0) &\qquad\text{for }\mathcal{P}_{1,1},\\
    (0,1),(0,1/2),(1,0) &\qquad\text{for }\mathcal{P}_{1,2},\\
    (0,1),(1,1/2),(1,0) &\qquad\text{for }\mathcal{P}_{2,2},\\
    (0,1),(1,1/2),(1,1) &\qquad\text{for }\mathcal{P}_{2,3}.
  \end{split}
\end{equation}
Then by Lemmas \ref{lm:simplex} and \ref{lm:simplex_expand}, we obtain 
\begin{equation}
  \begin{split}
    F(\mathbf{t};C_2)
    &=\biggl(\prod_{j=1}^4\frac{t_j}{e^{t_j}-1}\biggr)G(\mathbf{t};C_2)\\
    &=
1
+\frac{1}{2880} (
2 t_1 t_2 t_3^2 
-4 t_1 t_2 t_4^2 
-2 t_1 t_2^2 t_3 
+4 t_1 t_2^2 t_4 
-4 t_1 t_3 t_4^2 
-4 t_1 t_3^2 t_4 
\\
&\qquad\qquad\qquad\qquad
+t_1^2 t_2 t_3 
-4 t_1^2 t_2 t_4 
-4 t_1^2 t_3 t_4 
-t_2 t_3 t_4^2
-2 t_2 t_3^2 t_4
-2 t_2^2 t_3 t_4
)
\\
&\qquad
+\frac{1}{241920} (
3 t_1 t_2 t_3^2 t_4^2 
-3 t_1 t_2^2 t_3 t_4^2 
-3 t_1^2 t_2 t_3^2 t_4 
-3 t_1^2 t_2^2 t_3 t_4 
+2 t_1^2 t_2^2 t_3^2 
+8 t_1^2 t_2^2 t_4^2 
\\
&\qquad\qquad\qquad\qquad
+8 t_1^2 t_3^2 t_4^2
+2 t_2^2 t_3^2 t_4^2 
)
\\
&\qquad
+\frac{1}{9676800}t_1^2 t_2^2 t_3^2 t_4^2+\cdots,
\end{split}
\end{equation}
where
\begin{multline}
G(\mathbf{t};C_2)=\\
\begin{aligned}
&\frac{e^{t_3}}{(t_1+t_2-t_3) (t_1-2 t_3+t_4)}
-\frac{2 e^{\frac{t_1}{2}+\frac{t_4}{2}}}{(t_1+2 t_2-t_4) (t_1-2 t_3+t_4)}
+\frac{e^{t_1+t_2}}{(t_1+t_2-t_3) (t_1+2 t_2-t_4)}
\\
-&\frac{e^{t_2+t_3}}{(t_2+t_3-t_4) (t_1-2 t_3+t_4)}
+\frac{e^{t_4}}{(t_1+2 t_2-t_4) (t_2+t_3-t_4)}
+\frac{2 e^{\frac{t_1}{2}+t_2+\frac{t_4}{2}}}{(t_1+2 t_2-t_4) (t_1-2 t_3+t_4)}
\\
+&\frac{2 e^{\frac{t_1}{2}+t_3+\frac{t_4}{2}}}{(t_1+2 t_2-t_4) (t_1-2 t_3+t_4)}
+\frac{e^{t_1+t_2+t_3}}{(t_2+t_3-t_4) (t_1+2 t_2-t_4)}
-\frac{e^{t_1+t_4}}{(t_2+t_3-t_4) (t_1-2 t_3+t_4)}
\\
+&\frac{e^{t_3+t_4}}{(t_1+t_2-t_3) (t_1+2 t_2-t_4)}
-\frac{2 e^{\frac{t_1}{2}+t_2+t_3+\frac{t_4}{2}}}{(t_1+2 t_2-t_4) (t_1-2 t_3+t_4)}
+\frac{e^{t_1+t_2+t_4}}{(t_1+t_2-t_3) (t_1-2 t_3+t_4)}.
\end{aligned}
\end{multline}
Therefore
\begin{equation}
  \label{WVF-B2}
  \zeta_2(2,2,2,2;C_2)=(-1)^4\frac{(2\pi\sqrt{-1})^8}{2!\cdot 2^2}\frac{1}{9676800}=\frac{\pi^8}{302400},
\end{equation}
and so
\begin{equation}
  \zeta_W(2;C_2)=6^2\frac{\pi^8}{302400}=\frac{\pi^8}{8400}.
\end{equation}  
\end{example}

\begin{example}
\label{eg:exam_A3}
The set of positive roots of type $A_3$ consists of $\alpha_1$, $\alpha_2$, $\alpha_3$,
 $\alpha_1+\alpha_2$, $\alpha_2+\alpha_3$ and $\alpha_1+\alpha_2+\alpha_3$.
Let $t_1=t_{\alpha_1}$, $t_2=t_{\alpha_2}$, 
$t_3=t_{\alpha_3}$, 
$t_4=t_{\alpha_1+\alpha_2}$,
$t_5=t_{\alpha_2+\alpha_3}$ and
$t_6=t_{\alpha_1+\alpha_2+\alpha_3}$.
The vertices $(x_{\alpha_1+\alpha_2},x_{\alpha_2+\alpha_3},x_{\alpha_1+\alpha_2+\alpha_3})$
of the polytopes $\mathcal{P}_{m_1,m_2,m_3}$ in \eqref{eq:D_m} are given by
\begin{equation}
  \begin{split}
(0, 0, 0), (0, 0, 1), (0, 1, 0), (1, 0, 0)&\qquad\text{for }\mathcal{P}_{1, 1, 1},\\
(0, 0, 1), (0, 1, 0), (1, 0, 0), (1, 1, 0)&\qquad\text{for }\mathcal{P}_{1, 2, 1},\\
(0, 0, 1), (0, 1, 0), (0, 1, 1), (1, 1, 0)&\qquad\text{for }\mathcal{P}_{1, 2, 2},\\
(0, 0, 1), (1, 0, 0), (1, 0, 1), (1, 1, 0)&\qquad\text{for }\mathcal{P}_{2, 2, 1},\\
(0, 0, 1), (0, 1, 1), (1, 0, 1), (1, 1, 0)&\qquad\text{for }\mathcal{P}_{2, 2, 2},\\
(0, 1, 1), (1, 0, 1), (1, 1, 0), (1, 1, 1)&\qquad\text{for }\mathcal{P}_{2, 3, 2}.
\end{split}
\end{equation}
Then by Lemmas \ref{lm:simplex} and \ref{lm:simplex_expand}, we obtain 
\begin{equation}
  \begin{split}
    F(\mathbf{t};A_3)
    &=\biggl(\prod_{j=1}^6\frac{t_j}{e^{t_j}-1}\biggr)G(\mathbf{t};A_3)\\
    &=1+\cdots+\frac{23}{435891456000}t_1^2t_2^2t_3^2t_4^2t_5^2t_6^2+\cdots,
  \end{split}
\end{equation}
where
\begin{multline}
G(\mathbf{t};A_3)=\\
\begin{aligned}
&\frac{e^{t_3+t_4}}{(t_1+t_2-t_4) (t_1-t_3-t_4+t_5)(t_3+t_4-t_6)}
-\frac{e^{t_1+t_5}}{(t_2+t_3-t_5) (t_1-t_3-t_4+t_5)(t_1+t_5-t_6)}
\\
&
-\frac{e^{t_6}}{(t_1+t_2+t_3-t_6) (t_3+t_4-t_6)(t_1+t_5-t_6)}
+\frac{e^{t_1+t_2+t_3}}{(t_1+t_2-t_4) (t_2+t_3-t_5)(t_1+t_2+t_3-t_6)}
\\
&
-\frac{e^{t_4+t_5}}{(t_1+t_2-t_4) (t_2+t_3-t_5)(t_2-t_4-t_5+t_6)}
-\frac{e^{t_2+t_3+t_4}}{(t_2+t_3-t_5) (t_1-t_3-t_4+t_5)(t_3+t_4-t_6)}
\\
&
+\frac{e^{t_1+t_2+t_5}}{(t_1+t_2-t_4) (t_1-t_3-t_4+t_5)(t_1+t_5-t_6)}
+\frac{e^{t_2+t_6}}{(t_3+t_4-t_6) (t_1+t_5-t_6)(t_2-t_4-t_5+t_6)}
\\
&
+\frac{e^{t_3+t_4+t_5}}{(t_1+t_2-t_4) (t_3+t_4-t_6)(t_2-t_4-t_5+t_6)}
-\frac{e^{t_5+t_6}}{(t_2+t_3-t_5) (t_1+t_2+t_3-t_6)(t_3+t_4-t_6)}
\\
&
+\frac{e^{t_1+t_2+t_3+t_5}}{(t_1+t_2-t_4) (t_1+t_2+t_3-t_6)(t_1+t_5-t_6)}
-\frac{e^{t_2+t_3+t_6}}{(t_2+t_3-t_5) (t_1+t_5-t_6)(t_2-t_4-t_5+t_6)}
\\
&
+\frac{e^{t_1+t_4+t_5}}{(t_2+t_3-t_5) (t_1+t_5-t_6)(t_2-t_4-t_5+t_6)}
-\frac{e^{t_4+t_6}}{(t_1+t_2-t_4) (t_1+t_2+t_3-t_6)(t_1+t_5-t_6)}
\\
&
+\frac{e^{t_1+t_2+t_3+t_4}}{(t_2+t_3-t_5) (t_1+t_2+t_3-t_6)(t_3+t_4-t_6)}
-\frac{e^{t_1+t_2+t_6}}{(t_1+t_2-t_4) (t_3+t_4-t_6)(t_2-t_4-t_5+t_6)}
\\
&
-\frac{e^{t_1+t_3+t_4+t_5}}{(t_1+t_5-t_6) (t_2-t_4-t_5+t_6)(t_3+t_4-t_6)}
-\frac{e^{t_3+t_4+t_6}}{(t_1+t_2-t_4) (t_1-t_3-t_4+t_5)(t_1+t_5-t_6)}
\\
&
+\frac{e^{t_1+t_5+t_6}}{(t_2+t_3-t_5) (t_1-t_3-t_4+t_5)(t_3+t_4-t_6)}
+\frac{e^{t_1+t_2+t_3+t_6}}{(t_1+t_2-t_4) (t_2+t_3-t_5)(t_2-t_4-t_5+t_6)}
\\
&
-\frac{e^{t_4+t_5+t_6}}{(t_1+t_2-t_4) (t_2+t_3-t_5)(t_1+t_2+t_3-t_6)}
+\frac{e^{t_1+t_2+t_3+t_4+t_5}}{(t_1+t_2+t_3-t_6) (t_3+t_4-t_6)(t_1+t_5-t_6)}
\\
&
+\frac{e^{t_2+t_3+t_4+t_6}}{(t_2+t_3-t_5) (t_1-t_3-t_4+t_5)(t_1+t_5-t_6)}
-\frac{e^{t_1+t_2+t_5+t_6}}{(t_1+t_2-t_4) (t_1-t_3-t_4+t_5)(t_3+t_4-t_6)}.
\end{aligned}
\end{multline}
Therefore we obtain
\begin{equation}
  \label{WVF-A3}
  \zeta_3(2,2,2,2,2,2;A_3)=(-1)^6\frac{(2\pi\sqrt{-1})^{12}}{4!}\frac{23}{435891456000}=\frac{23}{2554051500}\pi^{12},
\end{equation}
which implies a formula of Gunnells-Sczech \cite{GS}:
\begin{equation}
  \zeta_W(2;A_3)=12^2\frac{23}{2554051500}\pi^{12}=\frac{92}{70945875}\pi^{12}.
\end{equation}  
\end{example}

In higher rank root systems, generating functions are more involved, 
since the polytopes are not simplicial any longer.
For instance,
we have
the generating function of type $G_2$, $A_4$, $B_3$ and $C_3$ with
1010 terms, 5040 terms, 19908 terms and 20916 terms respectively by
use of triangulation.
In \cite{KM5}, we will
improve Theorem
\ref{thm:gen_func}
and will give more compact forms of the generating functions $F(\mathbf{t},\mathbf{y};\Delta)$,
which do not depend on simplicial decompositions.
As a result, the numbers of terms in the above generating functions 
reduce to $15$, $125$, $68$ and $68$ respectively
(as for the $G_2$ case, see \cite{KM4}).

\begin{example}
\label{eg:funcrel}
In Theorem \ref{thm:FR}, we have already given general forms of
functional relations among zeta-functions of root systems. 
In previous
examples we observed generating functions and special values in
several cases, but here, we treat examples of explicit functional
relations which can be deduced from the general forms.   First
consider the $A_2$  case (see Example \ref{eg:exam_A2}).
Set
$$\Delta_{+}=\Delta_{+}(A_2)=\{ \alpha_1,\alpha_2,\alpha_1+\alpha_2\},$$
and ${\bf y}=0$, ${\bf s}=(2,s,2)$ for $s \in \mathbb{C}$ with $\Re
s>1$, $I=\{ 2\}$, that is, $\Delta_{I+}=\{ \alpha_2\}$. Then, from
\eqref{eq:def_S}, we can write the left-hand side of
\eqref{eq:func_eq} as
\begin{align*}
S({\bf s},{\bf y};I;\Delta) & =\sum_{m,n=1}^\infty \frac{1}{m^2 n^s (m+n)^2}+\sum_{m,n=1 \atop m\not=n}^\infty \frac{1}{m^2 n^s (-m+n)^2} \\
& \ =2\zeta_2(2,s,2;A_2)+\zeta_2(2,2,s;A_2). \notag
\end{align*}
On the other hand, the right-hand side of \eqref{eq:func_eq} is 
\begin{align*}
& \left( \frac{(2\pi\sqrt{-1})^2}{2!}\right)^2\sum_{m=1}^\infty \frac{1}{m^s} \int_{0}^{1}e^{-2\pi \sqrt{-1}mx}L_{2}(x,0)L_{2}(-x,0)dx \\
& =\left( \frac{(2\pi\sqrt{-1})^2}{2!}\right)^2\sum_{m=1}^\infty \frac{1}{m^s} \int_{0}^{1}e^{-2\pi \sqrt{-1}mx}B_{2}(x)B_{2}(1-x)dx, 
\end{align*}
by using \eqref{eq:L_is_B-2}.  From well-known properties of Bernoulli
polynomials, we can calculate the above integral (for details, see
Nakamura \cite{Na1}) and can recover from \eqref{eq:func_eq} the
formula
\begin{equation}
2\zeta_2(2,s,2;A_2)+\zeta_2(2,2,s;A_2)=4\zeta(2)\zeta(s+2)-6\zeta(s+4), \label{Fn-Mordell}
\end{equation}
proved in \cite{TsC} (see also \cite{Na1}).  The function
$\zeta_2(\mathbf{s};A_2)$ can be continued meromorphically to the
whole space $\mathbb{C}^3$ (\!\cite{Ma0}), so \eqref{Fn-Mordell} holds
for any $s \in \mathbb{C}$ except for singularities on the both sides.
In particular when $s=2$, we obtain \eqref{WVF-A2}. Similarly we can
treat the $C_2(\simeq B_2)$ case and give some functional relations from
\eqref{eq:func_eq} by combining the meromorphic continuation of
$\zeta_2(\mathbf{s};C_2)$ which has been shown in \cite{Ma2}.
\end{example}

\end{document}